\documentclass{report}
\usepackage{amsmath,amssymb,amsthm}
\usepackage{graphicx}
\usepackage[all]{xy}

\newtheorem{theorem}{Theorem}[section]
\newtheorem{proposition}[theorem]{Proposition}
\newtheorem{corollary}[theorem]{Corollary}
\newtheorem{lemma}[theorem]{Lemma}

\newtheorem{convention}[theorem]{Convention}
\newtheorem{observation}[theorem]{Observation}
\newtheorem{question}[theorem]{Question}
\newtheorem{claim}[theorem]{Claim}
\theoremstyle{definition}

\newtheorem{notes}[theorem]{Notes}
\newtheorem*{note*}{Note}
\newtheorem*{corollary*}{Corollary}
\newtheorem{definition}[theorem]{Definition}
\newtheorem{example}[theorem]{Example}
\newtheorem{remark}[theorem]{Remark}
\newtheorem{terminology}[theorem]{Terminology}
\newtheorem{notation}[theorem]{Notation}

\newcommand{\pair}[2]{{$\left[\scriptstyle\begin{matrix}#1\\#2\end{matrix}\right]$}}
\DeclareMathOperator{\PL}{\mathrm{PL}}
\DeclareMathOperator{\LL}{\mathbf{L}}
\DeclareMathOperator{\NN}{\mathbf{N}}
\DeclareMathOperator{\RR}{\mathbf{R}}
\DeclareMathOperator{\II}{\mathbf{I}}
\DeclareMathOperator{\EG}{\mathrm{EG}}
\DeclareMathOperator{\AG}{\mathrm{AG}}
\DeclareMathOperator{\NF}{\mathrm{NF}}
\DeclareMathOperator{\SG}{\mathrm{SG}}

\title{Thompson's Group $F$}
\author{James Michael Belk}
\date{Ph.D. Thesis \\ Cornell University \\ August 2004}

\begin{document}

\maketitle
\abstract{We introduce two new types of diagrams that aid in understanding elements of
Thompson's group $F$.

The first is the \emph{two-way forest diagram}, which represents an element of $F$
as a pair of infinite, bounded binary forests together with an order-preserving bijection
of the leaves.  These diagrams have the same relationship to a certain action of $F$
on the real line that the standard tree diagrams have to the action of $F$ on the unit
interval.  Using two-way forest diagrams, we derive a simple length formula for elements
of $F$ with respect to the finite generating set $\{x_0,x_1\}$.

We then discuss several applications of two-way forest diagrams and the length formula
to the geometry of $F$.  These
include a simplification of a result by S.~Cleary and J.~Taback that $F$ has dead ends
but no deep pockets;  a precise calculation of the growth function of the positive submonoid
with respect to the $\{x_0,x_1\}$ generating set;  a new upper bound on the isoperimetric
constant (a.k.a. Cheeger constant) of $F$; and a proof that $F$ is not minimally almost convex.

Next, we introduce \emph{strand diagrams} for elements of $F$.  These are similar to tree
diagrams, but they can be concatenated like braids.  Motivated by the fact that configuration
spaces are classifying spaces for braid groups, we present a classifying space for $F$ that is
the ``configuration space'' of finitely many points on a line, with the points allowed to split
and merge in pairs.

In addition to the new results, we present a thorough exposition of the basic theory of the
group $F$.  Highlights include a simplified proof that the commutator subgroup of $F$ is
simple, a discussion of open problems (with a focus on amenability), and a simplified
derivation of the standard presentation for $F$ and the normal form for elements
using \emph{one-way forest diagrams}.}

\newpage

\tableofcontents

\newpage

\chapter{The Group $F$}

Thompson's group $F$ is a certain group of piecewise-linear homeomorphisms of $[0,1]$.  We define
$F$ in section 1, and prove an important characterization of its elements.  In section 2 we
introduce \emph{tree diagrams}, an important tool for understanding elements of $F$.
We go on to prove some of the basic properties of $F$ in sections 3 and 4, and in section
5 we discuss various open problems, with a focus on amenability.
Finally, we discuss several other approaches toward $F$ in section 6, and show how they
relate to our homeomorphism approach.

None of the results in this chapter are new.  However, we have endeavored to make many of the
proofs simpler and clearer than in previously published versions.  Furthermore, this chapter
does not contain proofs for two primary results: the standard presentation for $F$
and the existence of normal forms.  Interested readers should look ahead to section 2.4
for a treatment using forest diagrams, or consult \cite{CFP} for the traditional tree-diagrams
approach.

\section{Dyadic Rearrangements}

Suppose we take the interval $\left[0,1\right]$, and cut it in half,
like this:
\begin{center}
\includegraphics{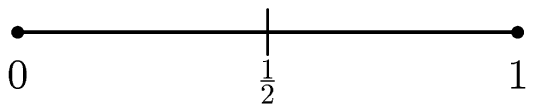}
\end{center}
We then cut each of the resulting intervals in half:
\begin{center}
\includegraphics{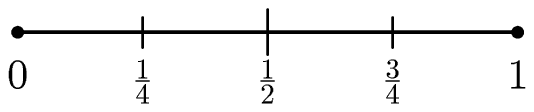}
\end{center}
and then cut some of the new intervals in half:
\begin{center}
\includegraphics{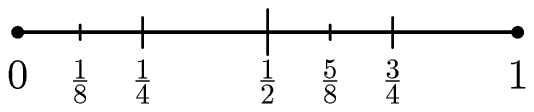}
\end{center}
to get a certain subdivision of $[0,1]$. Any subdivision of
$[0,1]$ obtained in this manner (i.e. by repeatedly cutting
intervals in half) is called a \emph{dyadic subdivision}.

The intervals of a dyadic subdivision are all of the form:
\begin{equation*}
\left[\frac{k}{2^n},\frac{k+1}{2^n}\right]\qquad k,n\in\mathbb{N}
\end{equation*}
These are the \emph{standard dyadic intervals}. We could alternatively define a
dyadic subdivision as any partition of $[0,1]$ into standard
dyadic intervals.

Given a pair $\mathcal{D},\mathcal{R}$ of dyadic subdivisions with the
same number of cuts, we can define a piecewise-linear homeomorphism
$f\colon[0,1]\rightarrow[0,1]$ by sending each interval
of $\mathcal{D}$ linearly onto the corresponding interval of $\mathcal{R}$.
This is called a \emph{dyadic rearrangement} of $[0,1]$

\begin{example}
Here are two dyadic rearrangements:
\begin{center}
\includegraphics{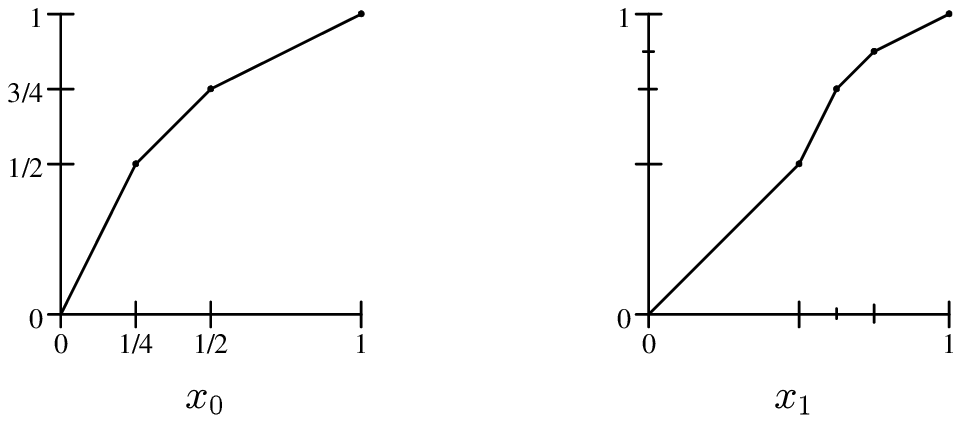}
\end{center}
\end{example}

\begin{theorem}
Let $f\colon\left[0,1\right]\rightarrow\left[0,1\right]$ be
a piecewise-linear homeomorphism. Then $f$ is a dyadic rearrangement if and
only if
\begin{enumerate}
\item All the slopes of $f$ are powers of $2$, and
\item All the breakpoints of $f$ have dyadic rational coordinates.
\end{enumerate}
\end{theorem}
\begin{proof}
Clearly every dyadic rearrangement satisfies conditions (1) and (2).
Suppose now that $f$ is a piecewise-linear homeomorphism satisfying these
two conditions. Choose $N$ sufficiently large so that:
\begin{enumerate}
\item $f$ is linear on each standard dyadic interval of width
$\displaystyle\frac{1}{2^N}$, and
\item The formula for each linear segment of $f$ can be written as:
\begin{equation*}
f(t)=2^m\left(t+\frac{k}{2^N}\right)\qquad m,k\in\mathbb{Z}
\end{equation*}
\end{enumerate}
Let $\mathcal{D}$ be the subdivision of $\left[0,1\right]$ into standard
dyadic intervals of width $1/2^N$. Then $f$ maps each interval of
$\mathcal{D}$ linearly to a standard dyadic interval, and therefore maps
$\mathcal{D}$ to some dyadic subdivision of $\left[0,1\right]$.
\end{proof}

\begin{corollary}
The set $F$ of all dyadic rearrangements forms a group under composition.
\end{corollary}

This group is called \emph{Thompson's Group $F$}.

\begin{theorem} $F$ is infinite and torsion-free. \end{theorem}
\begin{proof}  Let $f$ be any element of $F$ that is not the identity, and let:
\begin{equation*}
t_0=\inf\left\{t\in\left[0,1\right]:f(t)\neq t\right\}
\end{equation*}
Then $f(t_0)=t_0$, and $f$ has right-hand derivative $2^m$ at $t_0$
for some $m\neq 0$. By the chain rule, the right-hand derivative
of $f^n$ at $t_0$ is $2^{mn}$ for all $n\in\mathbb{N}$, hence all the
positive powers of $f$ are distinct.\end{proof}

The group $F$ was first defined by Richard J. Thompson in the 1960's, in connection
with his work on associativity.  It was later rediscovered by topologists
(Freyd and Heller, and independently Dydak)
who were researching the structure of topological spaces with homotopy idempotents.
(See section 1.6 for a discussion of these connections.)  Since then
$F$ has become an important object of study in geometric group theory, primarily
because of some long-standing problems regarding the geometric structure of its
Cayley graph (see section 1.5).

When Thompson originally defined $F$, he used ``backwards'' notation
for composition of functions. We will adopt this convention throughout:

\begin{convention} If $f$ and $g$ are functions, the expression $fg$ will denote ``$f$
followed by $g$''. In particular:
\begin{equation*}
(fg)(t)=g(f(t))
\end{equation*}
\end{convention}

\section{Tree Diagrams}

The standard dyadic intervals form a binary tree under inclusion:
\begin{center}
\includegraphics{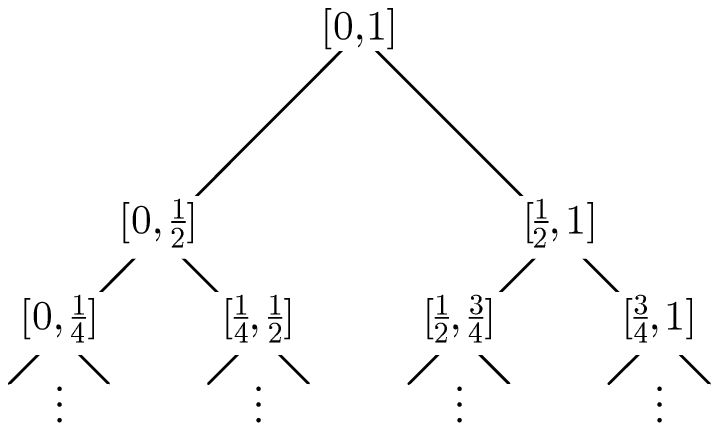}
\end{center}
Dyadic subdivisions of $\left[0,1\right]$ correspond to \emph{finite subtrees} of
this infinite binary tree. For example, the subdivision:
\begin{center}
\includegraphics{tpict03}
\end{center}
corresponds to the subtree:
\begin{center}
\includegraphics{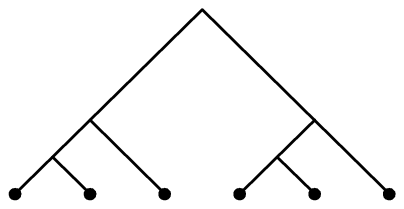}
\end{center}
Each leaf of this tree represents an interval of the subdivision, and the
root represents the interval $\left[0,1\right]$. The other nodes represent
intervals from intermediate stages of the dyadic subdivision.

Using this scheme, we can describe any element of $F$ by a pair of
finite binary trees. This is called a \emph{tree diagram}.

\begin{example}  Recall that the element $x_0$ sends intervals of the
subdivision:
\begin{center}
\includegraphics{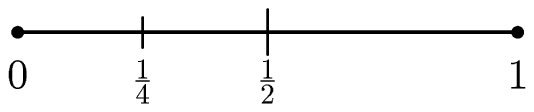}
\end{center}
linearly onto intervals of the subdivision:
\begin{center}
\includegraphics{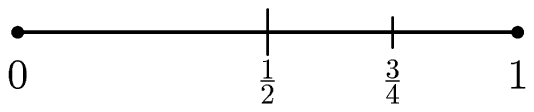}
\end{center}
Therefore, $x_0$ has tree diagram:
\begin{center}
\includegraphics{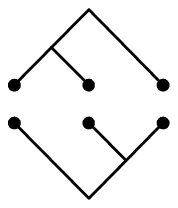}
\end{center}
We have aligned the two trees vertically so that corresponding leaves match
up. By convention, the \emph{domain} tree appears on the \emph{top}, and the
\emph{range} tree appears on the \emph{bottom}.
\end{example}

\begin{example} The tree diagram for $x_1$ is:
\begin{center}
\includegraphics{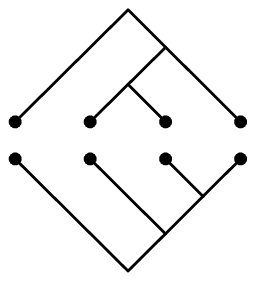}
\end{center}
\end{example}

Of course, the tree diagram for an element of $F$ is not unique. For
example, all of the following are tree diagrams for the identity:
\begin{center}
\includegraphics{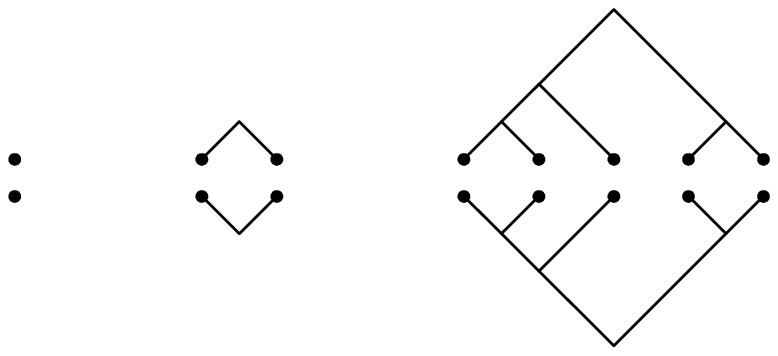}
\end{center}
In general, a \emph{reduction} of a tree diagram consists of removing an
opposing pair of carets, like this:
\begin{center}
\includegraphics{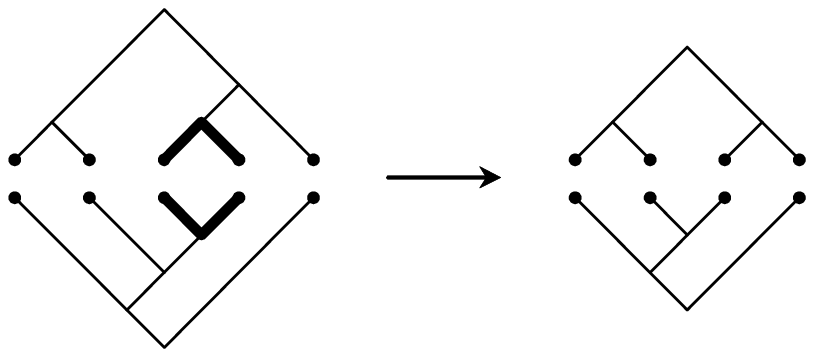}
\end{center}
Performing a reduction does not change the element of $F$ described by the
tree diagram --- it merely corresponds to removing an unnecessary ``cut''
from the subdivisions of the domain and range.

\begin{definition} A tree diagram is \emph{reduced} if it has no opposing pairs
of carets.
\end{definition}

\begin{theorem}Every element of $F$ has a unique reduced tree diagram.
\end{theorem}
\begin{proof} Note first that a tree diagram for a given $f\in F$ is determined
entirely by the domain tree. Furthermore, if $T \subset T'$ are possible
domain trees, then the tree diagram with domain tree $T$ is a reduction
of the tree diagram with domain tree $T'$. Therefore, it suffices to show that the set of
possible domain trees for $f$ has a minimum element under inclusion.

Define a standard dyadic interval to be \emph{regular} if $f$ maps it linearly
onto a standard dyadic interval. Then a tree $T$ is a possible domain tree for
$f$ if and only if its leaves are all regular. We conclude that the set of
possible domain trees is closed under intersections,
and therefore has a minimum element. \end{proof}

We will generally denote tree diagrams by column vectors \pair{T}{U},
where $T$ and $U$ are the component binary trees.

The following observation makes it possible to multiply two elements of
$F$ directly from the tree diagrams:

\begin{observation} Suppose that $f,g\in F$ have tree diagrams \pair{T}{U} and
\pair{U}{V}.  Then \pair{T}{V} is a tree diagram for $fg$.
\end{observation}

Therefore, to multiply two elements $f$ and $g$, we need only find a corresponding pair of tree diagrams
such that the bottom tree of $f$ is congruent to the top tree of $g$.  Such a pair can always be obtained
by expanding the reduced tree diagrams:

\begin{example} Let $f$ and $g$ be the elements:
\begin{center}
\includegraphics{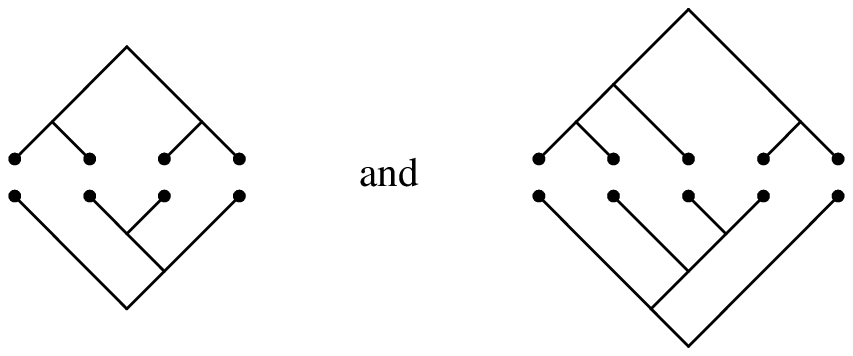}
\end{center}
We can expand the tree diagrams for $f$ and $g$ to get:
\begin{center}
\includegraphics{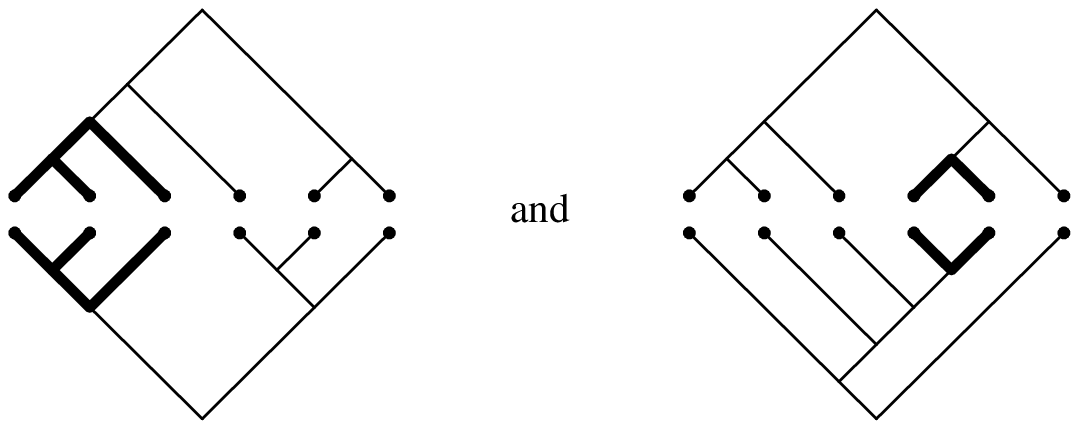}
\end{center}
Note that the bottom tree for $f$ is now the same as the top tree for $g$.
Therefore, $fg$ has tree diagram:
\begin{center}
\includegraphics{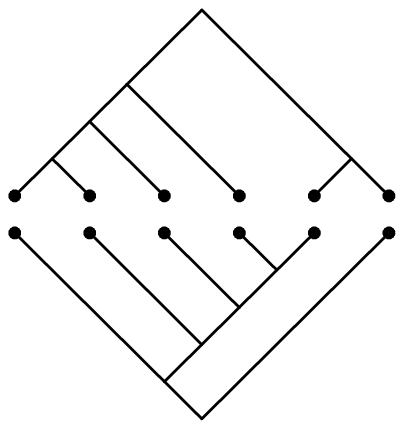}
\end{center}
\end{example}

It is also easy to find a tree diagram for the inverse of an element:

\begin{observation}
If \pair{T}{U} is a tree diagram for $f\in F$, then \pair{U}{T} is a tree diagram for $f^{-1}$.
\end{observation}

\section{Generators}

Let $x_0,x_1,x_2,\ldots$ be the elements of $F$ with tree diagrams:
\begin{center}
\includegraphics{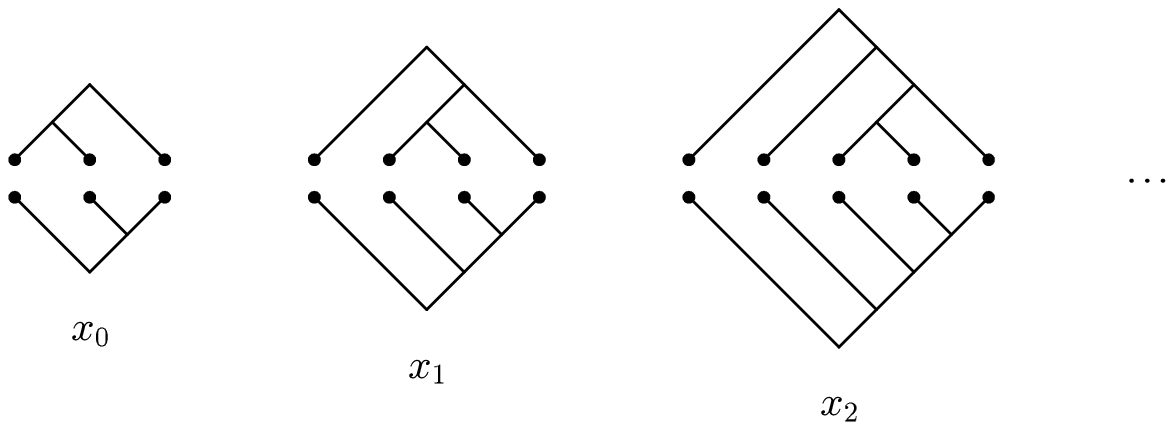}
\end{center}
In this section, we will show that these elements generate the group $F$.  We will also
state without proof an infinite presentation for $F$, and a normal form for words in the generators.
Finally, we will show that the elements $\{x_0,x_1\}$ alone generate $F$, and derive a finite presentation
for $F$ using these generators.

First, observe that the bottom trees of the elements $\{x_0,x_1,x_2,\ldots\}$ all have the same form:
a long edge on the right with left edges emanating from it.
Such a tree is called a \emph{right vine}:
\begin{center}
\includegraphics{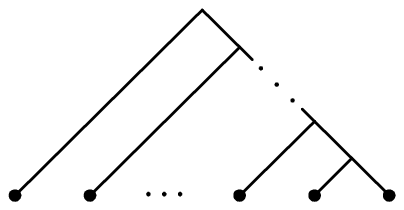}
\end{center}

\begin{definition} An element of $F$ is \emph{positive} if the bottom tree of
its reduced tree diagram is a right vine.
\end{definition}

Since any subtree of a right vine is a right vine, any tree diagram (reduced or not)
whose bottom tree is a right vine represents a positive element.

Given a binary tree $T$, let $[T]$ denote the positive element with top tree $T$.

\begin{proposition} Every element of $F$ can be expressed as $pq^{-1}$,
where $p$ and $q$ are positive. \end{proposition}
\begin{proof} If $f$ has tree diagram \pair{T}{U}, then $f=[T][U]^{-1}$.\end{proof}

Define the \emph{width} of a binary tree to be the number of leaves minus one.
Given a binary tree of width $w$, we number its leaves $0,1,\ldots ,w$ from left to right.

\begin{proposition} Let $T$ be a binary tree of width $w$. If $n<w$, then:
\begin{equation*}
[T]x_n=[T\wedge n]
\end{equation*}
where $T\wedge n$ is the binary tree obtained by attaching a caret to the
$n\!\text{'th}$ leaf of $T$.
\end{proposition}
\begin{proof}Let $V$ be a right vine of width $w$.  Then
\pair{T\wedge n}{V\wedge n} is a tree diagram for $[T]$,
and $x_n = [V \wedge n]$ (since $n<w$), so:
\begin{equation*}
[T]x_n=[T\wedge n]
\end{equation*}
by observation 1.2.5. \end{proof}

Define the \emph{right stalk} of a binary tree $T$ to be the right vine that grows
from the root of $T$:
\begin{center}
\includegraphics{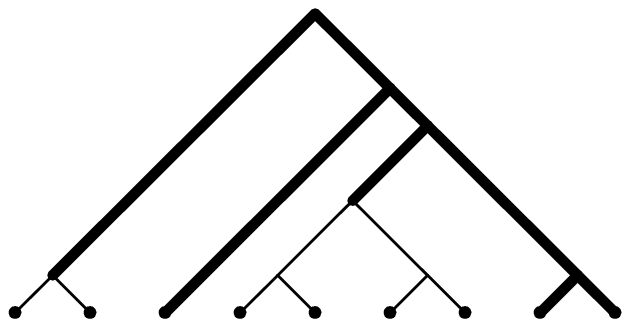}
\end{center}
Clearly any binary tree $T$ can be constructed by starting with its
right stalk and then attaching carets to the leaves one at a time.  Therefore:

\begin{corollary} The set of positive elements is precisely the submonoid
generated by $\left\{x_0,x_1,x_2,\ldots\right\}$.
\end{corollary}

\begin{corollary} The elements $\left\{x_0,x_1,x_2,\ldots\right\}$ generate Thompson's
    group $F$.
\end{corollary}

We will prove the following two theorems in section 2.4 using forest diagrams.  See
\cite{CFP} for a tree-diagram approach.

\begin{theorem}  $F$ has presentation:
\begin{equation*}
\langle x_0,x_1,x_2,\ldots\mid x_nx_k=x_kx_{n+1}\text{ for }k<n\rangle\tag*{\qedsymbol}
\end{equation*}
\end{theorem}

\begin{theorem}[Normal Form]Every element of $F$ can be expressed uniquely in the form:
\begin{equation*}
x_0^{a_0} \cdots x_n^{a_n} x_n^{-b_n} \cdots x_0^{-b_0}
\end{equation*}
where $a_0,\ldots,a_n,b_0,\ldots,b_n\in\mathbb{N}$, exactly one of $a_n,b_n$ is
nonzero, and:
\begin{equation*}
a_i\neq0\text{ and }b_i\neq0\quad\Rightarrow\quad a_{i+1}\neq0\text{ or }b_{i+1}\neq0
\end{equation*}
for all $i$.\quad\qedsymbol
\end{theorem}

It is possible to put any element of $F$ into normal form using the following four types
of moves, each of which follows from the relations given in theorem 1.3.6:
\begin{equation*}
\begin{aligned}
x_n^{-1}x_k&\rightarrow x_kx_{n+1}^{-1}\\
x_k^{-1}x_n&\rightarrow x_{n+1}x_k^{-1}\\
x_nx_k&\rightarrow x_kx_{n+1}\\
x_k^{-1}x_n^{-1}&\rightarrow x_{n+1}^{-1}x_k^{-1}
\end{aligned}
\qquad\text{(where }k<n\text{)}
\end{equation*}
These operations allow us to interchange any two generators that are in the
``wrong'' order, at the expense of incrementing the larger of the two subscripts.

It is also possible to use the inverses of these operations to put two
generators in the ``wrong'' order, but \emph{only if the subscripts differ by more
than one}. For example, we can switch each of the following generator pairs:
\begin{equation*}
x_3x_5\qquad x_3x_5^{-1}\qquad x_5x_3^{-1}\qquad x_5^{-1}x_3^{-1}
\end{equation*}
but we cannot apply an inverse operation to any of the following pairs:
\begin{equation*}
x_3x_4\qquad x_3x_4^{-1}\qquad x_4x_3^{-1}\qquad x_4^{-1}x_3^{-1}
\end{equation*}

\begin{example} Suppose we wish to put the word:
\begin{equation*}
x_0x_3x_6x_3^{-1}x_1x_4^{-1}x_0x_3^{-1}x_0^{-1}
\end{equation*}
into normal form. We start by applying operations of types (1) and (2).
(In each step, the generators about to be interchanged are indicated.)
\begin{equation*}
\begin{array}{l@{}l@{}l@{}l@{}l@{}l@{}l@{}l@{}l@{}l}
  & x_0\phantom{()} & x_3 \phantom{()}& x_6\phantom{)}(& x_3^{-1} & x_1)(    & x_4^{-1} & x_0)\phantom{(}     & x_3^{-1} & x_0^{-1}\\
=\phantom{()} & x_0 & x_3 & x_6 & x_1\phantom{)}(     & x_4^{-1} & x_0)     & x_5^{-1} & x_3^{-1} & x_0^{-1}\\
= & x_0 & x_3 & x_6 & x_1      & x_0      & x_5^{-1}\phantom{)} & x_5^{-1}\phantom{)} & x_3^{-1}\phantom{)} & x_0^{-1}
\end{array}
\end{equation*}
Next we apply operations of type (3) to arrange the left half of the word. (The right half is already arranged.)
\begin{equation*}
\begin{array}{l@{}l@{}l@{}l@{}l@{}l@{}l@{}l@{}l@{}l}
 &x_0&x_3&x_6\phantom{)}(&x_1&x_0)\phantom{(}&x_5^{-1}\phantom{)}&x_5^{-1}\phantom{)}&x_3^{-1}\phantom{)}&x_0^{-1}\phantom{)}\\
=\phantom{()}&x_0&x_3\phantom{)}(&x_6&x_0)\phantom{(}&x_2&x_5^{-1}&x_5^{-1}&x_3^{-1}&x_0^{-1}\\
=&x_0\phantom{)}(&x_3&x_0)(&x_7&x_2)&x_5^{-1}&x_5^{-1}&x_3^{-1}&x_0^{-1}\\
=&x_0&x_0\phantom{)}(&x_4&x_2)\phantom{(}&x_8&x_5^{-1}&x_5^{-1}&x_3^{-1}&x_0^{-1}\\
=&x_0&x_0&x_2&x_5&x_8&x_5^{-1}&x_5^{-1}&x_3^{-1}&x_0^{-1}
\end{array}
\end{equation*}
At this point, there are instances of $x_0$ and $x_0^{-1}$, but no instances
of $x_1$ or $x_1^{-1}$. We can therefore cancel an $x_0,x_0^{-1}$ pair:
\begin{equation*}
\begin{array}{l@{}l@{}l@{}l@{}l@{}l@{}l@{}l@{}l@{}l}
 &x_0\phantom{)}(&x_0&x_2)&x_5 &x_8&x_5^{-1}&x_5^{-1}(&x_3^{-1}&x_0^{-1})\\
=\phantom{()}&x_0&x_1\phantom{)}(&x_0&x_5)&x_8&x_5^{-1}(&x_5^{-1}&x_0^{-1})&x_2^{-1}\\
=&x_0&x_1&x_4\phantom{)}(&x_0&x_8)(&x_5^{-1}&x_0^{-1})&x_4^{-1}&x_2^{-1}\\
=&x_0&x_1&x_4&x_7\phantom{)}(&x_0&x_0^{-1})&x_4^{-1}&x_4^{-1}&x_2^{-1}\\
=&x_0&x_1&x_4&x_7&   &        &x_4^{-1}&x_4^{-1}&x_2^{-1}
\end{array}
\end{equation*}
Finally, we can also cancel one $x_4,x_4^{-1}$ pair, since there are no
instances of $x_5$ \mbox{or $x_5^{-1}$}:
\begin{equation*}
\begin{array}{l@{}l@{}l@{}l@{}l@{}l@{}l@{}l}
 &x_0&x_1\phantom{)}(&x_4&x_7)\phantom{)}&x_4^{-1}&x_4^{-1}&x_2^{-1}\\
=\phantom{()}&x_0\phantom{()}&x_1&x_6\phantom{(}(&x_4&x_4^{-1})&x_4^{-1}\phantom{(}&x_2^{-1}
\end{array}
\end{equation*}
This gives us the normal form for the element:
\begin{equation*}
x_0\phantom{()} x_1\phantom{()} x_6\phantom{()} x_4^{-1}\phantom{(} x_2^{-1}
\end{equation*}
\end{example}

Though the presentation for $F$ given in theorem 1.3.6 is infinite,
Thompson's group $F$ is actually finitely presented:

\begin{theorem} The elements $x_0$ and $x_1$ generate $F$, with
presentation:
\begin{equation*}
\langle x_0,x_1\mid x_2x_1=x_1x_3,x_3x_1=x_1x_4\rangle
\end{equation*}
where the symbol $x_n$ ($n\geq 2$) stands for the word
$(x_1)^{x_0^{n-1}}$.
\end{theorem}
\begin{proof} Since $x_{n+1}=(x_n)^{x_0}$ for each $n\geq 1$, we have:
\begin{equation*}
x_n=(x_1)^{x_0^{n-1}}
\end{equation*}
for $n\geq 2$. This proves that $x_0$ and $x_1$ generate $F$.

We must now show that the two relations suffice. For $n>k\geq 0$, let
$R_{n,k}$ denote the relation:
\begin{equation*}
x_nx_k=x_kx_{n+1}
\end{equation*}
where the symbol $x_n$ ($n\geq 2$) stands for the word
$(x_1)^{x_0^{n-1}}$. Our task is to show that the relations
$R_{2,1}$ and $R_{3,1}$ imply all the others.

First note that all the relations $R_{n,0}$ are trivially true by the
definition of the symbol $x_n$. Next,
if we know the relation $R_{n,k}$ for $n>k>0$, we can conjugate repeatedly by $x_0$
to prove all the relations $R_{n+i,k+i}$. This puts us in the
following situation:
\begin{equation*}
\xymatrix{
& & & & & \\
& & & R_{4,3} \ar@{=>}[ur] & R_{5,3} \ar@{=>}[ur] & \cdots \\
& & R_{3,2} \ar@{=>}[ur] & R_{4,2} \ar@{=>}[ur] & R_{5,2} \ar@{=>}[ur] & \cdots \\
& *+[F]{R_{2,1}} \ar@{=>}[ur] & *+[F]{R_{3,1}} \ar@{=>}[ur] & R_{4,1} \ar@{=>}[ur] & R_{5,1} \ar@{=>}[ur] & \cdots \\
*+[F]{R_{1,0}} & *+[F]{R_{2,0}} & *+[F]{R_{3,0}} & *+[F]{R_{4,0}} & *+[F]{R_{5,0}} & \cdots
}
\end{equation*}
Arrows in this diagram indicate implication, and boxes indicate relations we
already know.

Our strategy is to deduce the relations $R_{n,1}$ for $n>3$ by
induction, with base case $n=3$. Note that while proving $R_{n,1\text{,}}$
we may use any relation $R_{j,i}$ with the property that $j-i<n-1$. Here's
the calculation:
\begin{equation*}
\begin{array}{llllll}
 &x_2\phantom{(((}    &x_n\phantom{(((} &x_1     &\\
=&x_{n-1}&x_2  &x_1     &\text{(using }R_{n-1,2}\text{)}\\
=&x_{n-1}&x_1 &x_3     &\text{(using }R_{2,1}\text{)}\\
=&x_1    &x_n &x_3     &\text{(using }R_{n-1,1}\text{)}\\
=&x_1   &x_3 &x_{n+1} &\text{(using }R_{n,3}\text{)}\\
=&x_2    &x_1 &x_{n+1}\quad &\text{(using }R_{2,1}\text{)}
\end{array}
\end{equation*}
Cancelling the initial $x_2$'s yields $R_{n,1}$.\end{proof}

Brown and Geoghegan \cite{BrGe} have shown that $F$ has an Eilenberg-MacLane complex with exactly two cells
in each dimension.  Therefore, the group $F$ is infinite-dimensional and has type $F_\infty$.

By the way, the presentation above is one of two canonical finite presentations for
$F$. In the presentation above, the relations $R_{n,0}$ were true ``by
definition'', and the relations $R_{2,1}$ and $R_{3,1}$ were used to deduce
the rest. It is possible instead to assume the relations $R_{n,n-1}$ ``by
definition'', and then use the relations $R_{2,0}$ and $R_{3,0}$ to deduce
the rest. We state the resulting presentation without proof:

\begin{proposition} $F$ has presentation:
\begin{equation*}
\langle x_0,x_1\mid x_2x_0=x_0x_3,x_3x_0=x_0x_4\rangle
\end{equation*}
where the word $x_n$ for $n\geq 2$ is defined inductively by:
\begin{equation*}
x_{n+1}=x_{n-1}^{-1}x_nx_{n-1}\tag*{\qedsymbol}
\end{equation*}
\end{proposition}

\section{The Commutator Subgroup}

In this section, we will prove that the commutator subgroup $\left[F,F\right]$ is simple,
and that every proper quotient of $F$ is abelian.

\begin{proposition} The abelianization of $F$ is $\mathbb{Z}\oplus\mathbb{Z}$.
\end{proposition}
\begin{proof} Abelianizing the standard presentation for $F$ (see theorem 1.3.6) yields:
\begin{equation*}
\langle x_0,x_1,x_2,\ldots\mid x_n+x_k=x_k+x_{n+1}\text{ for }k<n\rangle
\end{equation*}
which is just:
\begin{equation*}
\langle x_0,x_1,x_2,\ldots\mid x_1=x_2=x_3=\cdots\rangle\tag*{\qedhere}
\end{equation*}
\end{proof}

In fact $H_n(F,\mathbb{Z})=\mathbb{Z}\oplus\mathbb{Z}$ for all $n\geq
2$.  See \cite{BrGe}.

There is a nice geometric description of the abelianization:

\begin{proposition}  Define a function $\varphi\colon
F\rightarrow\mathbb{Z}\oplus\mathbb{Z}$ by:
\begin{equation*}
\varphi(f)=\bigl(\log_2f'(0),\log_2f'(1)\bigr)
\end{equation*}
Then $\varphi$ is an epimorphism, and
$\ker(\varphi)=[F,F]$. \end{proposition}
\begin{proof} That $\varphi$ is a homomorphism follows from the chain rule. For
the rest, note that:
\begin{equation*}
\varphi(x_0)=(1,-1)\qquad\text{and}\qquad\varphi (x_1)=(0,-1)
\end{equation*}
Since these vectors generate $\mathbb{Z}\oplus\mathbb{Z}$, $\varphi$ is an
epimorphism.  Since they are linearly independent, the kernel of $\varphi$ is $[F,F]$.\end{proof}

\begin{corollary}  Let $f\in F$. Then $f\in\left[F,F\right]$ if and only
if $f$ is trivial in neighborhoods of $0$ and $1$.
\end{corollary}

Define a \emph{dyadic interval} to be any closed interval with dyadic rational
endpoints.

\begin{proposition} Let $I$ be any dyadic interval, and let
$\PL_2(I)$ be the subgroup of $F$ consisting of all elements with
support in $I$. Then $\PL_2(I)\cong F$.

In particular, there exists a piecewise-linear homeomorphism
$\gamma\colon\left[0,1\right]\rightarrow I$ such that:
\begin{enumerate}
\item  All slopes of $\gamma$ are powers of $2$, and
\item  All breakpoints of $\gamma$ have dyadic rational coordinates.
\end{enumerate}
Any such homeomorphism conjugates $F$ to $\PL_2(I)$.
\end{proposition}
\begin{proof} Clearly $I$ is a union of finitely many standard dyadic intervals.
By choosing a dyadic subdivision of $\left[0,1\right]$ with that same number
of intervals, we can construct the desired homeomorphism $\gamma$. \end{proof}

Note that $\PL_2(I)\subset [F,F]$ when $I\subset(0,1)$.

\begin{theorem} Any nontrivial subgroup of $F$ that is normalized by
$\left[F,F\right]$ contains $\left[F,F\right]$.
\end{theorem}
\begin{proof} Our argument is similar to Epstein's proof \cite{Eps} that various
large homeomorphism groups are simple.  In particular, Epstein proves that the
group of all piecewise-linear, orientation-preserving, compactly-supported
homeomorphisms of $(0,1)$ is simple, and that the group of all orientation-preserving,
compactly-supported diffeomorphisms of $(0,1)$ has simple commutator subgroup.
(It has since been shown that this diffeomorphism group is perfect, and therefore simple.  See \cite{Mat}.)

Let $N$ be a nontrivial subgroup of $F$ normalized by
$\left[F,F\right]$, and let $\eta$ be a nontrivial element of $N$. Since
$\eta$ is not the identity, there exists a sufficiently small dyadic interval
$I\subset(0,1)$ so that $\eta(I)$ is disjoint from
$I$.

Given any $f\in\PL_2(I)$, the commutator
$\left[\eta,f\right]=\eta^{-1}\left(f^{-1}\eta
f\right)=\left(\eta^{-1}f^{-1}\eta\right)f$ is in $N$, has support in
$I\cup\eta(I)$, and agrees with $f$ on $I$. Therefore:
\begin{equation*}
\left[f,g\right]=\left[\left[\eta,f\right],g\right]\in N
\end{equation*}
for any $f,g\in\PL_2(I)$, so $N$ contains every commutator with
support in the interior of $I$.

However, there exist elements $\gamma_1,\gamma_2,\gamma_3,\ldots$ of
$\left[F,F\right]$ such that:
\begin{equation*}
\gamma_1(I)=\left[\frac{1}{4},\frac{3}{4}\right],\quad\gamma_2%
(I)=\left[\frac{1}{8},\frac{7}{8}\right],\quad\gamma_3(I%
)=\left[\frac{1}{16},\frac{15}{16}\right],\quad\ldots
\end{equation*}
Conjugation by these elements shows that $N$ contains every commutator with
support in the interior of any of these intervals, and hence $N$ contains every
element of $[F,F]$.\end{proof}

\begin{corollary} Every proper quotient of $F$ is abelian.
\end{corollary}

\begin{corollary} The commutator subgroup of $F$ is simple.
\end{corollary}

Closely related to Thompson's group $F$ are Thompson's groups $V$ and $T$.  (See section
7.4 for a discussion of these groups.)  These groups are themselves simple (instead of
just having a simple commutator subgroup) and were the first known examples of
infinite, finitely-presented simple groups (see \cite{Hig}).

\section{Open Problems}

Recall the following definition:

\begin{definition} Let $G$ be a group, and let
$\mathcal{P}(G)$ be the collection of all subsets of $G$. We say
that $G$ is \emph{amenable} if there exists a function $\mu\colon
G\rightarrow\left[0,1\right]$ (called a \emph{measure}) with the following
properties:
\begin{enumerate}
\item  $\mu(G)=1$.
\item  If $S$ and $T$ are disjoint subsets of $G$, then
        $\mu(S\cup T)=\mu(S)+\mu(T)$.
\item  If $S\subset G$ and $g\in G$, then
        $\mu(gS)=\mu(S)$.
\end{enumerate}
\end{definition}

See \cite{Wag} for a lengthy discussion of amenability.

\begin{theorem} \
\begin{enumerate}
\item  All finite groups are amenable.
\item  All abelian groups are amenable.
\item  Subgroups and quotients of amenable groups are amenable.
\item  Any extension of an amenable group by an amenable group is amenable.
\item  Any direct union of amenable groups is amenable.\quad\qedsymbol
\end{enumerate}
\end{theorem}

\begin{theorem} Any group that contains a free subgroup of rank two is not amenable.\quad\qedsymbol
\end{theorem}

The amenability of $F$ has been an open problem for several decades:

\begin{question} Is $F$ amenable?
\end{question}

This question was originally motivated by the following considerations.
Let $\AG$ denote the class of all amenable groups, and let $\EG$ denote the
smallest class of groups containing all finite and abelian groups and closed
under subgroups, quotients, extensions, and direct unions (the \emph{elementary
amenable} groups). Let $\NF$ denote the class of all groups that do
not contain a free subgroup of rank two. According to the above theorems:
\begin{equation*}
\EG\quad\subset\quad\AG\quad\subset\quad\NF
\end{equation*}
The question then arises as to whether either of these inclusions is proper.
At the time that the amenability of $F$ was first investigated, there were
no known examples of groups in either $\AG\setminus\EG$ or $\NF\setminus\AG$.
However, $F$ was known to be in the class $\NF\setminus\EG$:

\begin{theorem}  $F$ is not elementary amenable.
\end{theorem}
\begin{proof}  Chou \cite{Chou} has proven the following result concerning the class
EG. Let $\EG_0$ be the class of all finite or abelian groups, and for each
ordinal $\alpha$, let $\EG_\alpha$ be the class of all groups that can be
constructed from elements of $\bigcup_{\beta<\alpha}\EG_\beta$ using
extensions and direct unions. Then each class $\EG_\alpha$ is closed under
subgroups and quotients, and hence:
\begin{equation*}
\EG=\bigcup_\alpha\EG_\alpha
\end{equation*}

Now, $F$ is certainly not in the class $\EG_0$. Furthermore, since $F$
is finitely generated, $F$ cannot be expressed as a nontrivial direct union.
Therefore, we need only show that $F$ cannot arise in some $\EG_\alpha$ as a
nontrivial group extension.

Suppose there were a nontrivial short exact sequence:
\begin{equation*}
N\hookrightarrow F \twoheadrightarrow Q
\end{equation*}
where $N,Q\in\EG_\beta$ for some $\beta<\alpha$. By theorem
1.4.5, $N$ must contain the commutator subgroup $[F,F]$ of
$F$, and therefore $N$ contains a copy of $F$. Since $\EG_\beta$ is closed under
taking subgroups, we conclude that $F\in\EG_\beta$, a contradiction since
$\beta < \alpha$.\end{proof}

\begin{theorem} $F$ does not contain the free group of rank two.
\end{theorem}
\begin{proof}Let $f,g\in F$.  We must show that $f$ and $g$ do not generate a free subgroup.

Assume first that $f$ and $g$ have no common fixed points in $(0,1)$.  Observe then that
any $t\in(0,1)$ can be sent arbitrarily close to $0$ using elements of $\langle f,g\rangle$
(otherwise the infimum of the orbit of $t$ would be a common fixed point of $f$ and~$g$).
In particular, there is an $h\in\langle f,g\rangle$ such that
$[f,g]^h$ has support disjoint from that of $[f,g]$.  Then $[f,g]^h$ and $[f,g]$ commute,
so $f$ and $g$ do not generate a free subgroup.

Now suppose that $f$ and $g$ have a common fixed point in $(0,1)$.  Then the support of
$\langle f,g\rangle$ is the union of the interiors of finitely many dyadic intervals
$I_1,\ldots,I_n$.  This gives us a monomorphism:
\begin{equation*}
\langle f,g\rangle\hookrightarrow \PL_2(I_1)\times\cdots\times
\PL_2(I_n)
\end{equation*}
By the above argument, the image of $\langle f,g\rangle$ in each $\PL_2(I_k)$
is not free, so each of the compositions:
\begin{equation*}
F_2 \twoheadrightarrow \langle f,g\rangle \rightarrow \PL_2(I_k)
\end{equation*}
has a nontrivial kernel.  Then the kernel of the projection
$F_2 \twoheadrightarrow \langle f,g\rangle$ is the intersection of finitely many
nontrivial normal subgroups of $F_2$, and is therefore nontrivial.
\end{proof}

In fact, every nonabelian subgroup of $F$ contains a free abelian group of infinite rank.
See \cite{BrSq} for details.

The status of the classes $\AG\setminus\EG$ and $\EG\setminus\NF$ was
resolved in the 1980's. In particular, Ol'shanskii \cite{Ol} constructed a
nonamenable torsion group, thereby supplying an element of $\NF\setminus\AG$;
and Grigorchuk \cite{Grig} constructed his famous group $\mathcal{G}$ of
intermediate growth, thereby supplying an element of $\AG\setminus\EG$.

However, the amenability of $F$ remains an interesting question, if
only because of its puzzling difficulty. $F$ is certainly the most
well-known group for which amenability is still an issue. It is hoped that
resolving this question will shed new light on either the structure of $F$,
or on the nature of amenability.

In addition to amenability, there are several other fundamental open
questions concerning the Cayley graph of $F$.

\begin{definition}  Let $G$ be a group with finite generating set
$\Sigma$. For each $n$, let $\gamma(n)$ be the number of
elements of $G$ which are products of at most $n$ elements of $\Sigma$.
Then $\gamma(n)$ is called the \emph{growth function} of
$(G,\Sigma)$.

It is known that the limit
$\lim_{n\rightarrow\infty}\sqrt[n]{\gamma(n)}$ always exists (see \cite{dlH}). We
say that $G$ has \emph{exponential growth} if this limit is positive, and
\emph{subexponential growth} if this limit is $0$.
\end{definition}

It turns out that the classification of $G$ as having exponential or
subexponential growth does not depend on the finite generating set $\Sigma$
(see \cite{dlH}).

\begin{proposition}  Any group with subexponential growth is amenable.
\end{proposition}

Letting $\SG$ denote all groups with subexponential growth, we have:
\begin{equation*}
\SG\quad\subset\quad\AG\quad\subset\quad\NF
\end{equation*}

\begin{proposition}  The submonoid of $F$ generated by
$\left\{x_0^{-1},x_1\right\}$ is free.
\end{proposition}
\begin{proof}  Consider any word in $x_0^{-1}$ and $x_1$:
\begin{equation*}
x_1^{a_1}x_0^{-1}x_1^{a_2}x_0^{-1}\cdots x_0^{-1}x_1^{a_n}
\end{equation*}
where $a_0,\ldots,a_n\geq 0$. We can put this element into normal form by
moving the $x_0^{-1}$'s to the right:
\begin{equation*}
x_1^{a_1}x_2^{a_2}\cdots x_n^{a_n}x_0^{-(n-1)}
\end{equation*}
Since the normal form is different for different values of $a_1,\ldots,a_n$,
every word in $x_0^{-1}$ and $x_1$ represents a different element of
$F$.\end{proof}

\begin{corollary}  The group $F$ has exponential growth.\end{corollary}

It would be interesting to determine the exact growth rate of $F$ with
respect to the $\left\{x_0,x_1\right\}$ generating set. Perhaps more
interesting is the following question:

\begin{question} Let $\gamma(n)$ be the growth rate of $F$
with respect to the $\left\{x_0,x_1\right\}$ generating set, and let:
\begin{equation*}
\Gamma(t)=\sum_{n=0}^\infty\gamma(n)t^n
\end{equation*}
Is $\Gamma(t)$ a rational function?
\end{question}

See \cite{dlH} for details about groups with rational growth functions.

Finally, it is not known whether $F$ is automatic.  Recall the following definitions:

\begin{definition}
Let $G$ be a group with finite generating set $\Sigma$, and let $\Gamma$ denote the corresponding
Cayley graph.  A \emph{combing} of $G$ is a choice, for each $g\in G$, of a path in $\Gamma$ from the identity
vertex to $g$ (i.e.~a word in the generators that multiplies to $g$).
\end{definition}

If $\sigma_1 \sigma_2 \cdots \sigma_\ell$ and $\tau_1 \tau_2 \cdots \tau_m$ are words in $\Sigma$, the
\emph{synchronous distance} between these words is:
\begin{equation*}
\max_{n\in\mathbb{N}} d(\sigma_1\cdots\sigma_n,\tau_1\cdots\tau_n)
\end{equation*}
where $d$ denotes the distance function in the Cayley graph, $\sigma_i=1$ for $i>\ell$, and $\tau_i=1$ for $i>m$.

\begin{definition}
Let $G$ be a group with finite generating set $\Sigma$.  A combing of $G$ has the \emph{fellow traveller
property} if, given any two elements of $G$ a distance one apart, the corresponding combing
paths have synchronous distance at most $k$, for some fixed $k\in\mathbb{N}$.
\end{definition}

\begin{definition}
Let $G$ be a group with finite generating set $\Sigma$.  We say that $G$ is \emph{automatic} if there exists a
combing of $G$ that has the fellow traveller property and whose set of combing paths is a regular
language over $\Sigma\cup\Sigma^{-1}$.
\end{definition}

It turns out that this definition does not depend on the finite generating set $\Sigma$.  See
\cite{ECH} for a thorough introduction to automatic groups.

\begin{question} Is $F$ automatic?
\end{question}

In \cite{Guba2}, V. Guba shows that the Dehn function of $F$ is quadratic.  (Any automatic group has
linear or quadratic Dehn function.)  In section 6.4, we will show that $F$ no \emph{geodesic} combing of
$F$ (with respect to the $\{x_0,x_1\}$ generating set) has the fellow traveller property.

\section{Alternate Descriptions}

Thompson's Group $F$ has arisen naturally in a variety of different
contexts, and this has led to several different ways of defining the group.
We have adopted the ``homeomorphism'' point of view: every element of $F$ is
a piecewise-linear homeomorphism of $\left[0,1\right]$. In this section, we
will describe three alternate definitions of $F$, and explain why they are
equivalent to the homeomorphism definition.

We begin with Thompson's original definition:

\subsubsection{Associative Laws}

An \emph{associative law} is any rule for rearranging a parenthesized
expression. For example, one associative law is the rule:
\begin{equation*}
x_0\colon\quad(ab)c\quad\rightarrow\quad a(bc)
\end{equation*}
This rule can be applied to any expression whose left part is nontrivial,
e.g.:
\begin{equation*}
\bigl((ab)(cd)\bigr)\bigl(e(fg)\bigr)%
\quad\rightarrow\quad(ab)\bigl((cd)(e(fg))\bigr)
\end{equation*}
However, $x_0$ can only be applied to the top level of an expression. In
particular, the rearrangement:
\begin{equation*}
x_1\colon\quad a\bigl((bc)d\bigr)\quad\rightarrow\quad a\bigl(b(cd)\bigr)
\end{equation*}
is not an application of $x_0$.

We can compose two associative laws by performing one and then the
other. For example, starting with the expression:
\begin{equation*}
\bigl(a(bc)\bigr)d
\end{equation*}
we can perform $x_0$:
\begin{equation*}
a\bigl((bc)d\bigr)
\end{equation*}
and then $x_1$:
\begin{equation*}
a\bigl(b(cd)\bigr)
\end{equation*}
This yields a new, composite law:
\begin{equation*}
x_0x_1\colon\quad\bigl(a(bc)\bigr)d\quad\rightarrow\quad a\bigl(b(cd)\bigr)
\end{equation*}

\begin{theorem}  The set of all associative laws forms a group under
composition, and this group is isomorphic with Thompson's group $F$.
\end{theorem}
\begin{proof}[Sketch of Proof]  Any parenthesized expression corresponds to a finite binary tree.
For example, the expression:
\begin{equation*}
\bigl(a(bc)\bigr)(de)
\end{equation*}
corresponds to the tree:
\begin{center}
\includegraphics{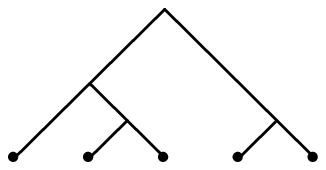}
\end{center}
A rearrangement of parentheses is really just a tree diagram for an element
of~$F$. Observe that two rearrangements are instances of the same
associative law if and only if the corresponding tree diagrams represent the
same element of $F$. Observe also that the rule for composing associative
laws is the same as the rule for multiplying tree diagrams.\end{proof}

\subsubsection{Automorphisms of Cantor Algebras}

This point of view was introduced by Galvin and Thompson, and used by Brown in \cite{Bro} to
show that $F$ has type F$_\infty$. Brown's paper actually considers a whole
class of groups defined by Higman \cite{Hig} using automorphisms of Cantor algebras (which Brown
refers to as J\a'onsson-Tarski algebras), and proves finiteness results for
all of them.

\begin{definition} A \emph{Cantor algebra} is a set $A$ together with a
bijection \linebreak $\alpha\colon A\rightarrow A\times A$.
\end{definition}

We will denote the components of $\alpha(a)$ by $(a)_0$ and $(a)_1$.  Also, if
$a,b\in A$, we will denote by $ab$ the element of $A$ satisfying
\begin{equation*}
(ab)_0=a\qquad\text{ and }\qquad(ab)_1=b
\end{equation*}

Using the theory of universal algebras, it is easy to show that there exists a free
Cantor Algebra $A(X)$ over any set $X$.  This algebra may be constructed explicitly as follows
(see \cite{Hig}):
\begin{enumerate}
\item[] \textbf{Base Case} Given any $x\in X$ and any word $\epsilon_1\cdots\epsilon_n$
in $\{0,1\}$, define a corresponding element $x_{\epsilon_1\cdots\epsilon_n}$.  By definition,
these elements satisfy:
\begin{equation*}
(x_{\epsilon_1\cdots\epsilon_n})_0=x_{\epsilon_1\cdots\epsilon_n 0}\qquad\text{ and }\qquad
(x_{\epsilon_1\cdots\epsilon_n})_1=x_{\epsilon_1\cdots\epsilon_n 1}
\end{equation*}
\item[] \textbf{Induction Step} Suppose we have defined $a,b\in A(X)$ but we have not yet defined
any element $c$ such that $(c)_0=a$ and $(c)_1=b$.  Then we define an element $ab$ satisfying:
\begin{equation*}
(ab)_0=a\qquad\text{ and }\qquad(ab)_1=b
\end{equation*}
\end{enumerate}

In the case where $X$ is a singleton set, it is helpful to think of the unique element $x$
as the interval $\{[0,1]\}$, and the elements $x_{\epsilon_1\cdots\epsilon_n}$ as standard dyadic
subintervals of $\{[0,1]\}$ (e.g. $x_0=[0,1/2]$ and $x_{011}=[3/8,4/8]$).  Each element
created during the induction step can be thought of as a mapping from some disjoint union
of these intervals linearly onto the intervals of some dyadic subdivision of $\{[0,1]\}$.
(For example, $(x_0 x_0 )x_1$ maps the disjoint union $[0,1/2] \uplus [0,1/2]
\uplus [1/2,1]$ onto $[0,1/4] \uplus [1/4,1/2] \uplus [1/2,1]$.)

\begin{definition}\emph{Thompson's group $V$} is the automorphism group of the free
Cantor algebra $A(\{x\})$.
\end{definition}

Observe that the free algebra generated by a singleton set $\{x\}$ is isomorphic to the free
algebra generated by a doubleton set $\{a,b\}$, under the mapping:
\begin{equation*}
a\rightarrow (x)_0\qquad\text{ and }\qquad b \rightarrow (x)_1
\end{equation*}
That is, \emph{the elements $\{x_0,x_1\}$ are a basis for the free Cantor algebra $A(\{x\})$}.
More generally, if $\{a_1\ldots a_n\}$ is a basis for a free Cantor algebra, a \emph{simple
expansion} of this basis is obtained by replacing the element $a_i$ with the elements
$(a_i)_0$ and $(a_i)_1$.  The inverse of an expansion (i.e. picking any two elements of a basis
and replacing them with their product) is called a \emph{simple contraction}.  Higman
\cite{Hig} proves that any basis for A($x$) can be obtained from $\{x\}$ by a sequence
of simple expansions followed by a sequence of simple contractions.
    Therefore, an element of $V$ can be represented by a pair of binary trees together with some bijection
of their leaves:
\begin{center}
\includegraphics{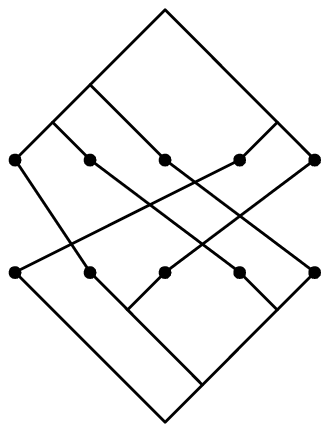}
\end{center}
The root of the top tree represents the basis $\{x\}$, while the root of the bottom tree
represents the image of $x$ under the automorphism.  The trees and permutation represent
the expansions and contractions necessary to get from $x$ to its image.

Thompson's group $F$ is the subgroup of $V$ consisting of tree diagrams whose permutation part is trivial.
In particular, given an ordered basis:
\begin{equation*}
\{a_1,\ldots,a_n\}
\end{equation*}
for a Cantor algebra, an \emph{simple ordered expansion} of this basis is any basis of the form:
\begin{equation*}
\{a_1,\ldots,(a_i)_0,(a_i)_1,\ldots,a_n\}
\end{equation*}
A \emph{simple ordered contraction} is the inverse of a simple ordered expansion.  An automorphism
$f\colon A(\{x\}) \rightarrow A(\{x\})$ is \emph{order-preserving} if the basis $\{f(x)\}$
can be obtained from the basis $\{x\}$ by a sequence of simple ordered expansions and simple ordered contractions.

\begin{theorem}The group of order-preserving automorphisms of $A(\{x\})$ is isomorphic
with Thompson's group $F$.\quad\qedsymbol
\end{theorem}

\subsubsection{$F$ as the Universal Conjugacy Idempotent}

The group $F$ was independently rediscovered by Freyd and Heller in 1969 (see \cite{FrHe}) during their
investigation into homotopy-idempotent homeomorphisms of topological spaces. They (and also
independently Dydak \cite{Dy}) developed $F$ as the universal example of a group with a
conjugacy-idempotent endomorphism. This point of view
motivated Brown and Geoghegan \cite{BrGe} to investigate the finiteness properties of $F$.

\begin{definition} Let $G$ be a group. An endomorphism $\varphi$ of $G$ is \emph{conjugacy idempotent}
if there exists a $c\in G$ such that:
\begin{equation*}
\varphi^2(g)=c^{-1}\varphi(g)c
\end{equation*}
for every $g\in G$.
\end{definition}

We shall refer to the element $c$ as a \emph{conjugator} for $\varphi$.

Let $\sigma\colon F\rightarrow F$ denote the shift endomorphism:
\begin{equation*}
\sigma(x_0)=x_1\text{,}\qquad\sigma(x_1)=x_2\text{,}%
\quad\sigma(x_2)=x_3\text{,}\quad\ldots
\end{equation*}
Then $\sigma$ is conjugacy idempotent, with conjugator $x_0$:
\begin{equation*}
\sigma^2(f)=x_0^{-1}\sigma(f)x_0
\end{equation*}
Furthermore, the triple $(F,\sigma,x_0)$ has the following universal property:

\begin{theorem}Suppose we are given any triple $(G,\varphi,c)$, where $G$ is a group and
$\varphi$ is a conjugacy idempotent on $G$ with conjugator $c$.
Then there exists a unique homomorphism $\pi\colon F\rightarrow G$ such that $\pi(x_0)=c$ and the following diagram
commutes:
\begin{equation*}
\def\labelstyle{\textstyle}
\xymatrix{F \ar[d]_{\pi} \ar[r]^{\sigma} & F \ar[d]^{\pi} \\ G \ar[r]_{\varphi} & G }
\end{equation*}
\end{theorem}
\begin{proof} Define $\pi$ as follows:
\begin{equation*}
\pi(x_0)=c,\quad\pi(x_1)=\varphi(c),\quad\pi%
(x_2)=\varphi^2(c),\quad\ldots
\end{equation*}
We must show that this respects the relations in $F$. Well, for $n>k$,
\begin{equation*}
\begin{aligned}[t]
\pi(x_k)^{-1}\,\pi(x_n)\,\pi(x_k)&=\varphi^k%
(c)^{-1}\,\varphi^n(c)\,\varphi^k(c)\\
&=\varphi^k\bigl(c^{-1}\varphi^{n-k}(c)\,c\bigr)\\
&=\varphi^k\bigl(\varphi^{n-k+1}(c)\bigr)\\
&=\varphi^{n+1}(c)\\
&=\pi(x_{n+1})
\end{aligned}
\end{equation*}
The homomorphism $\pi$ is therefore well-defined, and clearly it has the
required properties. Furthermore, $\pi$ is unique because it is required
to satisfy $\pi(x_0)=c$, and then the definitions of $\pi(x_1),\pi(x_2),\ldots$ follow from the
commutative diagram.\end{proof}

Freyd and Heller were interested in conjugacy idempotents because of their relationship with topology.
If $Y$ is a topological space with basepoint $y_0$, a \emph{homotopy idempotent} on $Y$ is a map
$g\colon(Y,y_0)\rightarrow(Y,y_0)$ such that $g^2$ is freely homotopic to $g$.
It is easy to check that a homotopy idempotent map induces a conjugacy idempotent endomorphism
on $\pi_1(Y,y_0)$.  Freyd and Heller were interested in the question of whether every homotopy idempotent
\emph{splits}, i.e. can be written as $hk$ where $kh \simeq \text{id}$ (see \cite{FrHe}).  They answered this
question in the negative by showing that the conjugacy idempotent $\sigma\colon F \rightarrow F$ does not split, and then
exporting this result back to the homotopy category.  (The same result was independently obtained by Dyadak \cite{Dy}.)

\chapter{One-Way Forest Diagrams}

In this chapter we introduce \emph{one-way forest diagrams} for elements of $F$.  These
diagrams have the same relationship to a certain action of $F$ on the positive real line
that tree diagrams have to the standard action of $F$ on the unit interval.  The advantage
is that forest diagrams are somewhat simpler, especially in their interactions with the
generating set $\{x_0,x_1,x_2,\ldots\}$.

The existence of one-way forest diagrams was noted by K. Brown in \cite{Bro}, but to our
knowledge no one has ever before used them to study $F$.

\section{The Group $\PL_2 (\mathbb{R}_+)$}

Let $\PL_2(\mathbb{R}_{+})$ be the group of all piecewise-linear
self-homeomorphisms $f$ of $[0,\infty)$ satisfying the following conditions:
\begin{enumerate}
\item  Each linear segment of $f$ has slope a power of $2$.
\item  $f$ has only finitely many breakpoints, each of which has dyadic
rational coordinates.
\item  The rightmost segment of $f$ is of the form:
\begin{equation*}
f(t)=t+m
\end{equation*}
for some integer $m$.
\end{enumerate}

\begin{proposition}  $\PL_2(\mathbb{R}_{+})$ is isomorphic with $F$.
\end{proposition}
\begin{proof}  Let $\psi\colon[0,\infty)\rightarrow[0,1)$ be the
piecewise-linear homeomorphism that maps the intervals:
\begin{center}
\includegraphics{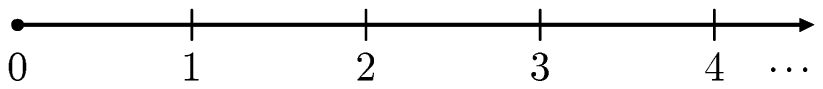}
\end{center}
linearly onto the intervals:
\begin{center}
\includegraphics{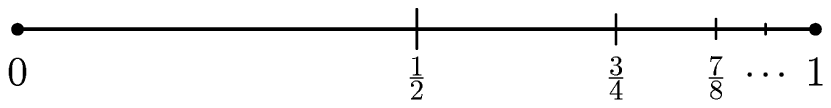}
\end{center}
Then $f\mapsto\psi f\psi^{-1}$ is the desired isomorphism
$F\rightarrow\PL_2(\mathbb{R}_{+})$. In particular, if $f$ has
slope $2^m$ at $t=1$, then the final linear segment of $\psi f\psi^{-1}$ will
be $t\mapsto t+m$.\end{proof}

Under this isomorphism, each generator $x_n$ of $F$ maps to the
piecewise-linear function $x_n\colon[0,\infty)\rightarrow[0,\infty)$ satisfying:
\begin{enumerate}
\item  $x_n$ is the identity on $[0,n]$.
\item  $x_n$ sends $\left[n,n+1\right]$ linearly onto $\left[n,n+2\right]$.
\item  $x_n(t)=t+1$ for $t\geq n+1$.
\end{enumerate}

\section{Forest Diagrams for Elements of $\PL_2 (\mathbb{R}_+)$}

We think of the positive real line as being pre-subdivided as follows:
\begin{center}
\includegraphics{tpict30}
\end{center}
A \emph{dyadic subdivision} of $[0,\infty)$ is any subdivision obtained by cutting
finitely many of these intervals in half, and then cutting finitely many of
the resulting intervals in half, etc.

\begin{proposition}  Let $f\in\PL_2(\mathbb{R}_{+})$. Then
there exist dyadic subdivisions $\mathcal{D},\mathcal{R}$ of $[0,\infty)$
such that $f$ maps each interval of $\mathcal{D}$ linearly onto an interval
of $\mathcal{R}$.$\quad\qedsymbol$
\end{proposition}

A \emph{binary forest} is a sequence $(T_0,T_1,\ldots)$ of finite
binary trees:
\begin{center}
\includegraphics{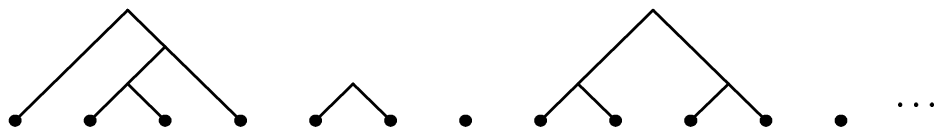}
\end{center}
A binary forest is \emph{bounded} if only finitely many of the trees $T_i$ are
nontrivial.

Every bounded binary forest corresponds to some dyadic subdivision of
the positive real line. For example, the forest above corresponds to the
subdivision:
\begin{center}
\includegraphics{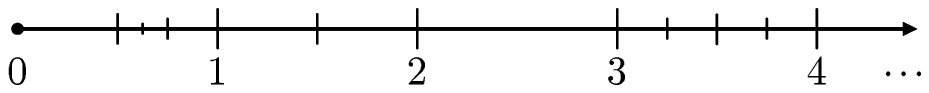}
\end{center}
Each tree $T_i$ represents an interval $\left[i,i+1\right]$, and each leaf
represents an interval of the subdivision.

Combining this with proposition 2.2.1, we see that any
$f\in\PL_2(\mathbb{R}_{+})$ can be represented by a pair of
bounded binary forests. This is called a \emph{(one-way) forest diagram} for $f$

\begin{example}  Let $f$ be the element of $\PL_2(\mathbb{R}_{+})$ with graph:
\begin{center}
\includegraphics{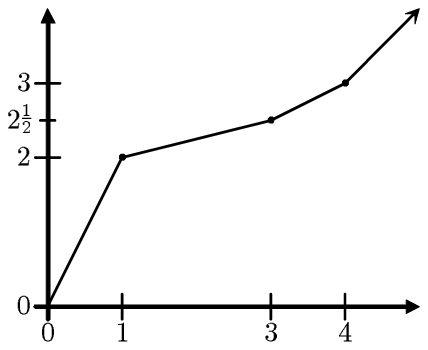}
\end{center}
Then $f$ has forest diagram:
\begin{center}
\includegraphics{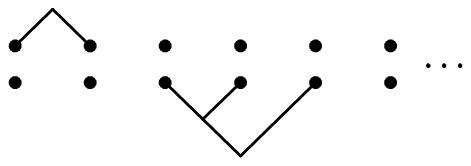}
\end{center}
Again, we have aligned the two forests vertically so that the corresponding
leaves match up.
\end{example}

\begin{example} Here are the forest diagrams for the generators $x_0,x_1,x_2,\ldots$:
\begin{center}
\includegraphics{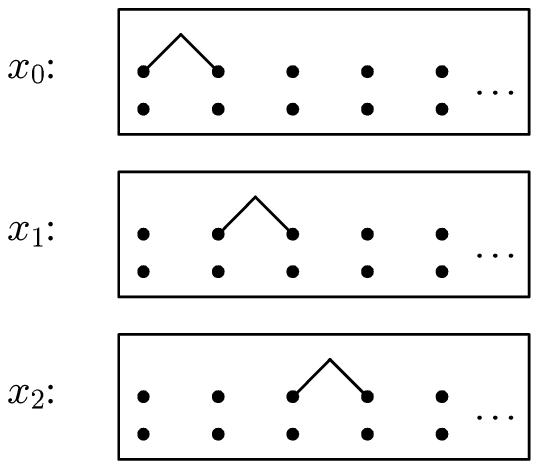}
\end{center}
\end{example}

Of course, there are several forest diagrams for each element of
$\PL_2(\mathbb{R}_{+})$. In particular, it is possible to delete
an opposing pair of carets:
\begin{center}
\includegraphics{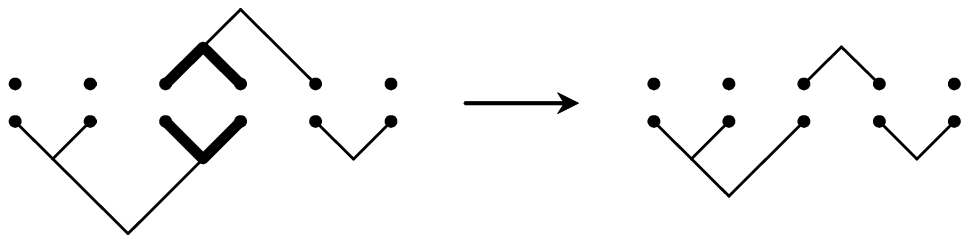}
\end{center}
without changing the resulting homeomorphism. This is called a \emph{reduction} of
a forest diagram. A forest diagram is \emph{reduced} if it does not have any
opposing pairs of carets.

\begin{proposition} Every element of $\PL_2(\mathbb{R}_{+})$
has a unique reduced forest diagram.\quad\qedsymbol\end{proposition}

\begin{remark} From this point forward, we will omit all the trivial trees
on the right side of a forest diagram, as well as the ``$\cdots$'' indicators.
\end{remark}

\begin{remark}  It is fairly easy to translate between tree diagrams and
forest diagrams. Given a tree diagram:
\begin{center}
\includegraphics{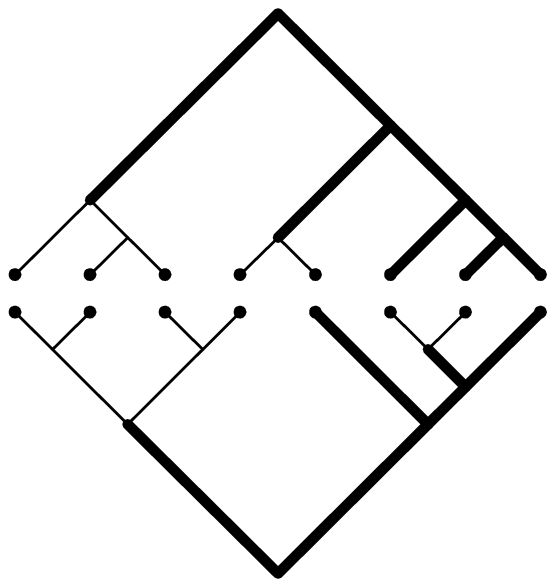}
\end{center}
we simply remove the right stalk of each tree to get the corresponding forest
diagram:
\begin{center}
\includegraphics{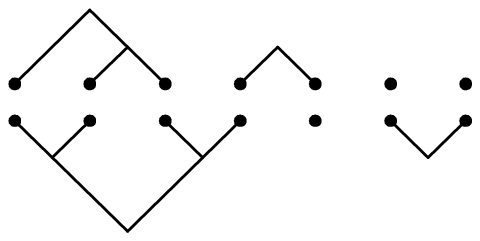}
\end{center}
\end{remark}

\section{The Action of the Generators}

The action of the generators $\left\{x_0,x_1,x_2,\ldots\right\}$ on forest
diagrams is particularly nice:

\begin{proposition}  Let $\mathfrak{f}$ be a forest diagram for some $f\in
F$. Then a forest diagram for $x_nf$ can be obtained by attaching a caret
to the roots of trees $n$ and $(n+1)$ in
the top forest of $\mathfrak{f}$.$\quad\qedsymbol$\end{proposition}

Note that the forest diagram given for $x_nf$ may not be reduced, even if we
started with a reduced forest diagram $\mathfrak{f}$. In particular, the
caret that was created could oppose a caret in the bottom forest. In this
case, left-multiplication by $x_n$ effectively ``cancels'' the bottom caret.

\begin{example}  Let $f\in F$ have forest diagram:
\begin{center}
\includegraphics{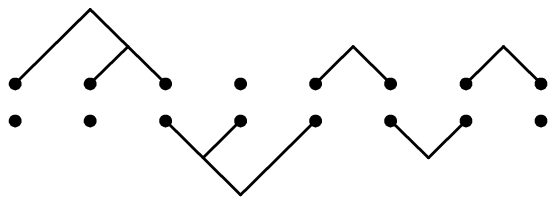}
\end{center}
Then $x_0f$ has forest diagram:
\begin{center}
\includegraphics{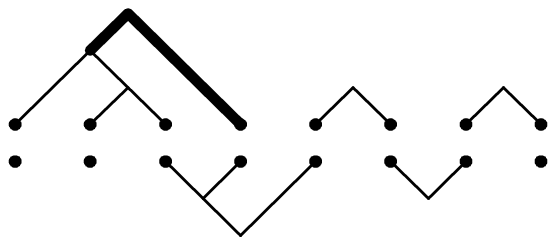}
\end{center}
$x_1f$ has forest diagram:
\begin{center}
\includegraphics{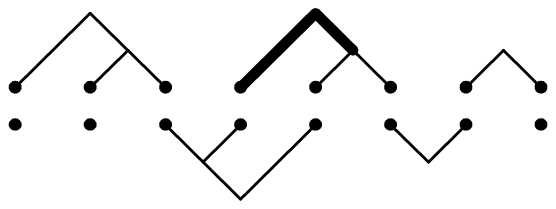}
\end{center}
and $x_2f$ has forest diagram:
\begin{center}
\includegraphics{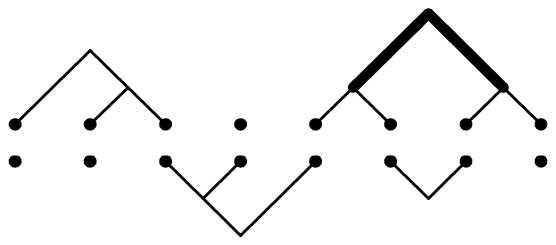}
\end{center}
\end{example}

\begin{example}  Let $f\in F$ have forest diagram:
\begin{center}
\includegraphics{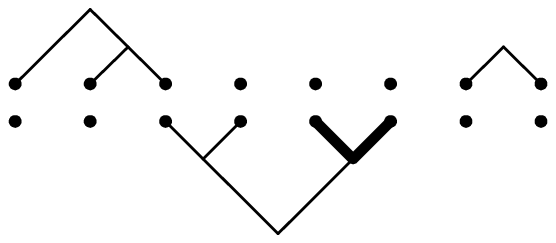}
\end{center}
Then $x_2f$ has forest diagram:
\begin{center}
\includegraphics{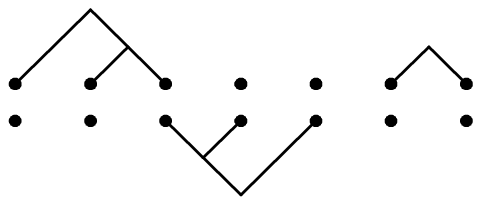}
\end{center}
Note that left-multiplication by $x_2$ cancelled the highlighted bottom caret.
\end{example}

\begin{proposition}  Let $\mathfrak{f}$ be a forest diagram for some $f\in
F$. Then a forest diagram for $x_n^{-1}f$ can be obtained by ``dropping a
negative caret'' at position $n$. If tree $n$ is nontrivial, the
negative caret cancels with the top caret of this tree. If the tree $n$
is trivial, the negative caret ``falls through'' to the bottom forest,
attaching to the specified leaf.\quad\qedsymbol
\end{proposition}

\begin{example} Let $f$ and $g$ be the elements of $F$ with forest
diagrams:
\begin{center}
\includegraphics{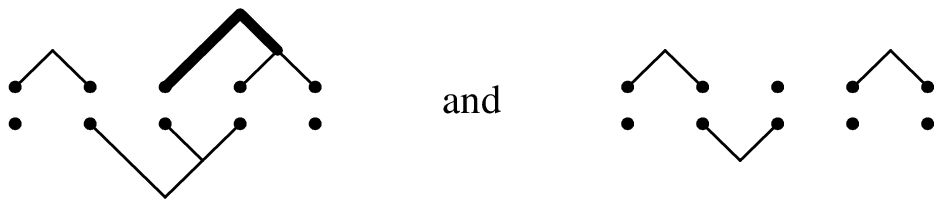}
\end{center}
Then $x_1^{-1}f$ and $x_1^{-1}g$ have forest diagrams:
\begin{center}
\includegraphics{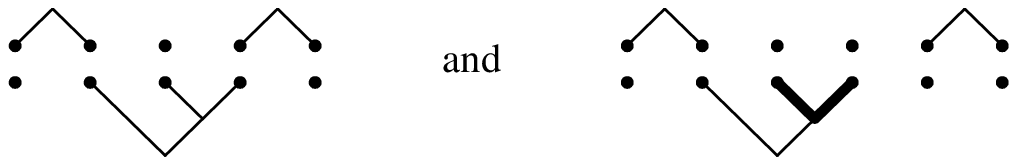}
\end{center}
In the first case, the $x_1^{-1}$ simply removed a caret from the top tree. In the second case, there was no caret on the top to remove, so a new caret
was attached to the leaf on the bottom. Note that this creates a new column
in the forest diagram.
\end{example}

\section{Positive Elements and Normal Forms}

In this section, we use forest diagrams to derive the standard presentation for Thompson's group $F$
(previously stated as theorem 1.3.6) and find a normal form for elements of $F$
(previously stated as theorem 1.3.7).  The proof involves first understanding the
structure of the positive monoid, and then extending this understanding to all of $F$.

Recall that the \emph{positive submonoid} of $F$ generated by
$\left\{x_0,x_1,x_2,\ldots\right\}$.

\begin{proposition}  Let $f\in F$, and let $\mathfrak{f}$ be its reduced
forest diagram. Then $f$ is positive if and only if the bottom forest of
$\mathfrak{f}$ is trivial.\quad\qedsymbol
\end{proposition}

The positive monoid can be thought of as a \emph{monoid of binary forests}.  Each element
corresponds to a binary forest, and two forests can multiplied by ``stacking them'',
i.e. by attaching the leaves of the first forest to the roots of the second.

\begin{example}Let $f$ have forest diagram:
\begin{center}
\includegraphics{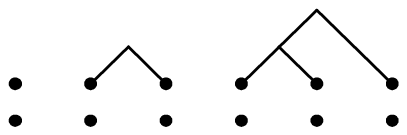}
\end{center}
and let $g$ have forest diagram:
\begin{center}
\includegraphics{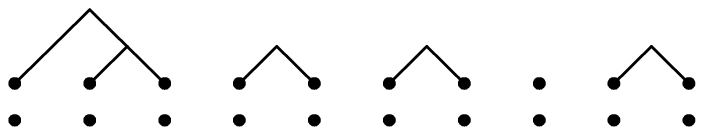}
\end{center}
Then $fg$ has forest diagram
\begin{center}
\includegraphics{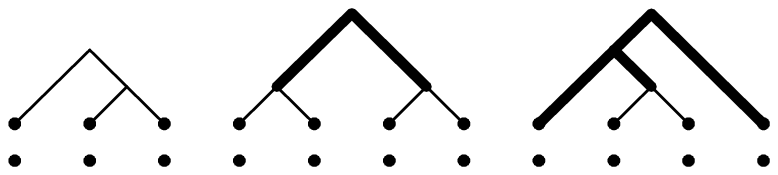}
\end{center}
\end{example}

Each word for a positive element corresponds to an ordering
of the carets of the forest diagram (namely, the order in which the carets are
constructed).

\begin{example}
Let $f$ be the element:
\begin{center}
\includegraphics{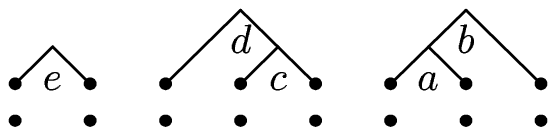}
\end{center}
If we build the carets of $f$ from right to left ($a$--$b$--$c$--$d$--$e$), we get the word:
\begin{equation*}
f = x_0 x_2 x_3 x_5 x_5
\end{equation*}
(Note that the order of the generators in the word is the \emph{opposite} of the order
in which the carets are created, since our primary operation is \emph{left}-multiplication.)

There are many other words for $f$ (30 in all), corresponding to different orderings of the carets.
For example, if we build carets from left to right ($e$--$c$--$d$--$a$--$b$), we get the word:
\begin{equation*}
f = x_2 x_2 x_1 x_2 x_0
\end{equation*}
If we build carets in the order $c$--$e$--$a$--$d$--$b$, we get the word:
\begin{equation*}
f = x_2 x_1 x_3 x_0 x_3
\end{equation*}
\end{example}

In general, we get the \emph{normal form} for a positive element by building the carets
of its forest diagram from right to left:

\begin{theorem}Every positive element can be expressed uniquely in the form:
\begin{equation*}
x_{i_1} \cdots x_{i_n}
\end{equation*}
where $i_1 \leq \cdots \leq i_n$.\quad\qedsymbol
\end{theorem}
\begin{proof}See the following section for a more rigorous proof of this theorem.\end{proof}

The normal form is usually written:
\begin{equation*}
x_0^{a_0}x_1^{a_1}\cdots x_n^{a_n}
\end{equation*}
The numbers $a_0,\ldots,a_n$ are called the \emph{exponents} of the element.  To determine
the exponents of a positive element from its forest diagram, it is helpful to draw the carets so
that left edges are vertical, like this:
\begin{center}
\includegraphics{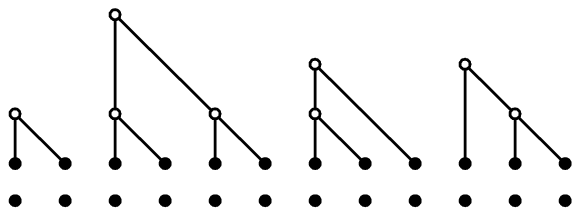}
\end{center}
When the forest diagram is drawn in this fashion, the exponent $a_i$ appears as the number
of nodes sitting directly above the $i\text{'th}$ leaf.
For example, the above element has normal form:
\begin{equation*}
x_0 \, x_2^2 \, x_4 \, x_6^2 \, x_9 \, x_{10}
\end{equation*}

\begin{proposition}The positive monoid has presentation:
\begin{equation*}
\langle x_0,x_1,x_2,\ldots \mid x_n x_k = x_k x_{n+1}\text{ for }k<n \rangle
\end{equation*}
\end{proposition}
\begin{proof}The given relations clearly hold: they arise from the two different ways
of constructing the element:
\begin{center}
\includegraphics{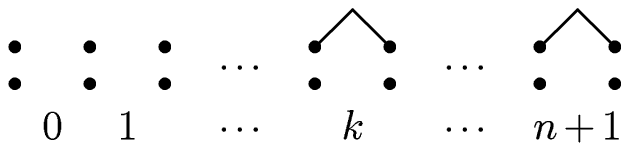}
\end{center}
To show that these relations suffice, we simply observe that any word can be put into
normal form by applying the operations:
\begin{equation*}
x_n x_k \; \rightarrow \; x_k x_{n+1} \qquad \text{(}k<n\text{)}\tag*{\qedhere}
\end{equation*}
\end{proof}

To show that the same presentation holds for $F$, we observe that $F$ is the
\emph{group of fractions} of its positive monoid.

\begin{definition}Let $M$ be any monoid.  A \emph{group of right fractions} for $M$ is a
group $G$ containing $M$ with the property that any element of $G$ can be written
as $pq^{-1}$ for some $p,q\in M$.
\end{definition}

It is easy to determine the presentation for a group of fractions:

\begin{proposition}Suppose that $M$ is a monoid with group of right fractions $G$.  Then
any presentation for $M$ is a presentation for $G$.
\end{proposition}
\begin{proof}This follows immediately from the fact that $G$ is the universal group to which $M$ maps
homomorphically.  See \cite{ClPr}.\end{proof}

\begin{theorem}Thompson's group $F$ has presentation:
\begin{equation*}
\langle x_0,x_1,x_2,\ldots \mid x_n x_k = x_k x_n+1\text{ for }k<n \rangle\tag*{\qedsymbol}
\end{equation*}
\end{theorem}

Next we wish to derive a normal form for elements of $F$. The key is to understand when the forest
diagram corresponding to an expression:
\begin{equation*}
x_0^{a_0} x_1^{a_1} \cdots x_n^{a_n} x_n^{-b_n} \cdots x_1^{-b_1} x_0^{-b_0}
\end{equation*}
is reduced.

An \emph{exposed caret} in a forest is a caret whose children are both leaves:
\begin{center}
\includegraphics{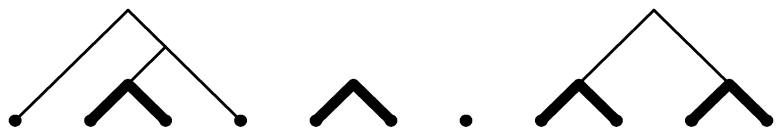}
\end{center}
A forest diagram is reduced if and only if it does not contain a matching pair
of exposed carets.

\begin{lemma}
Let $x_{i_1} \cdots x_{i_n}$ be the normal form for a given positive element $f$.  Then
the caret constructed by $x_{i_k}$ is exposed if and only if $k=n$ or $i_k < i_{k+1} - 1$.
\end{lemma}
\begin{proof}A caret is exposed if and only if it is allowed to be the first caret created
when constructing $f$.  Hence, $x_{i_n}$ is exposed if and only if we can move that generator
all the way to the right in the word for $f$, using operations of the form:
\begin{equation*}
x_k x_n \; \rightarrow x_{n-1} x_k \qquad \text{(}k<n-1\text{)}
\end{equation*}
Such movement is possible if and only if $i_k < i_{k+1} -1$.\end{proof}

\begin{theorem}[Normal Form]
Every element of $F$ can be expressed uniquely in the form:
\begin{equation*}
x_0^{a_0} \cdots x_n^{a_n} x_n^{-b_n} \cdots x_0^{-b_0}
\end{equation*}
where exactly one of $a_n,b_n \neq 0$ and
\begin{equation*}
a_i \neq 0 \text{ and } b_i \neq 0 \quad \Rightarrow \quad
a_{i+1} \neq 0 \text{ or } b_{i+1} \neq 0
\end{equation*}
\end{theorem}
\begin{proof}The top forest has an exposed caret in the $i$'th position if and only if
$a_i \neq 0$ and $a_{i+1} = 0$, and the bottom forest has an exposed caret in the $i$'th
position if and only if $b_i \neq 0$ and $b_{i+1} = 0$.  As long as these never happen
simultaneously, the above expression will represent a reduced forest diagram.\end{proof}

\section{Word Graphs and Anti-Normal Form}

The \emph{word graph} for a positive element $f$ is the directed graph whose vertices are
words for $f$ in the generators ${x_0,x_1,x_2,\ldots}$ and whose edges represent moves
of the form:
\begin{equation*}
x_nx_k\;\rightarrow\;x_kx_{n+1}\qquad\text{(}n>k\text{)}
\end{equation*}

\begin{example}  Let $f$ be the element:
\begin{center}
\includegraphics{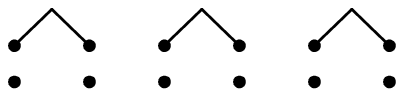}
\end{center}
Then $f$ has word graph:
\begin{center}
\includegraphics{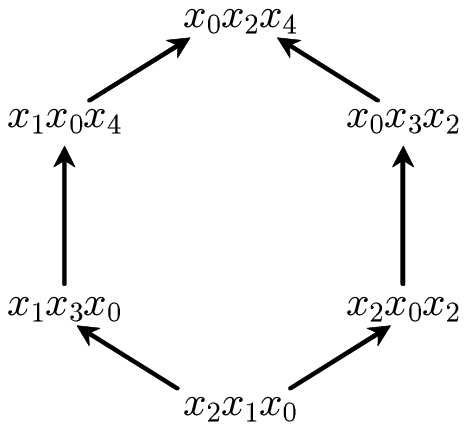}
\end{center}
\end{example}

\begin{example}  Let $f$ be the element:
\begin{center}
\includegraphics{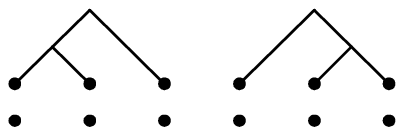}
\end{center}
Then $f$ has word graph:
\begin{center}
\includegraphics{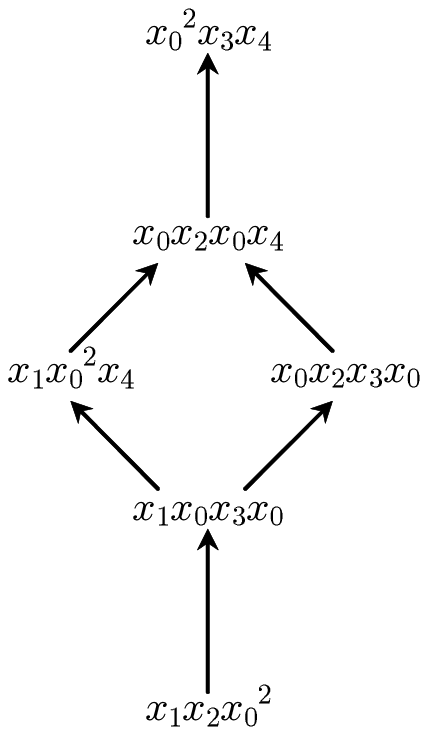}
\end{center}
\end{example}

\begin{example}  Let $f$ be the element:
\begin{center}
\includegraphics{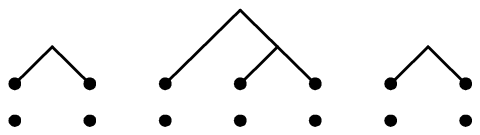}
\end{center}
Then $f$ has word graph:
\begin{center}
\includegraphics{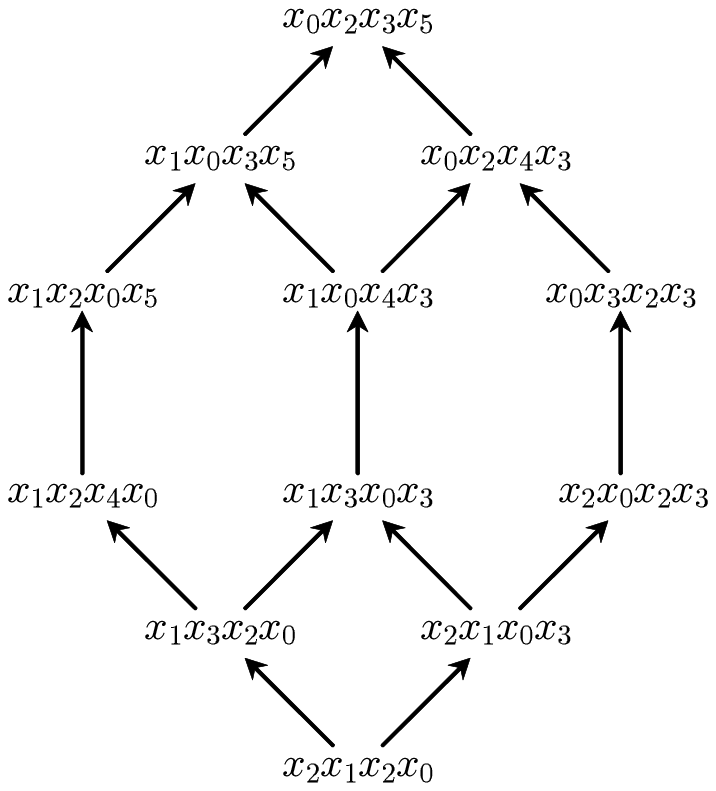}
\end{center}
\end{example}

It seems that each word graph has a unique terminal vertex, namely the normal form, and a unique
initial vertex, which we shall call the \emph{anti-normal form}.

\begin{definition}A word:
\begin{equation*}
x_{i_n}\cdots x_{i_1}
\end{equation*}
is in \emph{anti-normal form} if $i_{k+1} \geq i_k -1$ for all $k$.
\end{definition}

The anti-normal form corresponds to building the carets of a forest diagram from
\emph{left to right}.

Guba and Sapir used anti-normal forms in \cite{GuSa2} to prove that $F$ has a subexponential Dehn function.
(Guba \cite{Guba2} has since shown that the Dehn function of $F$ is quadratic.)  In section 4.1, we will show that the
anti-normal form describes a minimum-length word for any positive element with respect
to the $\{x_0,x_1\}$ generating set.

\begin{proposition}
Let $f$ be any positive element.  Then any terminal vertex in the word graph for $f$ is
in normal form, and any initial vertex is in anti-normal form.
\end{proposition}
\begin{proof}Consider a word $x_{i_n}\ldots x_{i_1}$ for $f$.  This word will have
an outgoing edge if and only if $i_{k+1} > i_k$ for some $k$ (so that we can apply
the move $x_{i_{k+1}} x_{i_k} \rightarrow x_{i_k} x_{i_{k+1}+1}$).  Similarly,
this word will have an incoming edge if and only if $i_{k+1} < i_k -1$ for some $k$
(so that we can apply the inverse move $x_{i_{k+1}} x_{i_k} \rightarrow x_{i_k-1} x_{i_{k+1}}$).
\end{proof}

This proposition gives a nice algorithm for putting a word into either normal form
or anti-normal form.  To put a word into normal form, repeatedly apply moves of the type:
\begin{equation*}
x_n x_k \rightarrow x_k x_{n+1}\quad\text{($n>k$)}
\end{equation*}
Similarly, to put a word into anti-normal form, repeatedly apply moves of the type:
\begin{equation*}
x_k x_n \rightarrow x_{n-1} x_k\quad\text{($k < n-1$)}
\end{equation*}
See example 4.1.9.

    We wish to show that the anti-normal form for an element is unique.  The idea, of course,
is that anti-normal form corresponds to the unique way of constructing the carets of a forest
diagram from left to right.  We shall now establish some notation that makes this idea
very precise.

    Let $\ll$ denote the linear order on the carets
of a forest diagram induced by the order of the spaces that the carets cover:
\begin{center}
\includegraphics{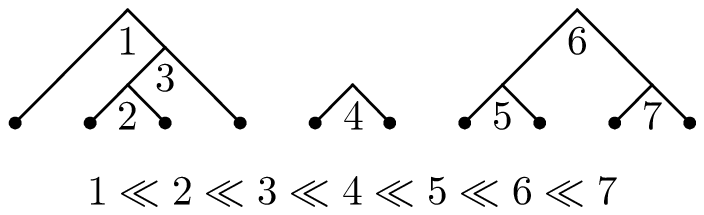}
\end{center}
and let $\prec$ denote the partial order defined by:
\begin{equation*}
c_1 \prec c_2 \quad \Leftrightarrow \quad c_1\text{ is a descendant of }c_2
\end{equation*}
(In the above forest, $2\prec 3\prec 1$ and $5,7\prec 6$.)  Define relations
$<_{\rm{N}}$ and $<_{\rm{AN}}$ as follows:
\begin{equation*}
\begin{array}{lll}
c_1 <_{\rm{N}} c_2 & \Leftrightarrow & c_1 \prec c_2 \;\;\text{or} \;\;
(c_1 \nsucc c_2\;\;\text{and}\;\;c_1 \gg c_2) \\
c_1 <_{\rm{AN}} c_2 & \Leftrightarrow & c_1 \prec c_2 \;\;\text{or} \;\;
(c_1 \nsucc c_2\;\;\text{and}\;\;c_1 \ll c_2)
\end{array}
\end{equation*}
It is not hard to check that $<_{\rm{N}}$ and $<_{\rm{AN}}$ are linear orders.  Furthermore:

\begin{proposition}
Let $x_{i_n} \cdots x_{i_1}$ be a word for a positive element $f$, and let $c_k$ denote the
caret in the reduced forest diagram for $f$ built by $x_{i_k}$.  Then:
\begin{enumerate}
\item The given word is in normal form if and only if $c_k <_{\rm{N}} c_{k+1}$ for each $k$.
\item The given word is in anti-normal form if and only if $c_k <_{\rm{AN}} c_{k+1}$ for each $k$.
\end{enumerate}
\end{proposition}
\begin{proof}
Observe that $c_k \prec c_{k+1}$ if and only if $i_{k+1} = i_k$ or $i_{k+1} = i_k - 1$. Further,
$c_k \ll c_{k+1}$ if and only if $i_k \leq i_{k+1}$.  Therefore:
\begin{equation*}
\begin{array}{lll}
c_k <_{\rm{N}} c_{k+1} & \Leftrightarrow & i_{k+1} \leq i_k \\
c_k <_{\rm{AN}} c_{k+1} & \Leftrightarrow & i_{k+1} \geq i_k - 1
\end{array}\tag*{\qedhere}
\end{equation*}
\end{proof}

\begin{corollary}
Every positive element has a unique normal form and a unique anti-normal form.
\end{corollary}

Next we would like to explain the structure of the word graph.  Given a finite
set $S$, the \emph{order graph} $\Gamma(S)$ on $S$ is the graph whose vertices are linear orders
on the elements of $S$, and whose edges correspond to transpositions of adjacent
elements.  For example, $\Gamma(\{1,2,3\})$ is:
\begin{center}
\includegraphics{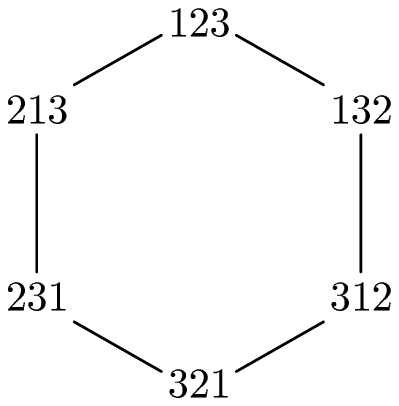}
\end{center}
The order graph of a set with $n$ elements is isomorphic to the Cayley graph of
the symmetric group $\Sigma_n$ with generating set $\{(1\;\;2),(2\;\;3),\ldots,(n-1\;\;n)\}$.
Recall that $\Sigma_n$ has presentation:
\begin{equation*}
\langle t_1,\ldots t_{n-1} \mid t_i^2=1, \; t_i t_{i+1} t_i = t_{i+1} t_i t_{i+1}, \; t_i t_j=t_j t_i
\text{ for $i-j>2$}\rangle
\end{equation*}
where $t_i = (i\;\;i+1)$.  Therefore, any order graph is the one-skeleton of a cell complex
whose two-skeleton is a union of hexagons and squares.

If $f$ is a positive element, the word graph for $f$ is a subgraph of $\Gamma(C)$, where
$C$ is the set of carets in the reduced forest diagram for $f$.  In particular, any word
corresponds to a certain order of building the carets, and any move of the form:
\begin{equation*}
x_n x_k \rightarrow x_k x_{n+1}\quad\text{($n>k$)}
\end{equation*}
corresponds to transposing the order of two adjacent carets.

For the following theorem, recall that a set of vertices $S$ in a graph is \emph{convex} if, whenever
$v,w\in S$, any vertex that appears on any geodesic from $v$ to $w$ is also in $S$.  The \emph{convex hull}
of $S$ is the intersection of all convex sets containing $S$.

\begin{theorem}Let $f$ be a positive element, and let $C$ be the set of carets in the
reduced forest diagram for $f$.  Then the word graph for $f$ is the convex hull in
$\Gamma(C)$ of the normal and anti-normal forms for $f$.
\end{theorem}
\begin{proof}
    Given a set $S$, a \emph{half-space} in $\Gamma(S)$ is a set of the form:
\begin{equation*}
\{s_1<s_2\} = \{\text{linear orders $<$ on $S$} \mid s_1<s_2\}
\end{equation*}
where $s$ and $t$ are fixed elements of $S$. It is a well-known fact that any
convex subset of $\Gamma(S)$ is an intersection of half-spaces.  (See \cite{Bro2}, section A.7.)

    The vertices in the word graph are precisely the linear orders on $C$ that are
extensions of the partial order $\prec$.  In particular, the word graph is precisely
the intersection of all half-spaces $\{c_1<c_2\}$ such that $c_1 \prec c_2$.  However,
for any carets $c_1,c_2\in C$:
\begin{equation*}
c_1 \prec c_2 \quad \Leftrightarrow \quad c_1 <_{\rm{N}} c_2\;\;\text{and}\;\;
c_1 <_{\rm{AN}} c_2
\end{equation*}
Therefore, a half-space contains the word graph if and only if the half-space contains
both the normal form $<_{\rm{N}}$ and the anti-normal form $<_{\rm{AN}}$.\end{proof}

\chapter{Two-Way Forest Diagrams}

In this chapter, we use an action of $F$ on the real line to construct \emph{two-way
forest diagrams}.  These forest diagrams interact very nicely with the
finite generating set $\{x_0,x_1\}$, and are therefore particularly well-suited for
studying the geometry of $F$.

We will generally refer to two-way forest diagrams simply as \emph{forest diagrams}.
We use the terminology ``one-way forest diagrams'' and ``two-way forest diagrams''
only when there is some ambiguity.

The material in this chapter represents joint work with my thesis advisor, Kenneth Brown.
It was originally published in \cite{BeBr}.

\section{The Group $\PL_2(\mathbb{R})$}

Let $\PL_2(\mathbb{R})$ be the group of all piecewise-linear, orientation-preserving
self-homeo\-morphisms $f$ of $\mathbb{R}$ satisfying the following conditions:
\begin{enumerate}
\item Each linear segment of $f$ has slope a power of $2$.
\item $f$ has only finitely many breakpoints, each of which has dyadic rational coordinates.
\item The leftmost linear segment of $f$ is of the form:
\begin{equation*}
f(t)=t+m
\end{equation*}
and the rightmost segment is of the form:
\begin{equation*}
f(t)=t+n
\end{equation*}
for some integers $m,n$.
\end{enumerate}

\begin{proposition}
$\PL_2(\mathbb{R})$ is isomorphic with $F$.
\end{proposition}

\begin{proof}
Let $\psi\colon\mathbb{R}\rightarrow(0,1)$ be the
piecewise-linear homeomorphism that maps the intervals:
\begin{center}
\includegraphics{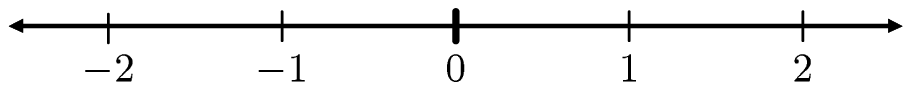}
\end{center}
linearly onto the intervals:
\begin{center}
\includegraphics{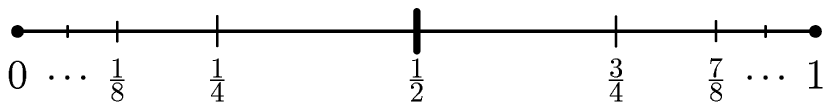}
\end{center}
Then $f\mapsto\psi f\psi^{-1}$ is the desired isomorphism $F\rightarrow\PL_2(\mathbb{R})$.
\end{proof}

Under this isomorphism, the generators $\{x_0,x_1\}$ of $F$ map to the functions:
\begin{equation*}
x_0(t)=t+1
\end{equation*}
and:
\begin{center}
$x_1(t)=
\begin{cases}
t&t\leq 0\\
2t&0\leq t\leq 1\\
t+1&t\geq 1
\end{cases}
\qquad\qquad\quad$\parbox{1.5in}{\includegraphics{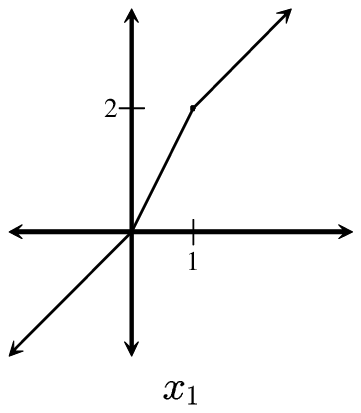}}
\end{center}

\

\section{Forest Diagrams for Elements of $\PL_2(\mathbb{R})$}

We think of the real line as being pre-subdivided as follows:
\begin{center}
\includegraphics{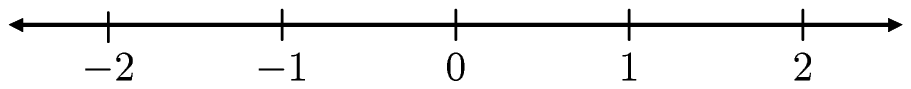}
\end{center}
A \emph{dyadic subdivision} of $\mathbb{R}$ is a subdivision obtained by cutting
finitely many of these intervals in half, and then cutting finitely many of the resulting
intervals in half, etc.

\begin{proposition}
Let $f\in\PL_2(\mathbb{R})$. Then there exist dyadic subdivisions
$\mathcal{D},\mathcal{R}$ of $\mathbb{R}$ such that $f$ maps each interval of
$\mathcal{D}$ linearly onto an interval of $\mathcal{R}.\quad\qedsymbol$
\end{proposition}

A \emph{two-way binary forest} is a sequence $(\ldots,T_{-1},T_0,T_1,\ldots)$ of
finite binary trees. We depict such a forest as a line of binary trees
together with a pointer at~$T_0$:
\begin{center}
\includegraphics{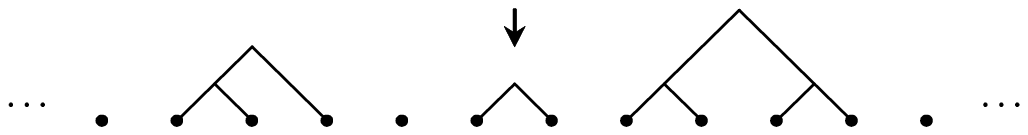}
\end{center}
Every bounded, two-way binary forest corresponds to some dyadic subdivision of the real line.
Therefore, any $f\in\PL_2(\mathbb{R})$ can be represented by a pair of bounded binary forests,
together with an order-preserving bijection of their leaves. This is called a
\emph{two-way forest diagram} for $f$.

\begin{example}
Here are the two-way forest diagrams for $x_0$ and $x_1$:
\begin{center}
\includegraphics{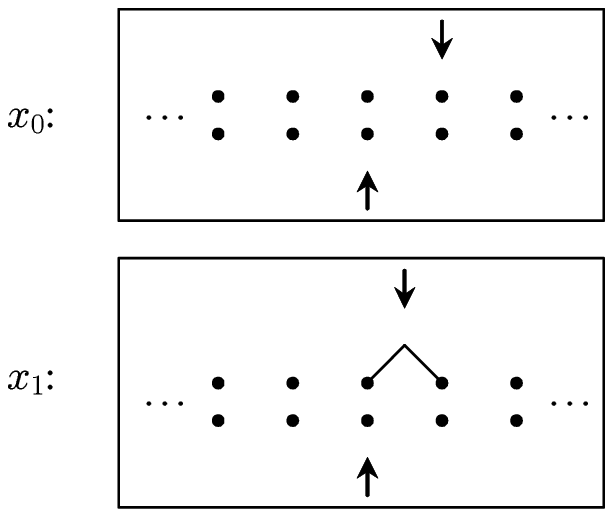}
\end{center}
\end{example}

\begin{proposition}
Every element of $\PL_2(\mathbb{R})$ has a unique reduced two-way forest diagram.\quad\qedsymbol
\end{proposition}

\begin{remark}
From this point forward, we will only draw the \emph{support} of the
two-way forest diagram (i.e. the minimum interval containing both pointers
and all nontrivial trees), and we will omit the ``$\cdots$'' indicators.

Also, the term ``forest diagram'' when used alone will always refer to two-way forest
diagrams.
\end{remark}

\begin{remark}
It is fairly easy to translate between tree diagrams, one-way forest diagrams, and
two-way forest diagrams.  Given a tree diagram:
\begin{center}
\includegraphics{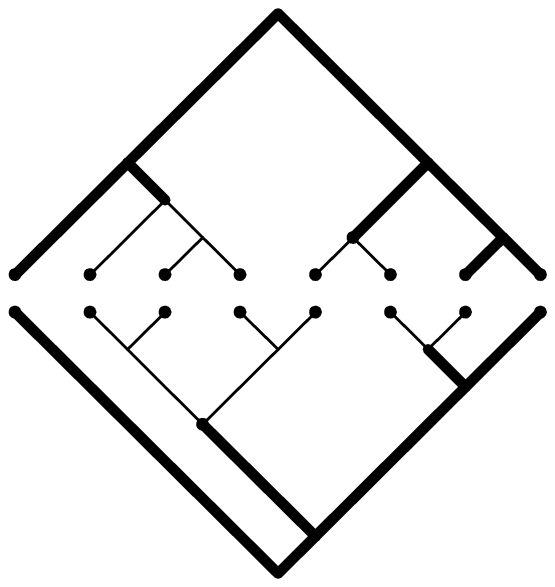}
\end{center}
we simply remove the outer layer of each tree to get the corresponding two-way forest diagram:
\begin{center}
\includegraphics{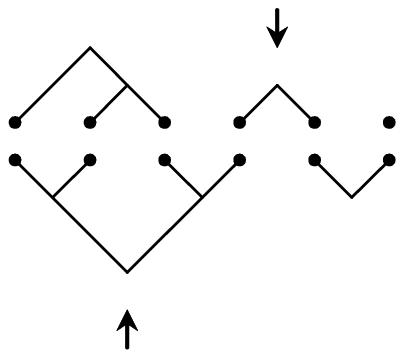}
\end{center}
The pointers of the forest diagram point to the first trees hanging to the right
of the roots in the original tree diagram.

Similarly, given a one-way forest diagram:
\begin{center}
\includegraphics{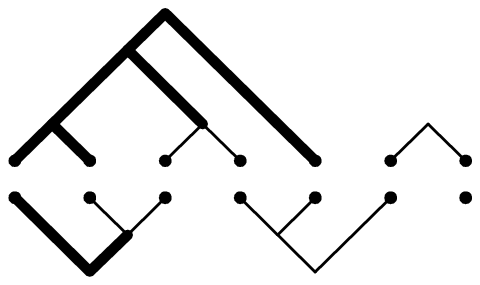}
\end{center}
simply remove the left stalk of $0$-tree on the top and bottom to get the corresponding
two-way forest diagram:
\begin{center}
\includegraphics{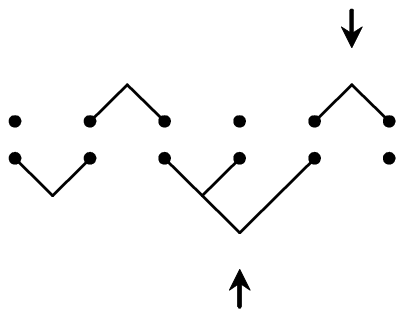}
\end{center}
The pointers of the two-way forest diagram point to the $1$-trees of the original one-way
forest diagram.
\end{remark}

\

\section{The Action of $\{x_0,x_1\}$}

Just as one-way forest diagrams interact well with the infinite generating set for~$F$,
two-way forest diagrams interact well with the $\{x_0,x_1\}$ generating set:

\begin{proposition}
Let $\mathfrak{f}$ be a forest diagram for some $f\in F$. Then:
\begin{enumerate}
\item A forest diagram for $x_0f$ can be obtained by moving the
top pointer of $\mathfrak{f}$ one tree to the right.
\item A forest diagram for $x_1f$ can be obtained by attaching a
caret to the roots of the $0$-tree and $1$-tree in the top forest of
$\mathfrak{f}$. Afterwards, the top pointer points to the new, combined tree.\quad\qedsymbol
\end{enumerate}
\end{proposition}

If $\mathfrak{f}$ is reduced, then the given forest diagram for $x_0f$ will
always be reduced. The forest diagram given for $x_1f$ will not be reduced,
however, if the caret that was created opposes a caret from the bottom tree.
In this case, left-multiplication by $x_1$ effectively ``cancels'' the
bottom caret.

\begin{example}
Let $f\in F$ have forest diagram:
\begin{center}
\includegraphics{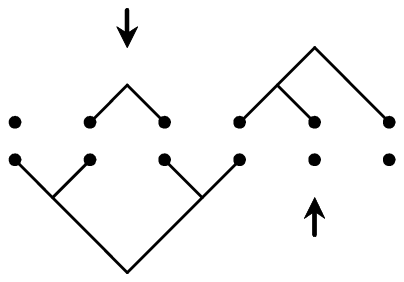}
\end{center}
Then $x_0f$ has forest diagram:
\begin{center}
\includegraphics{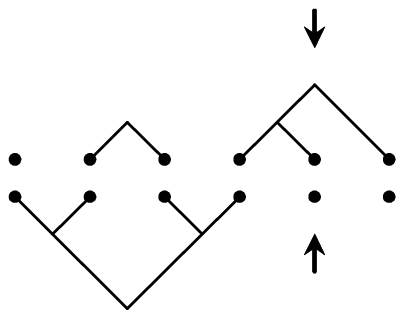}
\end{center}
and $x_1f$ has forest diagram:
\begin{center}
\includegraphics{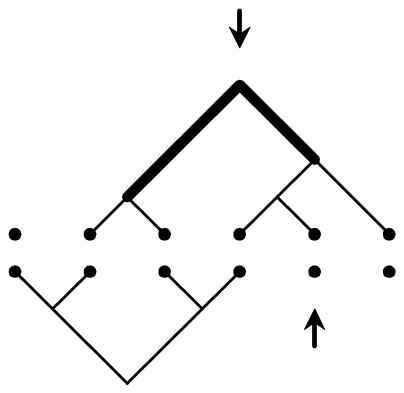}
\end{center}
\end{example}

\begin{example}
Let $f\in F$ have forest diagram:
\begin{center}
\includegraphics{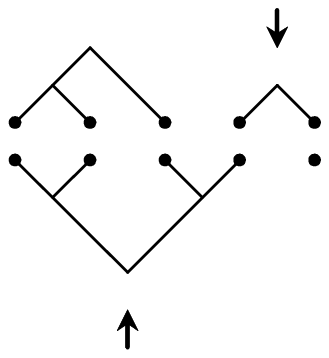}
\end{center}
Then $x_0f$ has forest diagram:
\begin{center}
\includegraphics{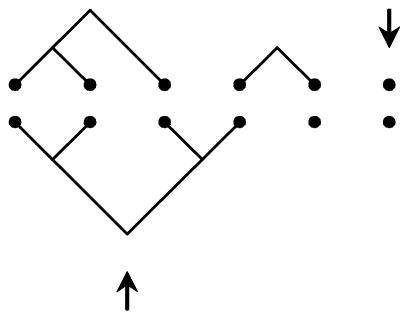}
\end{center}
and $x_1f$ has forest diagram:
\begin{center}
\includegraphics{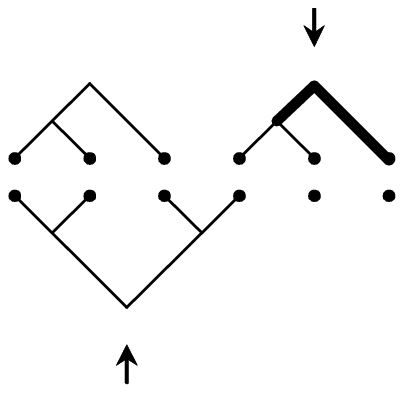}
\end{center}
Note that the forest diagrams for $x_0f$ and $x_1f$ both have larger support
than the forest diagram for $f$.
\end{example}

\begin{example}
Let $f\in F$ have forest diagram:
\begin{center}
\includegraphics{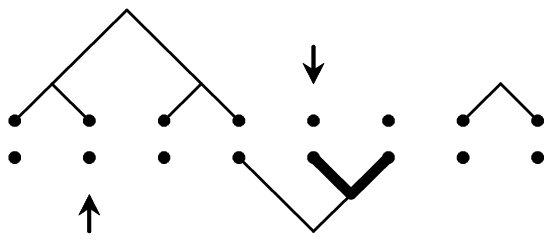}
\end{center}
Then $x_1f$ has forest diagram:
\begin{center}
\includegraphics{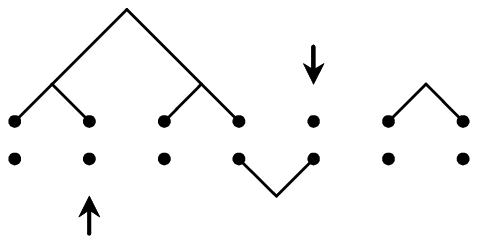}
\end{center}
Note that left-multiplication by $x_1$ cancelled the highlighted bottom caret.
\end{example}

\begin{proposition}
Let $\mathfrak{f}$ be a forest diagram for some $f\in F$. Then:
\begin{enumerate}
\item A forest diagram for $x_0^{-1}f$ can be obtained by moving the
top pointer of $\mathfrak{f}$ one tree to the left.
\item A forest diagram for $x_1^{-1}f$ can be obtained by ``dropping a
negative caret'' at the current position of the top pointer. If the current tree
is nontrivial, the negative caret cancels with the top caret of the current tree,
and the pointer moves to the resulting left child. If the current tree is trivial,
the negative caret ``falls through'' to the bottom forest, attaching to the specified leaf.
\quad\qedsymbol
\end{enumerate}
\end{proposition}

\begin{example}
\quad Let $f$ and $g$ be the elements of $F$ with forest diagrams:
\begin{center}
\includegraphics{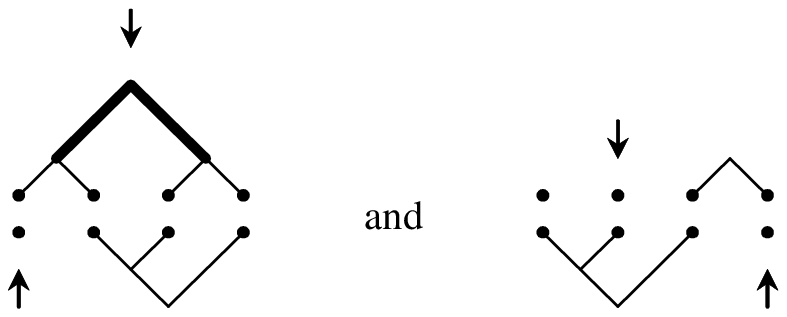}
\end{center}
Then $x_1^{-1}f$ and $x_1^{-1}g$ have forest diagrams:
\begin{center}
\includegraphics{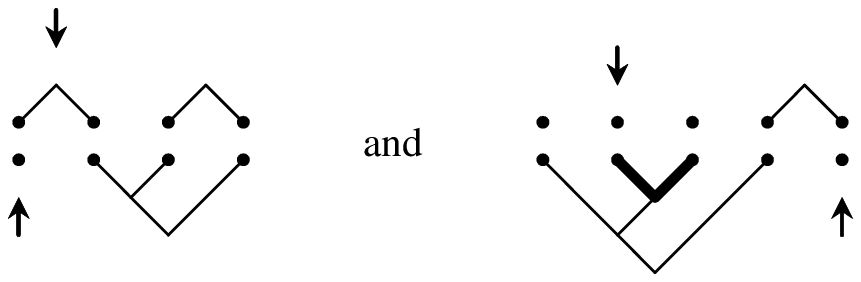}
\end{center}
In the first case, the $x_1^{-1}$ simply removed a caret from the top tree.
In the second case, there was no caret on top to remove, so a new caret was
attached to the leaf on the bottom. Note that this creates a new column
immediately to the right of the pointer.
\end{example}

\section{Normal Forms and Positive Elements}

The other generators of $F$ act on two-way forest diagrams in the following way:

\begin{proposition}
Let $\mathfrak{f}$ be the forest diagram for some $f\in F$, and let $n>1$.  Then a forest
diagram for $x_nf$ can be obtained by attaching a caret to the roots of $T_{n-1}$ and $T_n$
in the top forest of $\mathfrak{f}$.\quad\qedsymbol
\end{proposition}

Using this proposition, it is relatively easy to find the normal form from the forest diagram,
using a method similar to that given in section 2.4.

\begin{example}
Suppose $f\in F$ has forest diagram:
\begin{center}
\includegraphics{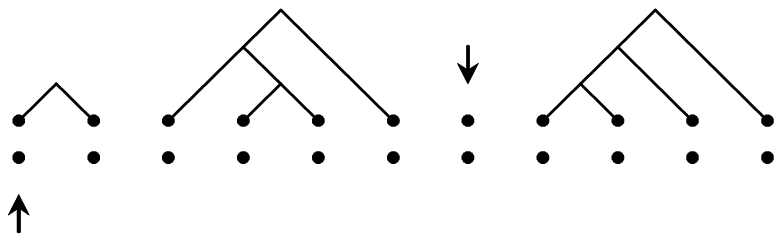}
\end{center}
Then:
\begin{equation*}
f = x_0^2 x_1 x_3^2 x_4 x_8^3
\end{equation*}
Since the top pointer of $f$ is two trees from the left, the normal form of $f$ has an $x_0^2$.
The powers of the other generators are determined by the number of carets built upon the
corresponding leaf.  Note that the carets are constructed from right to left.
\end{example}

\begin{example}
The element:
\begin{equation*}
x_0^3 x_2 x_5^2 x_7 x_6^{-1} x_5^{-1} x_1^{-2} x_0^{-1}
\end{equation*}
has forest diagram:
\begin{center}
\includegraphics{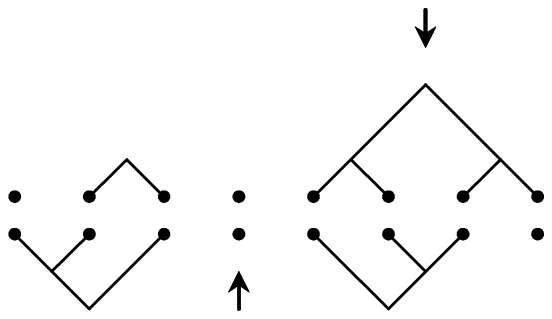}
\end{center}
\end{example}

Proposition 3.4.1 also yields a characterization of positive elements.

\begin{corollary}
Let $f\in F$, and let $\mathfrak{f}$ be its reduced forest diagram.  Then $f$ is positive
if and only if:
\begin{enumerate}
\item The entire bottom forest of $\mathfrak{f}$ is trivial, and
\item The bottom pointer is at the left end of the support of $\mathfrak{f}$.
\end{enumerate}
\end{corollary}

So a typical positive element looks like:
\begin{center}
\includegraphics{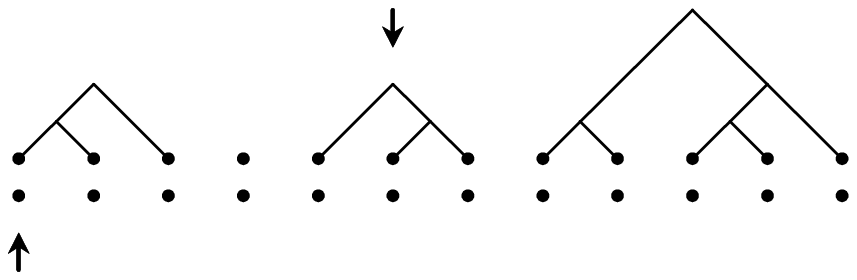}
\end{center}

An element of $f$ is \emph{right-sided} if \emph{both} pointers are at the left end of the
support. The monoid of right-sided elements is generated by
$\{x_1,x_2,\ldots ,x_1^{-1},x_2^{-1},\ldots\}$.

An element which is both positive and right-sided is called \emph{strongly positive}.
The monoid of strongly positive elements is generated by $\{x_1,x_2\ldots\}$.

\

\chapter{Lengths in $F$}

In this chapter, we derive a formula for the lengths of elements of $F$ with
respect to the $\{x_0,x_1\}$ generating set. This formula uses the two-way forest
diagrams introduced in chapter 3.

Lengths in $F$ were first studied by S.~B.~Fordham in his 1995 thesis (recently published,
see \cite{Ford}). Fordham gave a formula for the length of an element of $F$ based on its
tree diagram.  Our length formula can be viewed as a simplification of Fordham's work.

V. Guba has recently obtained another length formula for $F$ using the ``diagrams'' of Guba
and Sapir.  See \cite{Guba} for details.

The material in this chapter represents joint work with my thesis advisor, Kenneth Brown.
It was originally published in \cite{BeBr}.

\

\section{Lengths of Strongly Positive Elements}

We shall begin by investigating the lengths of strongly positive elements.
The goal is to develop some intuition for lengths before the statement of
the general length formula in section 4.2.

Recall that an element is \emph{strongly positive} if it lies in the submonoid
generated by $\{x_1,x_2,\ldots\}$.  Equivalently, $f$ is strongly positive if and only
if the entire bottom forest of $f$ is trivial and both pointers are at the left edge
of the support:
\begin{center}
\includegraphics{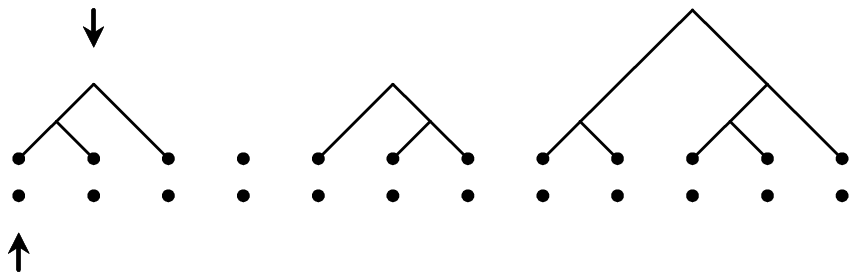}
\end{center}

Logically, the results of this section depend on the general length formula.
In particular, we need the following lemma:

\begin{lemma}
Let $f\in F$ be strongly positive. Then there exists a minimum-length
word for $f$ with no appearances of $x_1^{-1}$.
\end{lemma}

This lemma is intuitively obvious: there should be no reason to ever create
bottom carets, or to delete top carets, when constructing a strongly positive element.
Unfortunately, it would be rather tricky to supply a proof of this
fact.  Instead we refer the reader to corollary 4.3.8, from which the lemma follows
immediately.

From this lemma, we see that any strongly positive element $f\in F$ has a
minimum-length word of the form:
\begin{equation*}
x_0^{a_n}x_1\cdots x_0^{a_1}x_1x_0^{a_0}
\end{equation*}
where $a_0,\ldots,a_n\in\mathbb{Z}$. Since $f$ is strongly positive, we have:
\begin{equation*}
a_0+\cdots+a_n=0
\end{equation*}
and
\begin{equation*}
a_0+\cdots+a_i\geq 0
\parbox{0in}{\parbox{2in}{\qquad(for $i=0,\ldots,n-1$)}}
\end{equation*}
Such words can be represented by words in $\{x_1,x_2,\ldots\}$ via
the identifications $x_n=x_0^{1-n}x_1x_0^{n-1}$. For example, the word:
\begin{equation*}
x_0^{-5}\,x_1\,x_0^{-2}\,x_1\,x_0^4\,x_1\,x_0^{-3}\,x_1\,x_0^6
\end{equation*}
can be represented by:
\begin{equation*}
x_6\,x_8\,x_4\,x_7
\end{equation*}
More generally:

\begin{notation}
We will use the word:
\begin{equation*}
x_{i_n}\cdots x_{i_2}x_{i_1}
\end{equation*}
in $\{x_1,x_2,\ldots\}$ to represent the word:
\begin{equation*}
x_0^{1-i_n}\,x_1\,\cdots\,x_0^{i_3-i_2}\,x_1\,x_0^{i_2-i_1}\,x_1\,x_0^{i_1-1}
\end{equation*}
in $\{x_0,x_1\}$.
\end{notation}

Note then that $x_{i_n}\cdots x_{i_2}x_{i_1}$ represents a word with length:
\begin{equation*}
\left(|1-i_n|+\cdots+|i_3-i_2|+|i_2-i_1|+|i_1-1|\right)+n
\end{equation*}
We now proceed to some examples, from which we will derive a general theorem.

\begin{example}
Let $f\in F$ be the element with forest diagram:
\begin{center}
\includegraphics{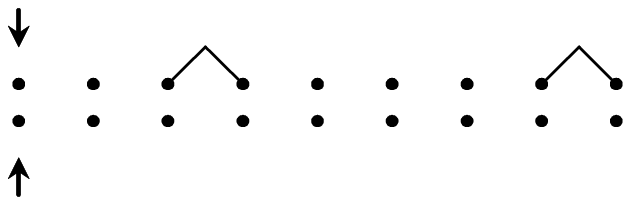}
\end{center}
There are only two candidate minimum-length words for $f$: $x_3x_8$ and
$x_7x_3$. Their lengths are:
\begin{equation*}
\begin{aligned}[t]
&&\quad&\left(2+5+7\right)+2=16&\quad&\quad\text{for the word }x_3x_8\\
&\text{and}\quad&\quad&\left(6+4+2\right)+2=14&\quad&\quad\text{for the word
}x_7x_3.
\end{aligned}
\end{equation*}
Let's see if we can explain this. The word $x_3x_8=x_0^{-2}x_1x_0^{-5}x_1x_0^7$
corresponds to the following construction of $f$:
\begin{enumerate}
\item Starting at the identity, move right seven times
and construct the right caret.
\item Next move left five times, and construct the left caret.
\item Finally, move left twice to position of the bottom pointer.
\end{enumerate}
This word makes a total of fourteen moves, crossing twice over each of seven spaces:
\begin{center}
\includegraphics{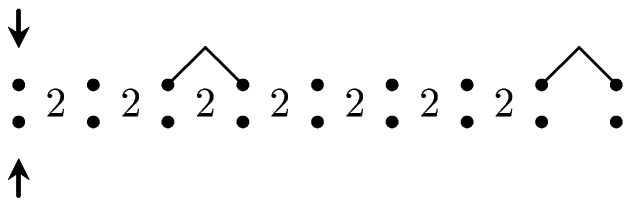}
\end{center}

On the other hand, the word $x_7x_3=x_0^{-6}x_1x_0^4x_1x_0^2$ corresponds to
the following construction:
\begin{enumerate}
\item Starting at the identity, move right twice and construct the
left caret.
\item Next move right four more times, and construct the right caret.
\item Finally, move left six times to the position of the bottom pointer.
\end{enumerate}
This word makes only twelve moves:
\begin{center}
\includegraphics{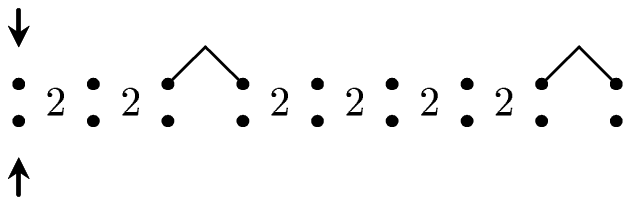}
\end{center}
In particular, this word never moves across the space under the left caret.
It avoids this by \emph{building the left caret early}.  Once the left caret
is built, the word can simply pass over the space under the left caret without
spending time to move across it.
\end{example}

\begin{terminology}
We call a space in a forest \emph{interior} if it lies under a tree (or
over a tree, if the forest is upside-down) and \emph{exterior} if it lies
between two trees.
\end{terminology}

\begin{example}
Let $f\in F$ be the element with forest diagram:
\begin{center}
\includegraphics{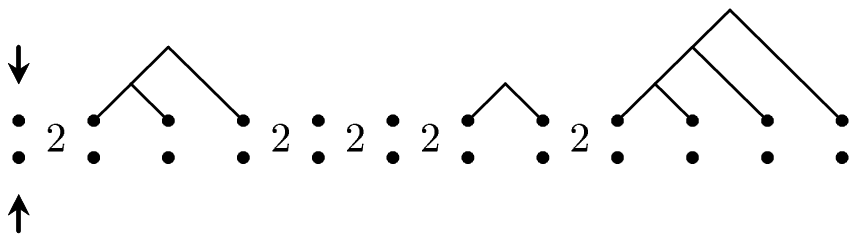}
\end{center}
Clearly, each of the five exterior spaces in the support of $f$ must be
crossed twice during construction. Furthermore, it is possible to avoid
crossing any of the interior spaces by \emph{constructing carets from left to
right}. In particular:
\begin{equation*}
x_6^3\,x_5\,x_2^2
\end{equation*}
is a minimum-length word for $f$. Therefore, $f$ has length:
\begin{equation*}
(5+1+3+1)+6=16
\end{equation*}
\end{example}

It is not always possible to avoid crossing all the interior spaces:

\begin{example}
Let $f\in F$ be the element with forest diagram:
\begin{center}
\includegraphics{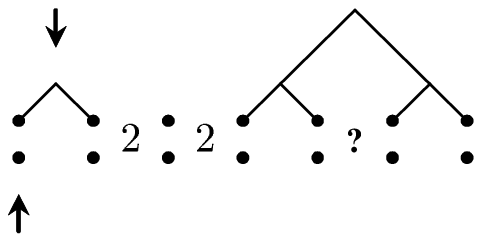}
\end{center}
Clearly, each of the two exterior spaces in the support of $f$ must be
crossed twice during construction. However, the space marked (?) must also
be crossed twice, since we must create the caret immediately to its right
before we can create the caret above it.

It turns out that these are the only spaces which must be crossed. For
example, the word:
\begin{equation*}
x_3\,x_4\,x_3\,x_1
\end{equation*}
crosses only these spaces.  Therefore, $f$ has length:
\begin{equation*}
(2+1+1+2+0)+4=10
\end{equation*}
\end{example}

Recall that a word:
\begin{equation*}
x_{i_n}\cdots x_{i_2}x_{i_1}
\end{equation*}
is in \emph{anti-normal form} if $i_{k+1} \geq i_k-1$ for all $k$.  While normal form
corresponds to building carets from right to left, anti-normal
form corresponds to building carets from \emph{left to right} (i.e. constructing the leftmost
possible caret at each stage).

As we have seen, the anti-normal form for a strongly positive element has minimum length, since
it crosses only those spaces that must be crossed.  We can get an explicit length
formula by counting these spaces:

\begin{theorem}
Let $f\in F$ be strongly positive. Then the length of $f$ is:
\begin{equation*}
2\,n(f)+c(f)
\end{equation*}
where:
\begin{enumerate}
\item $n\left(f\right)$ is the number of spaces in the support of $f$
that are either exterior or lie immediately to the left of some caret, and
\item $c\left(f\right)$ is the number of carets of $f$.\quad\qedsymbol
\end{enumerate}
\end{theorem}

\begin{example}
\quad Let $f\in F$ be the element with forest diagram:
\begin{center}
\includegraphics{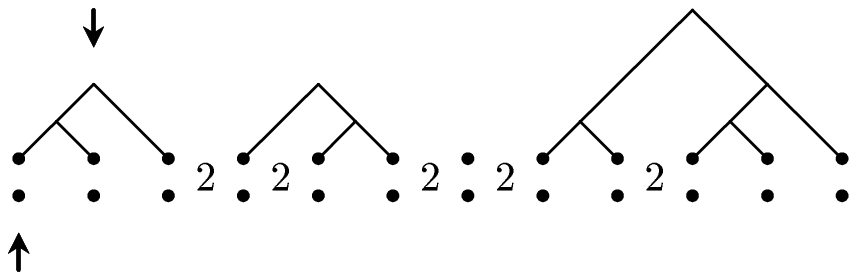}
\end{center}
Then $c(f)=8$ and $n(f)=5$, so $f$ has length $18$. The anti-normal form for $f$ is:
\begin{equation*}
x_4x_5^2x_4x_2x_3x_1^2
\end{equation*}
Therefore, a minimum-length $\{x_0,x_1\}$-word for $f$ is:
\begin{equation*}
x_0^{-3}x_1x_0^{-1}x_1^2x_0x_1x_0^2x_1x_0^{-1}x_1x_0^2x_1^2
\end{equation*}
\end{example}

Because the anti-normal form is the unique terminal vertex in the word graph (see section
2.5), we can put any positive element into anti-normal form by repeatedly applying the
operations:
\begin{equation*}
x_k x_n \; \rightarrow \; x_{n-1} x_k\parbox{0in}{\parbox{2in}{\qquad($k < n-1$)}}
\end{equation*}
This gives us an entirely algebraic algorithm for finding the length of an element.

\begin{example}
Let's find the length of the element:
\begin{equation*}
x_1\,x_3^3\,x_6\,x_7\,x_{10}
\end{equation*}
We put the word into anti-normal form:
\begin{equation*}
\begin{aligned}[b]
x_1\,&x_3^3\,x_6\,x_7\,x_{10}\\
= x_4\,&x_1\,x_3^3\,x_6\,x_7\\
= x_4\,&x_2^3\,x_5\,x_6\,x_1\\
= x_4\,&x_2\,x_3\,x_4\,x_2^2\,x_1
\end{aligned}
\parbox[b]{0in}{\parbox[b]{2in}{\qquad
$\begin{gathered}[b]
\\
\text{(}x_{10}\text{ moved left)}\\
\text{(}x_1\text{ moved right)}\\
\text{(}x_2^2\text{ moved right)}
\end{gathered}
$}}
\end{equation*}
Therefore, the length is:
\begin{equation*}
\left(3+2+1+1+2+1+0\right)+7=17
\end{equation*}
\end{example}

\

\section{The Length Formula}

We now give the length formula for a general element of $F$. Afterwards, we
will give several examples to illustrate intuitively why the formula works.
We defer the proof to section 4.3.

Let $f\in F$, and let $\mathfrak{f}$ be its reduced forest diagram. We label
the spaces of each forest of $\mathfrak{f}$ as follows:
\begin{enumerate}
\item Label a space $\LL$ (for \emph{left}) if it exterior and to the
left of the pointer.
\item Label a space $\NN$ (for \emph{necessary}) if it lies immediately
to the left of some caret (and is not already labeled $\LL$).
\item Label a space $\RR$ (for \emph{right}) if it exterior and
to the right of the pointer (and not already labeled $\NN$).
\item Label a space $\II$ (for \emph{interior}) if it interior
(and not already labeled $\NN$).
\end{enumerate}
We assign a \emph{weight} to each space in the support of $\mathfrak{f}$
according to its labels:

\begin{center}
\begin{tabular}{c}
\\
$\text{top}$\\
$\text{label}$
\end{tabular}%
\begin{tabular}{|c|cccc|}
\multicolumn{5}{c}{$\text{bottom label}$}\\
\hline\rule{0ex}{2.5ex}
&$\;\LL\;$&$\;\NN\;$&$\;\RR\;$&$\;\II\;$\\
\hline
$\;\LL\;$&$2$&$1$&$1$&$1$\rule{0ex}{2.5ex}\\
$\NN$&$1$&$2$&$2$&$2$\\
$\RR$&$1$&$2$&$2$&$0$\\
$\II$&$1$&$2$&$0$&$0$\\
\hline
\end{tabular}
\end{center}

\begin{example}
Here are the labels and weights for a typical forest diagram:
\begin{center}
\includegraphics{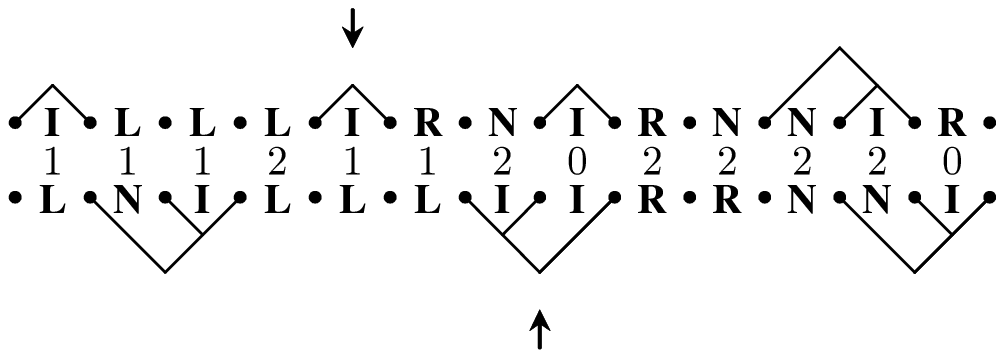}
\end{center}
\end{example}

\begin{theorem}[The Length Formula]
Let $f\in F$, and let $\mathfrak{f}$ be its reduced forest diagram. Then
the $\{x_0,x_1\}$-length of $f$ is:
\begin{equation*}
\ell(f)=\ell_0(f)+\ell_1(f)
\end{equation*}
where:
\begin{enumerate}
\item $\ell_0(f)$ is the sum of the weights of all spaces in the
support of\text{ }$\mathfrak{f}$, and
\item $\ell_1(f)$ is the total number of carets in $\mathfrak{f}$.
\end{enumerate}
\end{theorem}

\begin{remark}
Intuitively, the weight of a space is just the number of times it must
be crossed during the construction of $f$. Hence, there ought to exist
a minimum-length word for $f$ with $\ell_0(f)$ appearances of $x_0$ or $x_0^{-1}$
and $\ell_1(f)$ appearances of $x_1$ or $x_1^{-1}$. This will be established at the
end of the next section.
\end{remark}

\begin{example}
Let $f\in F$ be the element from example 4.1.8:
\begin{center}
\includegraphics{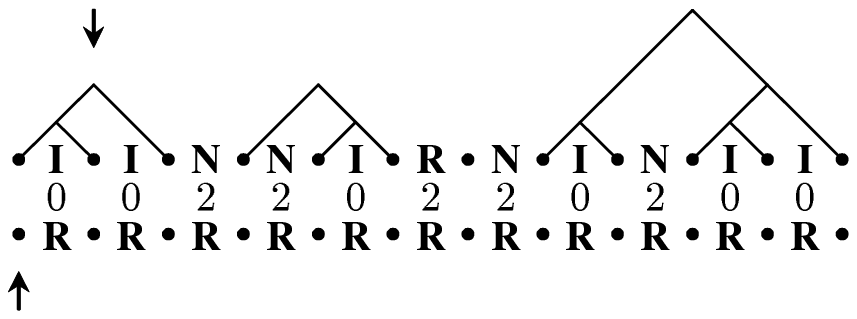}
\end{center}
Then $\ell_0(f)=10$ and $\ell_1(f)=8$, so $f$ has length $18$.
\end{example}

In general, suppose $f\in F$ is strongly positive, and let $\mathfrak{f}$ be
its reduced forest diagram. Then every space of $\mathfrak{f}$ is labeled \pair{\NN}{\RR},
\pair{\RR}{\RR}, or \pair{\II}{\RR}.  Each \pair{\II}{\RR} space has weight $0$, and each
\pair{\NN}{\RR} or \pair{\RR}{\RR} space has weight $2$, so that:
\begin{equation*}
\ell_0(f)=2{\,}n(f)
\end{equation*}
and hence:
\begin{equation*}
\ell_0(f)+\ell_1(f)=2\,n(f)+c(f)
\end{equation*}
Therefore, the length formula of theorem 4.2.2 reduces to theorem 4.1.7
for strongly positive elements.

\begin{example}
Let $f\in F$ be the right-sided element with forest diagram:
\begin{center}
\includegraphics{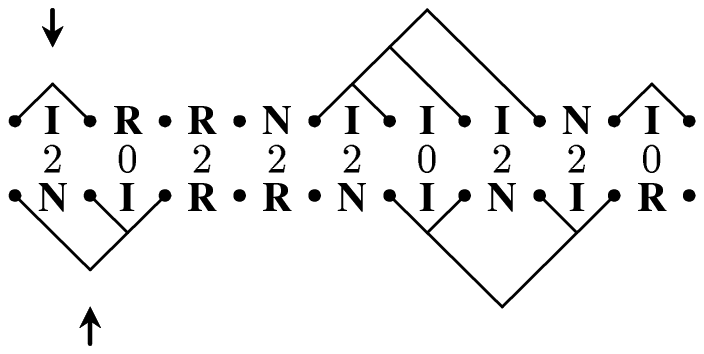}
\end{center}
Then $\ell_0(f)=12$ and $\ell_1(f)=10$, so $f$ has length $22$. One minimum-length word for $f$ is:
\begin{equation*}
x_1 x_0^{-1} x_1^{-1} x_0 x_1^{-1} x_0^{-3} x_1 x_0 x_1^3 x_0^{-1} x_1^{-1} x_0^{-1} x_1^{-1}
x_0 x_1^{-1} x_0^3
\end{equation*}
\end{example}

In general, every space in the forest diagram of a right-sided element is
labeled either $\NN$, $\RR$, or $\II$. The weight table for such spaces is:
\begin{center}
\begin{tabular}{c}
\\
top\\
label
\end{tabular}%
\begin{tabular}{|c|ccc|}
\multicolumn{4}{c}{bottom label}\\
\hline\rule{0ex}{2.5ex}
&$\;\NN\;$&$\;\RR\;$&$\;\II\;$\\
\hline
$\;\NN\;$&$2$&$2$&$2$\rule{0ex}{2.5ex}\\
$\RR$&$2$&$2$&$0$\\
$\II$&$2$&$0$&$0$\\
\hline
\end{tabular}
\end{center}
Observe that a space has weight $2$ if and only if:
\begin{enumerate}
\item It is exterior on both the top and the bottom, or
\item It lies immediately to the left of some caret, on either the top or the bottom.
\end{enumerate}
This can be viewed as a generalization of the length formula for strongly positive elements.
Specifically, if $f$ is right-sided, then:
\begin{equation*}
\ell(f)=2\,n(f)+c(f)
\end{equation*}
where $n(f)$ is the number of spaces satisfying condition (1) or (2), and $c(f)$ is the number
of carets of $f$.

As with strongly positive elements, it is intuitively obvious that this is a
lower bound for the length. Unfortunately, we have not been able to find an
analogue of the ``anti-normal form'' argument to show that it is an upper
bound.

\begin{example}
Let $f\in F$ be the left-sided element with forest diagram:
\begin{center}
\includegraphics{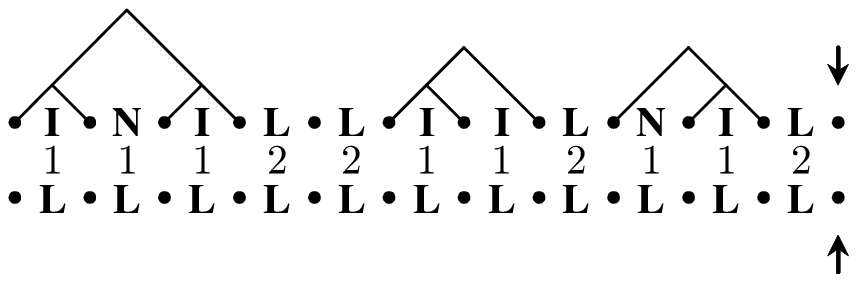}
\end{center}
Then $\ell_0(f)=15$ and $\ell_1(f)=7$, so $f$ has length $22$.

It is interesting to note that every interior space of $f$ has weight $1$:
for trees to the left of the pointer, one cannot avoid crossing interior
spaces at least once. Specifically, each caret is created from its
\emph{left} leaf, and we must move to this leaf somehow.

One minimum-length word for $f$ is
\begin{equation*}
x_0^4 x_1^2 x_0^{-2} x_1 x_0^{-3} x_1^2 x_0^{-3} x_1 x_0^{-1} x_1 x_0^{-2}
\end{equation*}
Note that this word creates carets \emph{right to left}.
\end{example}

\begin{example}
Let $f\in F$ be the element with forest diagram:
\begin{center}
\includegraphics{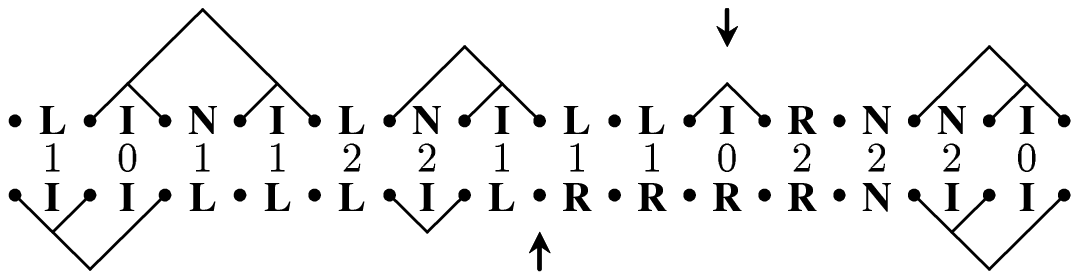}
\end{center}
Then $\ell_0(f)=16$ and $\ell_1(f)=13$, so $f$ has length $29$. One minimum-length word for $f$ is:
\begin{equation*}
x_0^{-2} x_1 x_0^{-1} x_1 x_0 x_1^{-2} x_0^2 x_1 x_0^3 x_1^2 x_0 x_1^{-2} x_0^{-1} x_1
x_0^{-2} x_1 x_0^{-1} x_1 x_0 x_1^{-1} x_0^{-1}
\end{equation*}
This is our first example with \pair{\LL}{\RR} pairs: note that they only need to be crossed once.
Also note how it affects the length to have bottom trees to the left of the pointer. In particular,
observe that the \pair{\NN}{\II} pair to the left of the pointers must crossed twice.
\end{example}

\

\section{The Proof of the Length Formula}

We prove the length formula using the same technique as Fordham \cite{Ford}:
\begin{theorem}
Let $G$ be a group with generating set $S$, and let $\ell\colon G\rightarrow\mathbb{N}$
be a function. Then $\ell$ is the length function for $G$ with respect to $S$ if and only if:
\begin{enumerate}
\item $\ell(e)=0$, where $e$ is the identity of $G$.
\item $\left|\ell(sg)-\ell(g)\right|\leq 1$ for all $g\in G$ and $s\in S$.
\item If $g\in G\setminus\{e\}$, there exists an $s\in S\cup S^{-1}$ such that $\ell(sg)<\ell(g)$.
\end{enumerate}
\end{theorem}
\begin{proof}
Conditions (1) and (2) show that $\ell$ is a lower bound for the
length, and condition (3) shows that $\ell$ is an upper bound for the
length.
\end{proof}

Let $\ell$ denote the function defined on $F$ specified by Theorem 4.2.2.
Clearly $\ell$ satisfies condition (1). To show that $\ell$ satisfies
conditions (2) and (3), we need only gather information about how
left-multiplication by generators affects the function~$\ell$.

\begin{terminology}
If $f\in F$, the \emph{current tree} of $f$ is the tree in forest diagram
indicated by the top pointer. The \emph{right space} of $f$ is the space
immediately to the right of the current tree, and the \emph{left space}
of $f$ is the space immediately to the left of the current tree.
\end{terminology}

Note that, if the top pointer is at the right edge of the support of $f$, then the right
space of $f$ has no label.  Similarly, if the top pointer is at the left edge of the support,
then the left space of $f$ has no label.

\begin{proposition}
If $f\in F$, then $\ell(x_0f)=\ell(f){\;\pm\;}1$. Specifically, $\ell(x_0f)<\ell(f)$
unless one of the following conditions holds:
\begin{enumerate}
\item $x_0f$ has larger support than $f$.
\item The right space of $f$ has bottom label $\LL$, and left-multiplication
by $x_0$ does not remove this space from the support.
\item The right space of $f$ is labeled \pair{\RR}{\II}.
\end{enumerate}
\end{proposition}

\begin{proof}
Clearly $\ell_1(x_0f)=\ell_1(f)$. As for $\ell_0$, note that the only space whose
label changes is the right space of $f$.

\emph{Case 1}:\quad Suppose $x_0f$ has larger support than $f$. Then the right
space of $f$ is unlabeled, and has label \pair{\LL}{\RR} in $x_0f$. Hence
$\ell_0(x_0f)=\ell_0(f)+1$.

\emph{Case 2}:\quad Suppose $x_0f$ has smaller support than $f$. Then the right
space of $f$ has label \pair{\RR}{\LL}, but becomes unlabeled in $x_0f$. Hence
$\ell_0(x_0f)=\ell_0(f)-1$.

\emph{Case 3}:\quad Suppose $x_0f$ has the same support as $f$. Then the right
space of $f$ has top label $\NN$ or $\RR$, but top label $\LL$ in $x_0f$.
The relevant rows of the weight table are:
\begin{center}
\begin{tabular}{c}
\\
top\\
label
\end{tabular}%
\begin{tabular}{|c|cccc|}
\multicolumn{5}{c}{bottom label}\\
\hline\rule{0ex}{2.5ex}
&$\;\LL\;$&$\;\NN\;$&$\;\RR\;$&$\;\II\;$\\
\hline
$\;\LL\;$&$\boldsymbol{2}$&$1$&$1$\rule{0ex}{2.5ex}&$\boldsymbol{1}$\\
$\NN$&$\boldsymbol{1}$&$2$&$2$&$2$\\
$\RR$&$\boldsymbol{1}$&$2$&$2$&$\boldsymbol{0}$\\
\hline
\end{tabular}%
\end{center}
Each entry of the $\NN$ and $\RR$ rows differs from the corresponding entry
of the $\LL$ row by exactly one. In particular, moving from an $\RR$ or $\NN$
row to an $\LL$ row only increases the weight when in the $\LL$ column or when
starting at~\pair{\RR}{\II}.
\end{proof}

\begin{corollary}
Let $f\in F$. Then $\ell(x_0^{-1}f)<\ell(f)$ if and only if one of the following conditions holds:
\begin{enumerate}
\item $x_0^{-1}f$ has smaller support than $f$.
\item The left space of $f$ has label \pair{\LL}{\LL}.
\item The left space of $f$ has label \pair{\LL}{\II}, and the current tree is trivial.
\end{enumerate}
\end{corollary}

\begin{proposition}
Let $f\in F$. If left-multiplying $f$ by $x_1$ cancels a caret from the
bottom forest, then $\ell(x_1f)=\ell(f)-1$.
\end{proposition}

\begin{proof}
Clearly $\ell_1(x_1f)=\ell_1(f)-1$. We must show that $\ell_0$ remains unchanged.

Note first that the right space of $f$ is destroyed. This space has label \pair{\RR}{\II},
and hence has weight $0$. Therefore, its destruction does not affect $\ell_0$.

The only other space affected is the left space of $f$. If this space is
not in the support of $f$, it remains unlabeled throughout. Otherwise,
observe that it must have top label $\LL$ in both $f$ and $x_1f$.
The relevant row of the weight table is:
\begin{equation*}
\begin{tabular}{|c|cccc|}
\hline\rule{0ex}{2.5ex}
&$\LL$&$\NN$&$\RR$&$\II$\\
\hline
$\LL$&$2$&$1$&$1$&$1$\rule{0ex}{2.5ex}\\
\hline
\end{tabular}%
\end{equation*}
In particular, the only important property of the bottom label is whether or
not it is an $\LL$. This property is unaffected by the deletion of the caret.
\end{proof}

\begin{proposition}
Let $f\in F$, and suppose that left-multiplying $f$ by $x_1$ creates a caret
in the top forest. Then $\ell(x_1f)=\ell(f)\pm 1$. Specifically, $\ell(x_1f)=\ell(f)-1$
if and only if the right space of $f$ has label \pair{\RR}{\RR}.
\end{proposition}

\begin{proof}
Clearly $\ell_1(x_1f)=\ell_1(f)+1$. As for $\ell_0$, observe that the only space whose
label could change is the right space of $f$.

\emph{Case 1}:\quad Suppose $x_1f$ has larger support than $f$. Then the right space of
$f$ is unlabeled, but has label \pair{\II}{\RR} in $x_1f$. This does not affect the value of $\ell_0$.

\emph{Case 2}:\quad Otherwise, note that the right space of $f$ has top label $\NN$ or
$\RR$. If the top label is an $\NN$, it remains and $\NN$ in $x_1f$. If it is an
$\RR$, then it changes to an $\II$. The relevant rows of the weight table are:
\begin{center}
\begin{tabular}{c}
\\
top\\
label
\end{tabular}%
\begin{tabular}{|c|cccc|}
\multicolumn{5}{c}{bottom label}\\
\hline\rule{0ex}{2.5ex}
&$\LL$&$\NN$&$\RR$&$\II$\\
\hline
$\RR$&$1$&$2$&$\boldsymbol{2}$&$0$\rule{0ex}{2.5ex}\\
$\II$&$1$&$2$&$\boldsymbol{0}$&$0$\\
\hline
\end{tabular}%
\end{center}
Observe that the weight decreases by two if the bottom label is an $\RR$,
and stays the same otherwise.
\end{proof}

We have now verified condition (2). Also, we have gathered enough
information to verify condition (3):

\begin{theorem}
Let $f\in F$ be a nonidentity element.
\begin{enumerate}
\item If current tree of $f$ is nontrivial, then either $\ell(x_1^{-1}f)<\ell(f)$,
or \mbox{$\ell(x_0f)<\ell(f)$}.
\item If left-multiplication by $x_1$ would remove a caret from the bottom tree, then
$\ell(x_1f)<\ell(f)$.
\item Otherwise, either $\ell(x_0f)<\ell(f)$ or $\ell(x_0^{-1}f)<\ell(f)$.
\end{enumerate}
\end{theorem}

\begin{proof}

\

\emph{Statement 1}:\quad If $\ell(x_1^{-1}f)>\ell(f)$, then the right space of $x_1^{-1}f$
has type~\pair{\RR}{\RR}. The right space of $f$ therefore has type \pair{\RR\text{ or }\NN}{\RR\text{ or }\NN},
so that \mbox{$\ell(x_0f)<\ell(f)$}.

\emph{Statement 2}:\quad See proposition 4.3.5.

\emph{Statement 3}:\quad Suppose $\ell(x_0f)>\ell(f)$. There are three cases:

\emph{Case 1}:\quad The right space of $f$ is not in the support of $f$. Then the left space
of $f$ has label \pair{\LL}{\RR}, \pair{\LL}{\LL}, or \pair{\LL}{\II}. In all three cases,
$\ell(x_0^{-1}f)<\ell(f)$.

\emph{Case 2}:\quad The right space of $f$ has bottom label $\LL$, and
right-multiplication by $x_0$ does not remove this space from the support.
Then the left space of $f$ must have label \pair{\LL}{\LL} or \pair{\LL}{\II},
and hence $\ell(x_0^{-1}f)<\ell(f)$.

\emph{Case 3}:\quad The right space of $f$ has label \pair{\RR}{\II}. Then the tree immediately to the right
of the top pointer is trivial, and the bottom leaf under it is a right leaf. If the bottom leaf
under the top pointer were a left leaf, then left-multiplying $f$ by $x_1$ would cancel a caret.
Hence, it is also a right leaf, so the left space of $f$ has label \pair{\LL}{\II}. We conclude that
$\ell(x_0^{-1}f)<\ell(f).$
\end{proof}

\begin{corollary}
Let $f\in F$, and let $\mathfrak{f}$ be the reduced forest diagram for $f$. Then there exists
a minimum-length word $w$ for $f$ with the following properties:
\begin{enumerate}
\item Each instance of $x_1$ in $w$ creates a top caret of $\mathfrak{f}$.
\item Each instance of $x_1^{-1}$ in $w$ creates a bottom caret of $\mathfrak{f}$.
\end{enumerate}
In particular, $w$ has $\ell_1(f)$ instances of $x_1$ or $x_1^{-1}$, and $\ell_0(f)$
instances of $x_0$~or~$x_0^{-1}$.
\end{corollary}

\begin{proof}
By the previous theorem, it is always possible to travel from $f$ to the identity in such
a way that each left-multiplication by $x_1$ deletes a bottom caret and each left-multiplication
by $x_1^{-1}$ deletes a top caret.
\end{proof}

Of course, not every minimum-length word for $f$ is of the given form. We will
discuss this phenomenon in the next section.

\

\section{Minimum-Length Words}

In principle, the results from the last section specify an algorithm for
finding minimum-length words.  (Given an element, find a generator which
shortens it.  Repeat.)  In practice, though, no algorithm is necessary:
one can usually guess a minimum-length word by staring at the forest diagram.
Our goal in this section is to convey this intuition.

\begin{example}
Let $f$ be the element of $F$ with forest diagram:
\begin{center}
\includegraphics{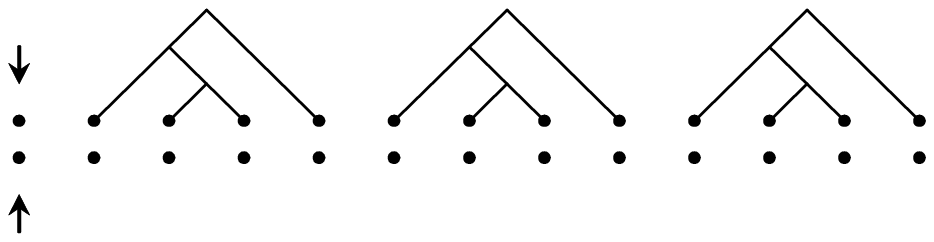}
\end{center}
Then there is exactly one minimum-length word for $f$, namely:
\begin{equation*}
x_0^{-3} u x_0 u x_0 u x_0
\end{equation*}
where $u = x_1^2 x_0^{-1} x_1 x_0$. Note that the trees of $f$ are constructed
from \emph{left to right}.

Similarly, $f^{-1}$ has forest diagram:
\begin{center}
\includegraphics{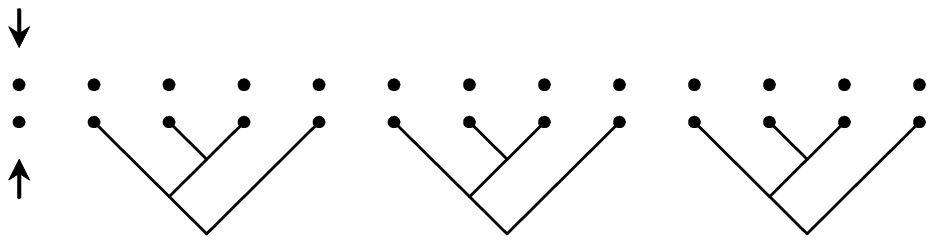}
\end{center}
and the only minimum-length word for $f^{-1}$ is:
\begin{equation*}
x_0^{-1} u^{-1} x_0^{-1} u^{-1} x_0^{-1} u^{-1} x_0^3
\end{equation*}
Note that the trees of $f^{-1}$ are constructed from \emph{right to left}.
\end{example}

\begin{example}
\quad Let $f$ be the element of $F$ with forest diagram:
\begin{center}
\includegraphics{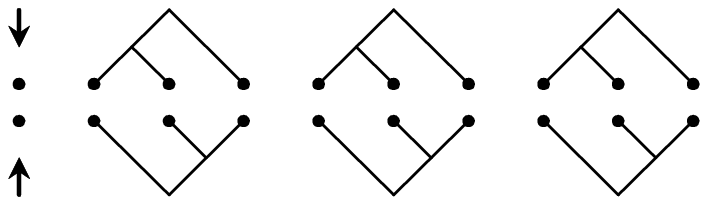}
\end{center}
There are precisely four minimum-length words for $f$:
\begin{equation*}
\begin{aligned}[t]
&x_0^{-3}\,v\,x_0\,v\,x_0\,v\,x_0\\
&x_0^{-1}\,v\,x_0^{-2}\,v\,x_0\,v\,x_0^2\\
&x_0^{-2}\,v\,x_0^{-1}\,v\,x_0^2\,v\,x_0\\
&x_0^{-1}\,v\,x_0^{-1}\,v\,x_0^{-1}\,v\,x_0^3
\end{aligned}
\end{equation*}
where $v = x_1^2 x_0^{-1} x_1^{-1} x_0 x_1^{-1}$. In particular, each of the
first two components can be constructed either when the pointer is moving
right, or later when the pointer is moving back left.
\end{example}

\begin{example}
Let $f$ be the element of $F$ with forest diagram:
\begin{center}
\includegraphics{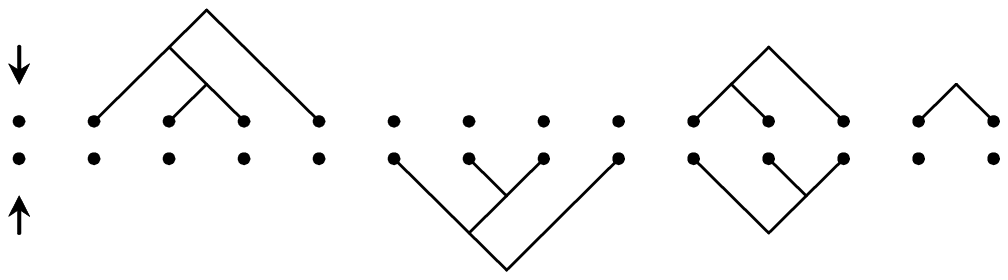}
\end{center}
There are precisely two minimum-length words for $f$:
\begin{equation*}
\begin{aligned}[t]
&x_0^{-2}\,u^{-1}\,x_0^{-2}\,x_1\,x_0\,v\,x_0^2\,u\,x_0\\
&x_0^{-2}\,u^{-1}\,x_0^{-1}\,v\,x_0^{-1}\,x_1\,x_0^3\,u\,x_0
\end{aligned}
\end{equation*}
where $u = x_1^2 x_0^{-1} x_1 x_0$ and $v = x_1^2 x_0^{-1} x_1^{-1} x_0 x_1^{-1}$. Note
that the first component must always be constructed on the journey right, and
the second component must always be constructed on the journey left. The
only choice lies with the construction of the third component: should it be
constructed when moving right, or should it be constructed while moving back
left?
\end{example}

In general, certain components act like ``top trees'' while others act like
``bottom trees'', while still others are ``balanced''. For example, the
forest diagram:
\begin{center}
\includegraphics{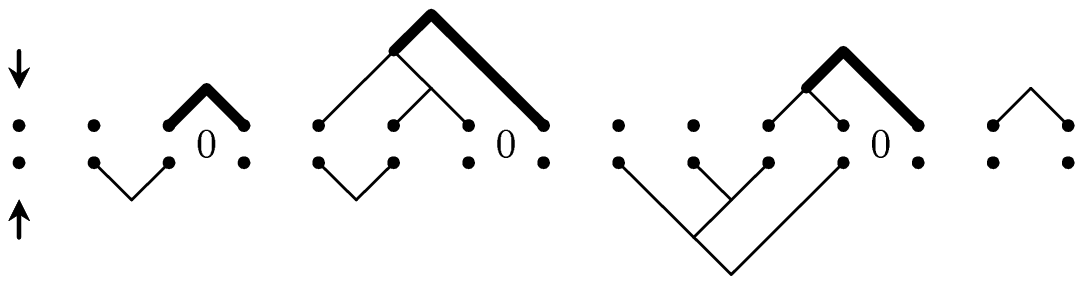}
\end{center}
must be constructed from left to right (so all the components act like ``top
trees''). The reason is that the three marked spaces each have weight $0$,
so that each of the three highlighted carets must be constructed \emph{before}
the pointer can move farther to the right. Essentially, the highlighted carets
are acting like \emph{bridges} over these spaces.

The idea of the ``bridge'' explains two phenomena we have already observed.
First, consider the following contrapositive of proposition 4.3.6:

\begin{proposition}
Let $f\in F$, and suppose that the top pointer of $f$ points at a
nontrivial tree. Then $\ell\left(x_1^{-1}f\right)<\ell(f)$
unless the resulting uncovered space has type \pair{\RR}{\RR}.\quad\qedsymbol
\end{proposition}

This proposition states conditions under which the destruction of a top caret
decreases the length of an element. Essentially, the content of the
proposition is that it makes sense to delete a top caret \emph{unless that caret is
functioning as a bridge}. (Note that the deletion of any of the bridges in
the example above would result in an \pair{\RR}{\RR} space.) It makes no sense to
delete a bridge, since the bridge is helping you access material further to the right.

Next, recall the statement of corollary 4.3.8: every $f \in F$ has a minimum-length word
with $\ell_1(f)$ instances of $x_1$ or $x_1^{-1}$ and $\ell_0(f)$ instances of $x_0$ or $x_0^{-1}$.
After the corollary, we mentioned that not every minimum-length word for $f$ is necessarily
of this form. The reason is that it sometimes makes sense to build bridges
during the creation of an element:

\begin{example}
Let $f$ be the element of $F$ with forest diagram:
\begin{center}
\includegraphics{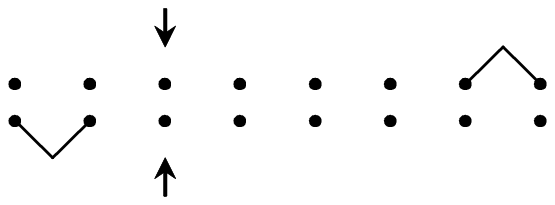}
\end{center}
Then one minimum-length word for $f$ is:
\begin{equation*}
x_0^2 x_1^{-1} x_0^{-5} x_1 x_0^4
\end{equation*}
This word corresponds to the instructions ``move right, create the top caret,
move left, create the bottom caret, and then move back to the origin''.
However, here is another minimum-length word for $f$:
\begin{equation*}
x_0^2 x_1^{-1} (x_0^{-1}x_1^{-3}x_0^{-1}) x_1 (x_0x_1^3)
\end{equation*}
In this word, the ``move right'' is accomplished by building three temporary
bridges:
\begin{center}
\includegraphics{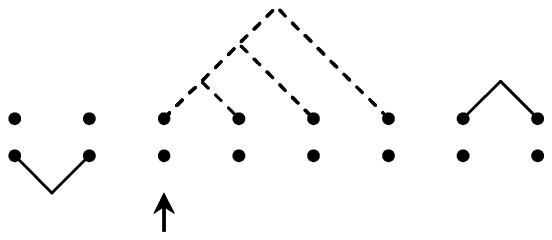}
\end{center}
These bridges are torn down during the ``move left''.

Finally, here is a third minimum-length word for $f$:
\begin{equation*}
x_1^{-3} x_0^2 x_1^{-1} x_0^{-2} x_1 (x_0x_1^3)
\end{equation*}
In this word, bridges are again built during the ``move right'', but they
aren't torn down until the very end of the construction.
\end{example}

We now turn our attention to a few examples with some more complicated
behavior.

\begin{example}
\quad Let $f$ be the element of $F$ with forest diagram:
\begin{center}
\includegraphics{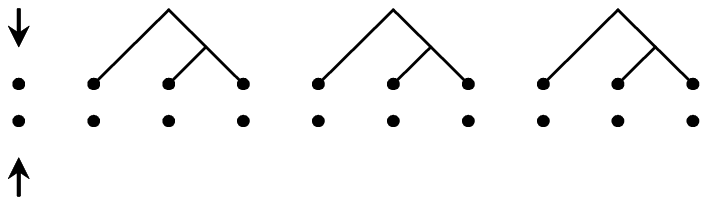}
\end{center}
There are four different minimum-length words for $f$:
\begin{equation*}
\begin{aligned}[t]
&x_0^{-3}x_1x_0^{-1}x_1x_0^2x_1x_0^{-1}x_1x_0^2x_1x_0^{-1}x_1x_0^2\\
&x_0^{-1}x_1x_0^{-3}x_1x_0^{-1}x_1x_0^2x_1x_0^{-1}x_1x_0^2x_1x_0^2\\
&x_0^{-2}x_1x_0^{-2}x_1x_0^{-1}x_1x_0^2x_1x_0^2x_1x_0^{-1}x_1x_0^2\\
&x_0^{-1}x_1x_0^{-2}x_1x_0^{-2}x_1x_0^{-1}x_1x_0^2x_1x_0^2x_1x_0^2
\end{aligned}
\end{equation*}
Note that each of the first two components may be either partially or fully
constructed during the move to the right. This occurs because the trees in
this example do not end with bridges. (Compare with example 4.4.1.)
\end{example}

\begin{example}
\quad Let $f$ be the element of $F$ with forest diagram:
\begin{center}
\includegraphics{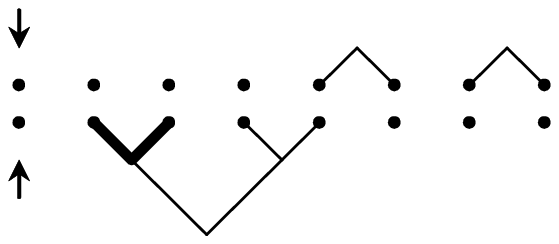}
\end{center}
There is exactly one minimum-length word for $f$:
\begin{equation*}
x_0^{-1} x_1^{-1} x_0^{-3} x_1 x_0 x_1 x_0 x_1^{-1} x_0 x_1^{-1} x_0
\end{equation*}
Note that the highlighted caret must be constructed \emph{last}, since the space it
spans should not be crossed. However, we must begin by partially constructing the
first component, because of the bridge on its right end.
\end{example}

\

\chapter{Applications}

This chapter contains various applications of forest diagrams and the length
formula. Some of the material in this chapter represents joint work with my thesis advisor,
Kenneth Brown. Sections 1 and 2 were originally published in \cite{BeBr}.

\section{Dead Ends and Deep Pockets}

In \cite{ClTa1}, S.~Cleary and J.~Taback prove that $F$ has ``dead
ends'' but no ``deep pockets''. In this section, we show how forest
diagrams can be used to understand these results.

\begin{definition}
A \emph{dead end} is an element $f\in F$ such that
\mbox{$\ell(xf)<\ell(f)$} for all $x\in\left\{x_0, x_1, x_0^{-1}, x_1^{-1}\right\}$.
\end{definition}

\begin{example}
Consider the element $f\in F$ with forest diagram:
\begin{center}
\includegraphics{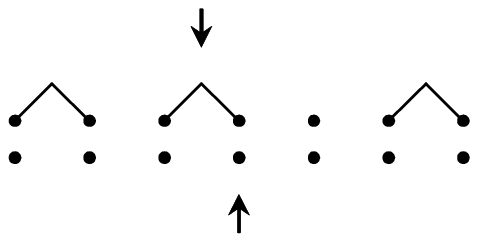}
\end{center}
Left-multiplying by $x_0^{-1}$ decreases the length since the left space of
$f$ is of type \pair{\LL}{\LL}. Left-multiplying by $x_0$ or $x_1$ decreases the
length since the right space of $f$ is of type \pair{\RR}{\RR}. Finally, left-multiplying
by $x_1^{-1}$ decreases the length since it deletes a top caret and the right
space of $x_1^{-1}f$ is not of type~\pair{\RR}{\RR}.
\end{example}

This example is typical:

\begin{proposition}
Let $f\in F$. Then $f$ is a dead end if and only if:
\begin{enumerate}
\item The current tree of $f$ is nontrivial,
\item The left space of $f$ has label \pair{\LL}{\LL},
\item The right space of $f$ has label \pair{\RR}{\RR}, and
\item The right space of $x_1^{-1}f$ does not have label \smash[t]{\pair{\RR}{\RR}}.
\end{enumerate}
\end{proposition}

\begin{proof}
The ``if'' direction is trivial. To prove the ``only if''
direction, assume that $f$ is a dead end. Then:
\begin{description}
\item[Condition (1)] follows from the fact that
$\ell\left(x_1^{-1}f\right)<\ell(f)$ (see proposition 4.3.5).
\item[Condition (2)] now follows from the fact that
$\ell\left(x_0^{-1}f\right)<\ell(f)$ (see corollary 4.3.4).
\item[Condition (3)] now follows from the fact that
$\ell(x_1f)<\ell(f)$ (see proposition 4.3.6).
\item[Condition (4)] now follows from the fact that
$\ell\left(x_1^{-1}f\right)<\ell(x_1f)$ (see proposition 4.3.6).\hfill\qedsymbol
\end{description}\renewcommand{\qedsymbol}{}
\end{proof}

Note that there are several ways to meet condition (4): the right space of
$x_1^{-1}f$ could be of type \pair{\RR}{\LL} (as in example 5.1.2), or it could
be of type \pair{\RR}{\II}:
\begin{center}
\includegraphics{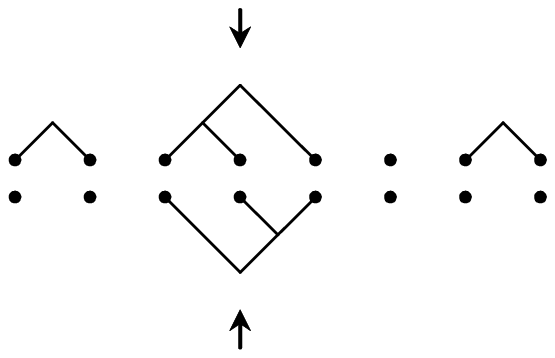}
\end{center}
or it could just have an $\NN$ on top:
\begin{center}
\includegraphics{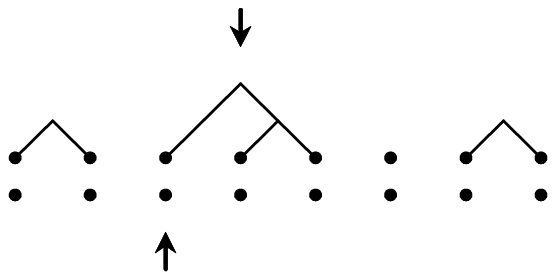}
\end{center}

Notice also that the proof of proposition 5.1.3 did not use the hypothesis that
\mbox{$\ell\left(x_0f\right)<\ell(f)$}.  In particular, if the length of $f$ increases
when you left-multiply by $x_1$, $x_1^{-1}$, and $x_0^{-1}$, then $f$ must be a dead end.

\begin{definition}
Let $k\in\mathbb{N}$. A \emph{$k$-pocket} of $F$ is an element
$f\in F$ such that:
\begin{equation*}
\ell(s_1\cdots s_kf)\leq\ell(f)
\end{equation*}
for all $s_1,\ldots,s_k\in\left\{x_0,x_1,x_0^{-1},x_1^{-1},1\right\}$.
\end{definition}

Note that a 2-pocket is just a dead end. S.~Cleary and J.~Taback demonstrated
that $F$ has no $k$-pockets for $k\geq 3$. We give an alternate proof:

\begin{proposition}
$F$ has no $k$-pockets for $k\geq 3$.
\end{proposition}

\begin{proof}
Let $f\in F$ be a dead-end element. Then the right space of $f$
has label \pair{\RR}{\RR}, so the tree to the right of the top pointer is trivial.
Therefore, repeatedly left-multiplying $x_0f$ by $x_1^{-1}$ will create
negative carets:
\begin{center}
\includegraphics{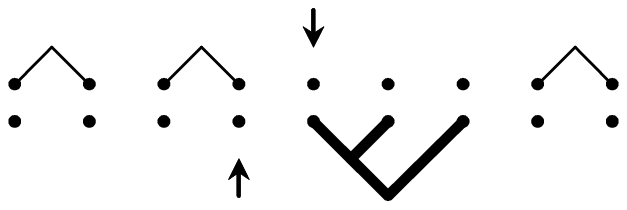}
\end{center}
In particular, $x_1^{-1}x_1^{-1}x_0f$ has length $\ell(f)+1$.
\end{proof}

\

\section{Growth}

We can use forest diagrams to calculate the growth function of the positive monoid with respect
to the $\{x_0,x_1\}$ generating set.  Burillo \cite{Bur} recently arrived at the same result using
tree diagrams and Fordham's length formula:

\begin{theorem}Let $p_n$ denote the number of positive elements of length $n$, and let:
\begin{equation*}
p(x) = \sum_{n=0}^\infty p_n x^n
\end{equation*}
Then:
\begin{equation*}
p(x) = \frac{1 - x^2}{1 - 2x - x^2 + x^3}
\end{equation*}
In particular, $p_n$ satisfies the recurrence relation:
\begin{equation*}
p_n = 2 p_{n-1} + p_{n-2} - p_{n-3}
\end{equation*}
for all $n \geq 3$.
\end{theorem}

\begin{proof}
Let $P_n$ be the set of all positive elements of length $n$.  Define four subsets of $P_n$ as follows:
\begin{enumerate}
\item $A_n = \{f\in P_n :$ the current tree of $f$ is trivial and is not the leftmost tree$\}$
\item $B_n = \{f\in P_n :$ the current tree of $f$ is nontrivial, but its right subtree is trivial$\}$
\item $C_n = \{f\in P_n :$ the current tree of $f$ is trivial and is the leftmost tree.$\}$
\item $D_n = \{f\in P_n :$ the current tree of $f$ is nontrivial, and so is its right subtree.$\}$
\end{enumerate}

Given an element of $A_n$, we can remove the current tree and move the pointer left, like this:
\begin{center}
\includegraphics{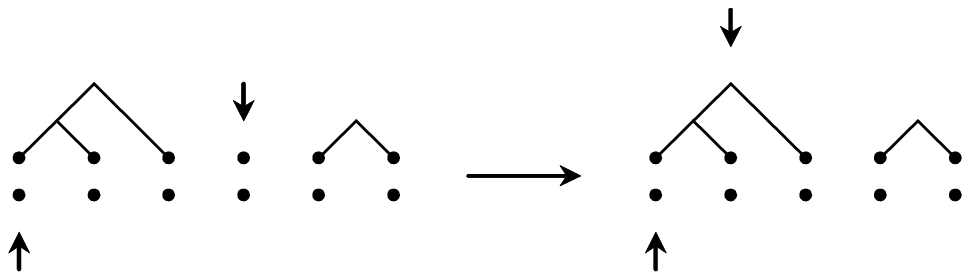}
\end{center}
This defines a bijection $A_n\rightarrow P_{n-1}$, so that:
\begin{equation*}
|A_n|=|P_{n-1}|
\end{equation*}

Given an element of $B_n$, we can remove the top caret together with the resulting trivial tree, like this:
\begin{center}
\includegraphics{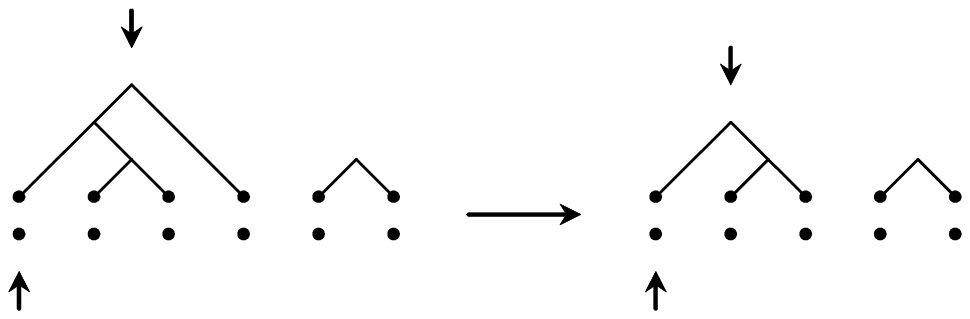}
\end{center}
This defines a bijection $B_n\rightarrow P_{n-1}$, so that:
\begin{equation*}
|B_n|=|P_{n-1}|
\end{equation*}

Given an element of $C_n$, we can move both the top and bottom arrows one space to the right, like this:
\begin{center}
\includegraphics{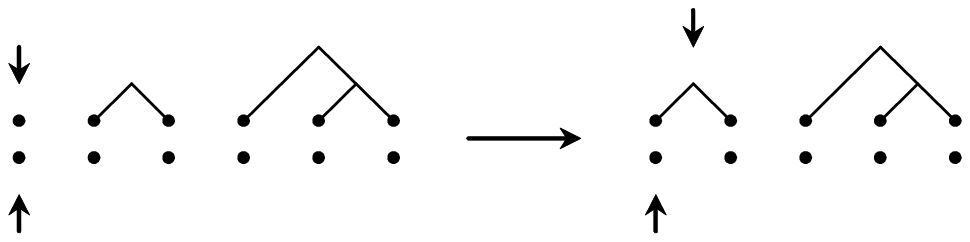}
\end{center}
When $n \geq 2$, this defines an injection $\varphi\colon C_n\rightarrow P_{n-2}$.  The image of $\varphi$ is all elements of
$P_{n-2}$ whose current tree is the first tree.

Finally, given an element of $D_n$, we can remove the top caret and move the pointer to the right subtree, like this:
\begin{center}
\includegraphics{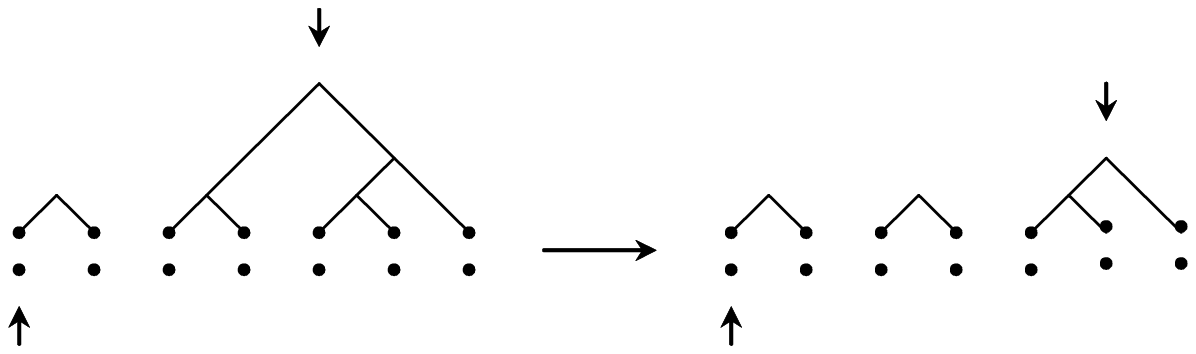}
\end{center}
This defines an injection $\psi\colon D_n\rightarrow P_{n-2}$.  The image of $\psi$ is all elements of
$P_{n-2}$ whose current tree is nontrivial, and is not the first tree.  In particular:
\begin{equation*}
(\text{im}\, \varphi) \cup (\text{im}\, \psi) = P_{n-2}-A_{n-2}
\end{equation*}
so that:
\begin{equation*}
|C_n|+|D_n|=|P_{n-2}|-|A_{n-2}|=|P_{n-2}|-|P_{n-3}|
\end{equation*}
This proves that $p_n$ satisfies the given recurrence relation for $n \geq 3$.
It is not much more work to verify the given expression for $p(x)$.
\end{proof}

\section{The Isoperimetric Constant}

Let $G$ be a group with finite generating set $\Sigma$, and let $\Gamma$
denote the Cayley graph of $G$ with respect to $\Sigma$. If $S\subset G$,
define:
\begin{equation*}
\delta S=\left\{\text{edges in }\Gamma\text{ between }S\text{ and }S^c\right\}
\end{equation*}
The \emph{isoperimetric constant} of $G$ is defined as follows:
\begin{equation*}
\iota\left(G,\Sigma\right)=\inf\left\{\frac{\left|\delta
S\right|}{\left|S\right|}:S\subset G\text{ and }\left|S\right|<\infty\right\}
\end{equation*}

\begin{theorem}[F\o lner]
The group $G$ is amenable if and only if $\iota\left(G,\Sigma\right)=0$.
\end{theorem}
\begin{proof}See \cite{Wag}. \end{proof}

Guba \lbrack Guba\rbrack\ has shown that $\iota\bigl(F,\left\{x_0,x_1\right\}\bigr)\leq 1$.
In this section, we shall prove a slightly better estimate:

\begin{theorem}  $\iota\bigl(F,\left\{x_0,x_1\right\}\bigr)\leq 1/2$.
\end{theorem}

The proof will occupy the remainder of this section.

Define the \emph{height} of a binary tree to be the length of the longest
descending path starting at the root and ending at a leaf. Define the \emph{width}
of a binary forest to be the number of spaces in its support. For each
$n,k\in\mathbb{N}$, let $S_{n,k}$ denote all positive elements whose forest
diagram has width at most $n$ and all of whose trees have height at most $k$.
We shall show that:
\begin{equation*}
\lim_{k\rightarrow\infty}\lim_{n\rightarrow\infty}\frac{\left|\delta
S_{n,k}\right|}{\left|S_{n,k}\right|}=\frac{1}{2}
\end{equation*}

First of all, observe that each element of $S_{n,k}$ can be represented by
a finite, $n$-space binary forest together with a pointer pointing
to one of the trees:
\begin{center}
\includegraphics{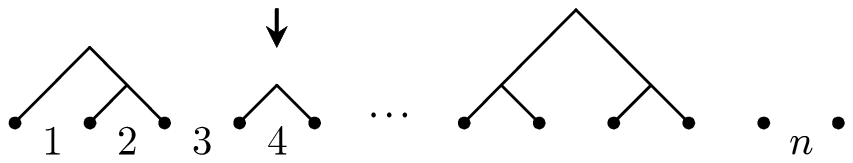}
\end{center}
We shall refer to such an object as a \emph{pointed forest}.
Note that the trivial trees on the right of this picture count as part of this
pointed forest, even though they are not in the support of the standard forest
diagram.

\begin{claim} Let $k\in\mathbb{N}$. If $f$ is a randomly chosen element of
$S_{n,k}$, then:
\begin{equation*}
\frac{|\delta S_{n,k}|}{|S_{n,k}|} - 2\lim_{n\rightarrow\infty}P\left(x_1^{-1}f\notin S_{n,k}\right)
\rightarrow 0
\end{equation*}
as $n\rightarrow \infty$.
\end{claim}
\begin{proof} Observe that:
\begin{equation*}
\frac{\left|\delta S_{n,k}\right|}{\left|S_{n,k}\right|}=P(x_0f\notin
S_{n,k})+P\bigl(x_0^{-1}f\notin S_{n,k}\bigr)+P(x_1f\notin
S_{n,k})+P\bigl(x_1^{-1}f\notin S_{n,k}\bigr)
\end{equation*}
Now, $S_{n,k}$ must have the same number of incoming and outgoing edges of
each type, so both $x_0$ terms and both $x_1$ terms must be equal.
Therefore:
\begin{equation*}
\frac{\left|\delta
S_{n,k}\right|}{\left|S_{n,k}\right|}=2\,P\bigl(x_0^{-1}f\notin
S_{n,k}\bigr)+2\,P\bigl(x_1^{-1}f\notin S_{n,k}\bigr)
\end{equation*}

Next, note that $x_0^{-1}f\notin S_{n,k}$ if and only if the current
tree of $f$ is the leftmost tree. However, as $n\rightarrow\infty$ the
minimum number of trees in each element of $S_{n,k}$ goes to $\infty$, and
hence the probability that the current tree is the leftmost tree goes to $0$.
Therefore:
\begin{equation*}
\lim_{n\rightarrow\infty}P\bigl(x_0^{-1}f\notin
S_{n,k}\bigr)=0\tag*{\qedhere}
\end{equation*}
\end{proof}

Now, if $f\in S_{n,k}$, then $x_1^{-1}f\notin S_{n,k}$ if and only if
the current tree of $f$ is trivial. Therefore, we must determine the
probability that the current tree of a random pointed forest of
width $n$ and height at most $k$ is the trivial tree.  (Here and elsewhere the word \emph{random}
means randomly chosen with respect to the uniform distribution on pointed forests.)

\begin{theorem}Suppose we choose a random pointed forest $f$ with $n$ leaves and height at most $k$.
Then the limit:
\begin{equation*}
\lim_{n\rightarrow\infty} P(\text{the current tree of $f$ is trivial})
\end{equation*}
exists and is the unique positive root of the polynomial equation:
\begin{equation*}
t_{1,k}p + t_{2,k}p^2+t_{3,k}p^3 + \cdots = 1
\end{equation*}
where $t_{\ell,k}$ is the number of binary trees with $\ell$ leaves and height at most $k$.
\end{theorem}
\begin{proof}Let $f_n$ denote the number of binary forests with $n$ leaves and height at most $k$.
Then $f_n$ satisfies the following recurrence relation:
\begin{equation*}f_n = t_{1,k} f_{n-1} + t_{2,k} f_{n-2} + \cdots
\end{equation*}
Observe that $t_{n,k} \neq 0$ for $0<n\leq 2^k$ and $t_{n,k} = 0$ for $n > 2^k$.  Using the standard theory
of linear recurrence relations (see \cite{Bru}), we deduce that:
\begin{equation*}
\lim_{n \rightarrow \infty} \frac{f_{n-1}}{f_n} = p
\end{equation*}
where $p$ is the unique positive root of the polynomial equation above.

Now let $R_n$ be the number of pointed forests with $n$ leaves and height at most $k$, and let $R_n^*$ be the number
of such pointed forests whose current tree is trivial.  Then:
\begin{equation*}
R_n = f_1 f_{n-1} + f_2 f_{n-2} + \cdots + f_n f_0
\end{equation*}
and:
\begin{equation*}
R_n^* = f_0 f_{n-1} + f_1 f_{n-2} + \cdots + f_{n-1} f_0
\end{equation*}
Therefore, the probability that the current tree is trivial is given by:
\begin{equation*}
\frac{R_n^*}{R_n} =
\frac{f_0 f_{n-1} + f_1 f_{n-2} + \cdots + f_{n-1} f_0}{f_1 f_{n-1} + f_2 f_{n-2} + \cdots + f_n f_0}
\end{equation*}
It is not hard to show that this approaches $p$ as $n\rightarrow \infty$.  In particular, if we ignore
the first term of the numerator and the middle term of the denominator, then each of the remaining terms
in the numerator is equal to $f_{m-1} / f_m$ times the corresponding term in the denominator for some $m>n/2$.
\end{proof}

It is interesting to note that the probability that the current tree is a single caret approaches $p^2$ as
$n \rightarrow \infty$, since the probability that the current tree is a single caret should be equal to the
probability that both the current tree and the right tree are trivial.  More generally, if $\sigma$ is a fixed
binary tree with $\ell$ leaves, the probability that the current tree is $\sigma$ approaches $p^\ell$ as
$n \rightarrow \infty$.  This gives us a nice intuitive understanding of the polynomial equation in the
preceding theorem.

Now, let $p_k$ denote the unique positive root of the equation:
\begin{equation*}
t_k(p_k) = 1
\end{equation*}
where $t_k$ is the polynomial from theorem 5.3.4:
\begin{equation*}
t_k(x) = t_{1,k}x + t_{2,k}x^2+t_{3,k}x^3 + \cdots = 1
\end{equation*}
All that remains is to show that $\displaystyle\lim_{k\rightarrow\infty}p_k=\frac{1}{4}$.

Note first that a binary tree has height at most $k$ if and only if its
left and right subtrees both have height at most $k-1$. Hence:
\begin{equation*}
t_k(x)=t_{k-1}(x)^2+x
\end{equation*}
where the $x$ term corresponds to the trivial tree. This lets us derive the
polynomials $t_k(x)$ iteratively, starting at $t_{-1}(x)=0$.

Therefore, to solve the equation:
\begin{equation*}
t_k(p_k)=1
\end{equation*}
we must investigate iteration of the map:
\begin{equation*}
t\mapsto t^2+c
\end{equation*}
In particular, $c=p_k$ if and only if we arrive at $1$ after $k+1$ iterations,
starting at $t=0$.

A graph of the equation $y=x^2+c$ is shown below for
$c=\displaystyle\frac{-1+\sqrt{5}}{2}$:
\begin{center}
\includegraphics{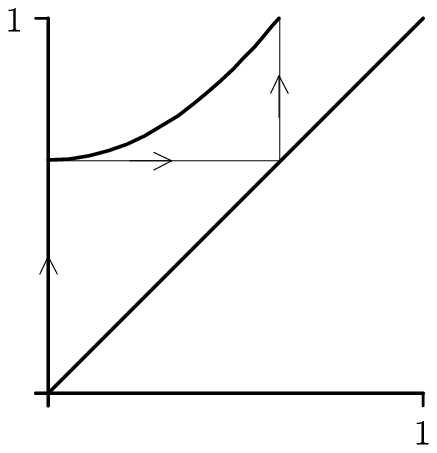}
\end{center}
Since this quadratic arrives at $1$ after two iterations,
$p_1=\displaystyle\frac{-1+\sqrt{5}}{2}$.

By decreasing $c$ (i.e. moving the parabola down), we can increase the
number of iterations that it takes to get to $1$, and hence find $p_k$ for
larger values of $k$. Here's a graph of $y=x^2+p_4$:
\begin{center}
\includegraphics{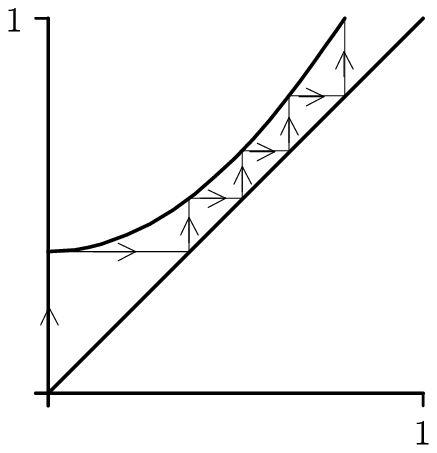}
\end{center}
As $k\rightarrow\infty$, we must lower the parabola $y=x^2+c$ arbitrarily
close to the line $y=x$. They become tangent at
$c=\displaystyle\frac{1}{4}$, so:
\begin{equation*}
\lim_{k\rightarrow\infty}p_k=\frac{1}{4}
\end{equation*}
This concludes the proof of the theorem.\quad\qedsymbol

\

The following corollary explains one reason that it is difficult to
improve upon this result:

\begin{corollary}  Let $\mathcal{T}$ be any finite set of binary trees
which is closed under the taking of subtrees, and let $S_{n,\mathcal{T}}$
denote all positive elements whose forest diagram has width at most $n$ and
all of whose trees are from $\mathcal{T}$. Then:
\begin{equation*}
\displaystyle\frac{\left|\delta
S_{n,\mathcal{T}}\right|}{\left|S_{n,\mathcal{T}}\right|}>\frac{1}{2}
\end{equation*}
\end{corollary}
\begin{proof}  Let $a_i$ be the number of trees in $\mathcal{T}$ with $i$
leaves, and let:
\begin{equation*}
a\left(x\right)=a_1x+a_2x^2+a_3x^3+\cdots
\end{equation*}
By the argument above,
\begin{equation*}
\frac{\left|\delta
S_{n,\mathcal{T}}\right|}{\left|S_{n,\mathcal{T}}\right|}=2p
\end{equation*}
where $p$ is the unique positive root of the polynomial equation
$a\left(p\right)=1$.

Let $k$ be the maximum height of the trees in $\mathcal{T}$. Then
$a_i\leq t_{i,k}$ for each $i$, so $a\left(x\right)\leq t_k\left(x\right)$
for all $x\geq 0$. \ Since $a\left(x\right)$ and $t_k\left(x\right)$ are both
increasing, we deduce that $p\geq p_k$, and so $p>1/4$.\end{proof}

\

\chapter{Convexity}

In this chapter, we prove that $F$ is not minimally almost convex
with respect to the generating set $\{x_0, x_1\}$. This improves upon a recent result
of S.~Cleary and J.~Taback \cite{ClTa2}.  The results in this chapter represent joint work with
Kai-Uwe Bux, and were originally published in \cite{BeBu}.

\section{Convexity Conditions}

A group $G$ is \emph{convex} (with respect to a given finite generating set) if the $n$-ball $B^n(G)$ is a convex subset
of the Cayley graph of $G$ for each $n$.  Very few groups are convex, but Cannon \cite{Can} has introduced
the following weaker property:

\begin{definition}
A group $G$ is \emph{almost convex} (with respect to a given finite generating set) if there exists an integer $L$
with the following property:  given any $g,h\in B^n(G)$ a distance two apart, there exists a path from
$g$ to $h$ in $B^n(G)$ of length at most $L$.
\end{definition}

There exist examples of groups which are almost convex with respect to one finite generating set, but not
with respect to another.

In \cite{Can}, Cannon gave an algorithm to construct arbitrarily large sections of the Cayley
graph of an almost convex group, thereby solving the word problem.
He also proved that groups of hyperbolic isometries, groups of Euclidean isometries, and small-cancellation groups
are almost convex. Coxeter groups are also almost convex \cite{DaSh}, as are all discrete groups
based on seven of the eight three-dimensional geometries \cite{ShSt}. Groups based on the Sol
geometry are not almost convex \cite{CFGT}, however, and neither are solvable
Baumslag-Solitar groups \cite{MiSh}.

The convexity of $F$ was first investigated by S.~Cleary and J.~Taback \cite{ClTa2}.
Using tree diagrams and Fordham's length formula, they proved the following:

\begin{theorem}$F$ is not almost convex with respect to the $\{x_0,x_1\}$ generating set.
\end{theorem}
\begin{proof}Let $f$ and $g$ be the following two elements, differing only in the position of
the top pointer:
\begin{center}
\includegraphics{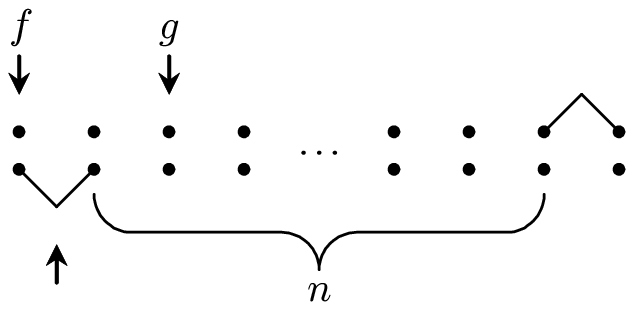}
\end{center}
Clearly $f$ and $g$ are a distance $2$ apart, and they both have length $2n+2$:
\begin{center}
\includegraphics{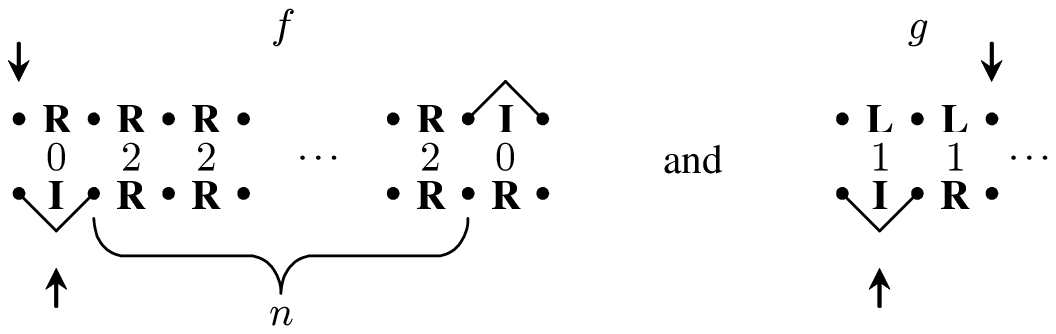}
\end{center}
However, $x_0f$ has length $2n+3$:
\begin{center}
\includegraphics{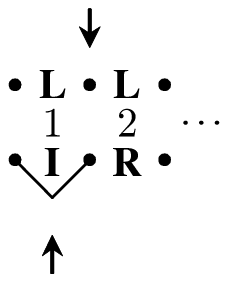}
\end{center}
Therefore, the geodesic:
\begin{equation*}
f \;\;\rule[3pt]{0.5in}{0.5pt}\;\; x_0 f \;\;\rule[3pt]{0.5in}{0.5pt}\;\; g
\end{equation*}
leaves the ball of radius $2n+2$.

In particular, if one wishes to go from $g$ to $f$ inside
the ball of radius $2n+2$, one must first move all the way to the right and delete the top
caret.  Therefore, the shortest path from $f$ to $g$ in $B^{2n+2}(F)$ has length at least $n$.
\end{proof}

The notion of almost convexity can be generalized as follows:

\begin{definition}Let $G$ be a finitely-generated group, and let $c\colon\mathbb{N}\rightarrow\mathbb{N}$
be any function.  We say that $G$ satisfies a
\emph{weak almost-convexity condition} with respect to $c$ if, given any $g,h\in B^n(G)$
a distance two apart, there is a path from $g$ to $h$ in $B^n(G)$ of length at most $c(n)$.
\end{definition}

Since there is always a path from $g$ to $h$ of length $2n$, the
weakest nontrivial convexity condition occurs when $c(n) = 2n-1$.  If $G$ satisfies this
condition (with respect to some finite generating set), we say that $G$
is \emph{minimally almost convex}.

I. Kapovich \cite{Kap} has shown that any minimally almost convex group is finitely
presented, and T. Riley \cite{Riley} derives upper bounds for the area function, the
isodiametric function, and the filling length for minimally almost convex groups.

In the next two sections, we will show that $F$ is not minimally almost convex.
In particular, we will prove the following:

\begin{theorem}For any even $n \geq 4$, there exist elements $l,r\in F$ of length $n$ such that:
\begin{enumerate}
\item $l$ and $r$ are distance two apart in the Cayley graph of $F$, and
\item The shortest path from $l$ to $r$ inside $B^n(F)$ has length $2n$.
\end{enumerate}
\end{theorem}

\section{$F$ is not Minimally Almost Convex}

Let $l$ and $r$ be the following two elements, differing only in the position of the top
pointer:
\begin{center}
\includegraphics{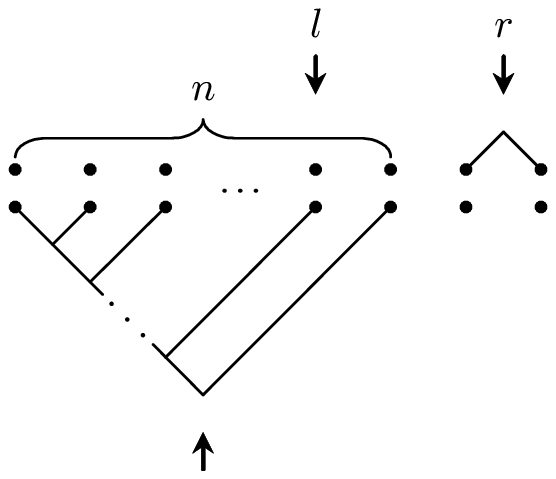}
\end{center}

Then $l$ and $r$ are clearly a distance two apart in the Cayley graph of $F$.

In this section and the next, we shall prove the following:

\begin{theorem} \
\begin{enumerate}
\item Both $l$ and $r$ have length $2n + 2$.
\item Any path from $l$ to $r$ inside the $(2n + 2)$-ball has length at least $4n + 4$.
\end{enumerate}
\end{theorem}

Condition (1) is trivial to verify:

\begin{lemma}  The elements $l$ and $r$ both have length $2n + 2$.
\end{lemma}
\begin{proof}  Note that the forest diagram for $l$ has exactly $n + 1$
carets. Furthermore, its forest diagram has the following weights:
\begin{center}
\includegraphics{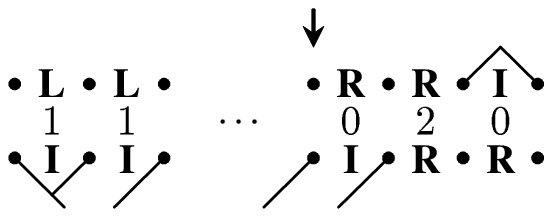}
\end{center}
Therefore, $l$ has length $(n+1) + \underset{n-1}{\underbrace{1+\cdots +1}} + 0 + 2 + 0 = 2n+2$.

Similarly, $r$ has exactly $n + 1$ carets. Its weights are:
\begin{center}
\includegraphics{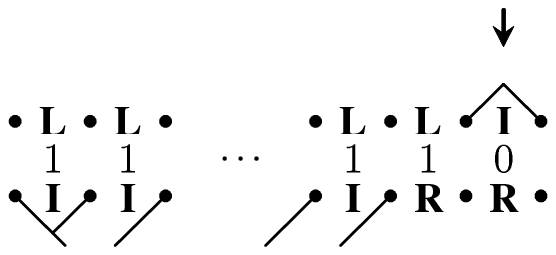}
\end{center}
Therefore, $r$ has length $(n+1) + \underset{n-1}{\underbrace{1+\cdots +1}} + 1 + 1 + 0 = 2n+2$.
\end{proof}

\begin{remark}  Note that the element $x_0 l$ has weights:
\begin{center}
\includegraphics{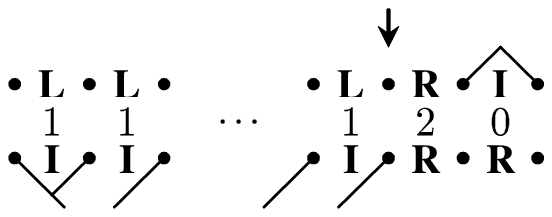}
\end{center}
and hence has length $2n + 3$. Therefore, the geodesic path
\begin{equation*}
l \;\;\rule[3pt]{0.5in}{0.5pt}\;\; x_0 l \;\;\rule[3pt]{0.5in}{0.5pt}\;\; r
\end{equation*}
leaves the ball of radius $2n + 2$.
\end{remark}

The proof of condition (2) is rather technical, so we postpone it until the next section.
For the remainder of this section, we shall attempt to convey the intuitive ideas behind
the proof, particularly in the choice of $l$ and $r$.

The main idea is as follows.  The forest diagram for $l$ and $r$ has a ``critical line'',
pictured below:
\begin{center}
\includegraphics{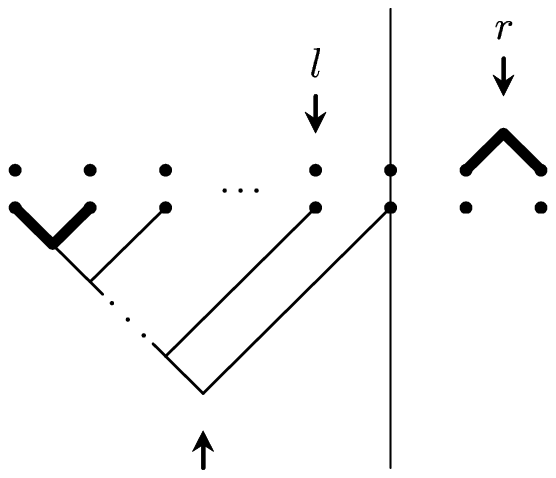}
\end{center}
This line has the following crucial property: if one wishes to remain in the $(2n + 2)$-ball,
one cannot cross the critical line as long as both of the highlighted carets are in place.
Therefore, any path in the $(2n + 2)$-ball from $l$ to $r$ must go through
the following four stages:
\begin{enumerate}
\item  Move to the left, and delete the leftmost caret.
\item  Move to the right (crossing the critical line), and delete the rightmost caret.
\item  Move back left (crossing again), and re-create the leftmost caret.
\item  Move back right (crossing the critical line for a third time), and
re-create the rightmost caret.
\end{enumerate}

\begin{example}  The word:
\begin{equation*}
\left( x_1 x_0^{n+1} \right) \left( x_1^{-1} x_0^{-n} \right)
\left( x_1^{-1} x_0^n \right)  \left( x_1 x_0^{1-n} \right)
\end{equation*}
describes a path in $B^{2n+2}(F)$ from $l$ to $r$ of length $4n + 4$.  Note that
the bulk of the bottom tree remains intact throughout this path. In
particular, this path does not pass through the identity vertex.
\end{example}

\begin{example}  The word
\begin{equation*}
\left( x_1 x_0^{n+1} \right) \left( x_1^{-n} x_0^{-1} \right)
\left( x_1^{-1} x_0 x_1^{n-1} \right)  \left( x_1 x_0^{1-n} \right)
\end{equation*}
represents a minimum-length path from $l$ to $r$ that passes through the identity vertex.
This time, we "travel to the right" by destroying the bottom tree, and "travel to the
left" be re-creating it.
\end{example}

\begin{example}  For $n = 8$, here is another minimum-length path from $l$ to $r$:
\begin{equation*}
\left( x_1 x_0^7 \underline{x_1^{-4}} x_0^2 \right) \left( x_1^{-1} x_0^{-4} \right)
\left( x_1^{-1} x_0^4 \right)  \left( x_1 x_0^{-2} \underline{x_1^4} x_0^{-5} \right)
\end{equation*}
In this path, we build carets in the top forest while moving to the left, and destroy
them later during the final move to the right. (Note that we have underlined the
segments of the word under discussion.) The resulting bridge saves us travel time
during the middle two stages, but its construction and demolition cost the same
amount of time during the first and last stages.
\end{example}

Finally, we would like to give some indication of how the elements $l$ and $r$ were chosen.
To do so, we give an example of two elements that would not work, despite having a similar
structure:

\begin{example}  Consider the elements $f$ and $g$ used during the proof of theorem 6.1.2:
\begin{center}
\includegraphics{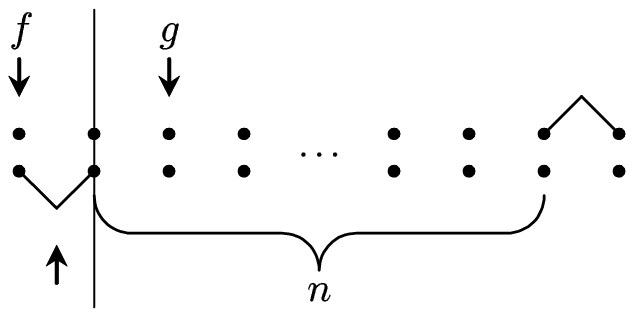}
\end{center}
We have previously observed that $f$ and $g$ are a distance $2$ apart in the Cayley graph
of $F$, and that they both have length $2n+2$.  Furthermore, the path
\begin{equation*}
f \;\;\rule[3pt]{0.5in}{0.5pt}\;\; x_0 f \;\;\rule[3pt]{0.5in}{0.5pt}\;\; g
\end{equation*}
leaves the $(2n+2)$-ball.  This suggests a ``critical line'' in the forest diagram
(already shown), which might lead one to believe that the word:
\begin{equation*}
\left( x_0^{1-n} x_1 \right) \left( x_0^{n+1} x_1^{-1} \right)
\left( x_0^{-n} x_1^{-1} \right)  \left( x_0^n x_1 \right)
\end{equation*}
of length $4n+4$ is a minimum-length path from $l'$ to $r'$ in $B^{2n+2}(F)$.

However, this turns out not to be the case. For example, when $n = 6$,
\begin{equation*}
\left( \underline{x_1^{-4}} x_0^{-1} x_1 \right) \left( x_0^3 x_1^{-1} \right)
\left( x_0^{-2} x_1^{-1} \right)  \left( x_0 \underline{x_1^4} x_0 x_1 \right)
\end{equation*}
is a path from $l'$ to $r'$ of length $24$ in $B_{18}$.  This path saves time
by building a bridge during the initial move to the right:
\begin{center}
\includegraphics{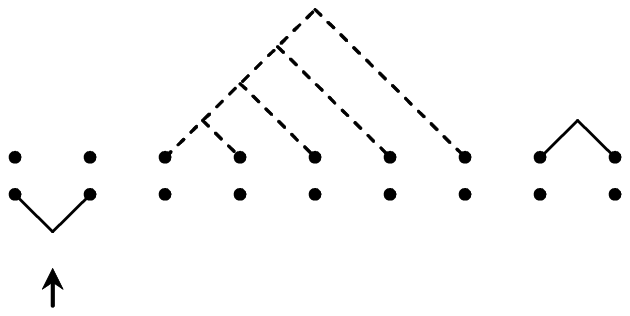}
\end{center}
This does not cost any extra time during the initial move to the right, but it saves
an enormous amount time during the subsequent two stages.  The bridge is torn down
during the final move to the left, which again does not cost any extra time.

The key difference between the elements $f,g$ and the element $l,r$ is that the bulk
of the support of $l$ and $r$ is to the \emph{left} of the critical line.
This means that one must travel left-right-left-right to get from $l$ to $r$, so
that one cannot save time by building bridges.

\end{example}

\section{Proof of Condition 2}

Fix a path $p$ from $l$ to $r$ that does not leave $B^{2n+2}(F)$. We wish to show that its length
$L(p)$ is at least $4n + 4$.

We claim that it suffices to show:

\begin{lemma}  On the path $p$, there are two vertices, $h_l$ and $h_r$, such that:
\begin{equation*}
d(h_l,h_r) \geq 2n+3
\end{equation*}
\end{lemma}

Why is this sufficient? Well clearly,
\begin{equation*}
L(p) \geq d(l,h_l) + d(h_l,h_r) + d(h_r,r)
\end{equation*}
However, by the triangle inequality:
\begin{equation*}
d(h_l,h_r) \leq d(h_l,l) + d(l,r) + d(r,h_r)
\end{equation*}
Since $d(l,r)=2$, we conclude that:
\begin{equation*}
L(p) \geq 2\, d(h_l,h_r)-2 \geq 4n+4
\end{equation*}

It remains to prove lemma 6.3.1.  We begin by formalizing the notion of ``crossing
the critical line'' from the previous section:

\begin{definition}Suppose $f\in F$.
\begin{enumerate}
\item Define the \emph{right foot} of $f$ to be the rightmost leaf of the current
    tree of $f$.
\item Define the \emph{critical leaf} of $f$ to be the rightmost leaf of the bottom tree
of $f$ currently indicated by the bottom pointer.
\end{enumerate}
\end{definition}

Note that the right foot of $l$ is to the left of the critical leaf and
the right foot of $r$ is to the right of the critical leaf. Let $h_l$ be the
first vertex of $p$ whose right foot coincides with the critical leaf, and let
$h_r$ be the last vertex of $p$ with this property.

\begin{remark} Note that left-multiplication by a generator can change the position of the right
foot by more than one unit.  However, since the tree directly above the critical leaf is trivial,
the right foot is guaranteed to not jump over the critical leaf.
\end{remark}

\begin{lemma}  The path $p$ ends with
\begin{equation*}
x_0^{-1} x_1^{-1} r \longrightarrow x_1^{-1} r \longrightarrow r
\end{equation*}
In particular, $h_r = x_0^{-1} x_1^{-1} r$.
\end{lemma}
\begin{proof}  It is easy to check that every path of length three emanating from $r$ either
passes through $x_0^{-1} x_1^{-1} r$ or leaves the $(2n+2)$-ball. \end{proof}

Recall that an element of $F$ is \emph{left-sided} if both pointers point to the rightmost
trees in the support.  Recall also that the \emph{width} $w(f)$ of an element of $F$
is the number of spaces in the support of its forest diagram.

\begin{lemma}  If $f \in F$ is left-sided, then:
\begin{equation*}
\ell(f) \geq 2 \, w(f)
\end{equation*}
\end{lemma}
\begin{proof}  We can associate to each caret in a forest diagram the interior space that it
covers. In this way, the interior spaces of $f$ each contribute $1$ to the length.
However, since $f$ is left-sided, every exterior space of $f$ is of type $\LL$.
The claim now follows, since the weight of any space is greater than or equal
to the number of $\LL$'s in its label pair. \end{proof}

\begin{remark} Note that this lemma \emph{fails} for right-sided elements:
An \pair{\RR}{\II}-space pair has weight $0$, and therefore only contributes $1$ to the length.
Hence, the best available estimate for right-sided elements is $\ell(f) \geq w(f)$.

This difference is related to the fact that one can ``move right'' by dropping carets,
but one cannot simultaneously build a structure and move left. (Compare with example 6.2.7.)
\end{remark}

\begin{proof}[Proof of Lemma 6.3.1]
Note that $h_r^{-1}$ has a trivial bottom forest:
\begin{center}
\includegraphics{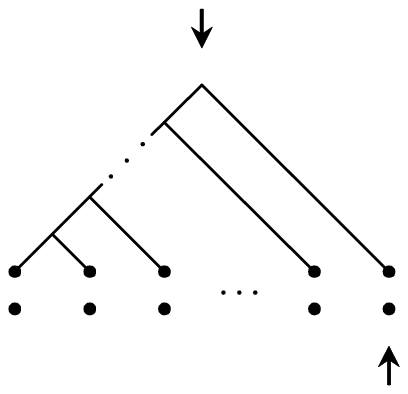}
\end{center}
Therefore, for any $f \in F$, a forest diagram for $h_r^{-1} f$ can be obtained
by stacking $h_r^{-1}$ on top of $f$.  Moreover, this diagram will be reduced
unless the bottom forest of $f$ has an exposed caret in exactly the right position
(namely, $n$ spaces to the left of the critical leaf) to cancel with the unique
exposed caret of $h_r^{-1}$.

Consider the element $x_2$.  Observe that:
\begin{enumerate}
\item Every left-sided element commutes with $x_2$. In particular, $h_r^{-1}$ and $x_2$ commute.
\item $\ell (x_2 f) = \ell(f) + 3$ for any left-sided $f \in F$.
\end{enumerate}

Now, $h_l$ is the first vertex of $p$ whose right foot hits the critical leaf. Therefore,
when we get to $h_l$ in $p$, we have not yet modified any material to the right of the
critical leaf. In particular, there is some left-sided $h_l' \in F$ satisfying:
\begin{equation*}
h_l = x_2 h_l'
\end{equation*}
Observe that $\ell (h_l') = \ell (h_l) - 3 \leq 2n-1$, and hence $h_l'$ has width strictly
less than $n$.  Then the stacked diagram for $h_r^{-1} h_l'$ must already be reduced,
since no caret of $h_l'$ is far enough to the left to oppose the grounded caret of
$h_r^{-1}$.  From this, we conclude that $h_r^{-1} h_l'$ has width at least $n$.
Since $h_r^{-1} h_l'$ is strongly negative, lemma 6.3.4 implies that:
\begin{equation*}
\ell (h_r^{-1} h_l') \geq 2n
\end{equation*}
and hence:
\begin{equation*}
d(h_l,h_r) = \ell \left(h_r^{-1} h_l \right)
= \ell \left( x_2 h_r^{-1} h_l' \right) = \ell \left(h_r^{-1} h_l' \right) + 3
\geq 2n + 3\tag*{\qedhere}
\end{equation*}
\end{proof}

\section{Consequences}

Theorem 6.1.4 has some interesting consequences for the Cayley graph of $F$.
The first holds in any group that is not minimally almost convex, but we
state and prove it for $F$:

\begin{corollary} The Cayley graph for $F$ with respect to $\{x_0,x_1\}$ contains isometrically
embedded loops of arbitrary large circumference.
\end{corollary}
\begin{proof}
Fix an even $n \geq 4$, and let $l$ and $r$ be the two elements from Theorem 6.1.4.
Choose geodesics $p_l$ and $p_r$ connecting $l$ and $r$ to the identity vertex, and extend
these arcs to a closed loop using a path of length two from $l$ to $r$. We claim that this
loop $\gamma$ (of length $2n + 2$) is isometrically embedded.

Let $x\in p_l$ and $y\in p_r$, and suppose there were a
path $q$ from $x$ to $y$ shorter than both arcs connecting $x$ and $y$
inside $\gamma$.
\begin{center}\includegraphics{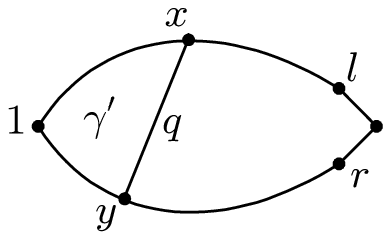}\end{center}
Then the loop $\gamma '$ (indicated in the picture) is shorter than
$\gamma$, and hence lies entirely inside the $n$-ball. In particular, the arc $q$ lies entirely
inside $B^n(F)$, a contradiction since this provides a shortcut from $l$ to $r$.

Since the two arcs in $\gamma$ from the identity to its antipode are also geodesics,
the loop $\gamma$ is isometrically embedded.
\end{proof}

Next, recall that a \emph{combing} of a group $G$ is, for each $g\in G$, a choice of a path in the Cayley graph
from the identity to $g$ (see section 1.5).  A \emph{geodesic combing} is a combing whose paths are geodesic segments.

The following result holds in any group that is not almost convex, but again we state and prove it just for $F$:

\begin{corollary}No geodesic combing of $F$ (with respect to $\{x_0,x_1\}$) has the fellow traveller property.
\end{corollary}
\begin{proof} Suppose we are given a geodesic combing of $F$.  Given any even $n\geq 4$, consider
the elements $l$ and $r$ from Theorem 6.1.4. Since the combing paths $p_l$ and $p_r$ are
geodesics, they can be closed to form an isometrically embedded loop as in Corollary 6.4.1.
In particular, the midpoint of $p_l$ has distance $n/2$ from the path~$p_r$.
\end{proof}

S.~Cleary and J.~Taback \cite{ClTa2} have also obtained Corollary 6.4.2.

\

\chapter{Strand Diagrams}

We already possess one algorithm for multiplying elements of $F$ using tree
diagrams (see example 1.2.6). In this chapter, we develop a much simpler,
more geometric understanding of multiplication, in the form of \emph{strand
diagrams}. These diagrams are closely related to a description of $F$ as the
fundamental group of the groupoid of fractions of a certain category of
finite binary forests.

The strand diagrams introduced in this chapter are ``dual'' to the diagrams
of Guba and Sapir (see \cite{GuSa}).  Matt Brin uses strand diagrams in
\cite{Brin} to represent elements of the braided Thompson group $BV$.

\section{Strand Diagrams}

The material in this section will be relatively informal, and many of the
proofs will be omitted. In section 7.2 we will develop a rigorous,
algebraic viewpoint towards many of the ideas introduced here.

A \emph{strand diagram} is any picture of the form:
\begin{center}
\includegraphics{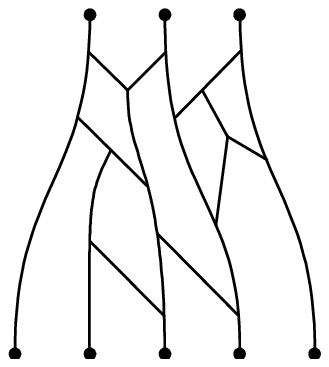}
\end{center}
A strand diagram is similar to a braid, except that instead of twists, there
are \emph{splits} and \emph{merges}:
\begin{center}
\includegraphics{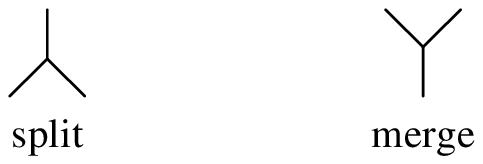}
\end{center}
Because of these splits and merges, a strand diagram may begin and end with
different numbers of strands.

\begin{notes} \
\begin{enumerate}
\item  As with a braid, the strands of a strand diagram are required
to have nonzero slope at all times. Hence, each strand has an ``up''
direction and a ``down'' direction.
\item  Isotopic strand diagrams are considered equal; that is, a
``strand diagram'' is really an isotopy class of strand diagram pictures.
The starting points and endpoints of the strand diagram are allowed to move
horizontally during these isotopies.
\end{enumerate}
\end{notes}

A \emph{reduction} of a strand diagram is one of the following two types of
moves:
\begin{center}
\includegraphics{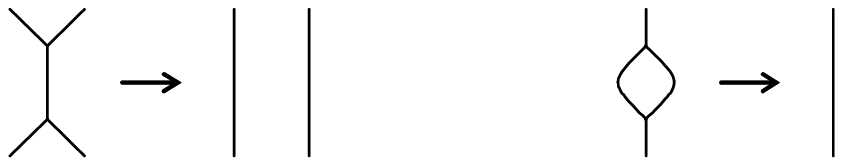}
\end{center}
Two strand diagrams are \emph{equivalent} if one can be transformed into the other
using reductions and inverse reductions. A strand diagram is \emph{reduced} if it
is not subject to any reductions.

\begin{proposition}  Every strand diagram is equivalent to a unique
reduced strand diagram.
\end{proposition}

\begin{notation}  If $i,j\geq 1$, the notation $f\colon i\rightarrow j$
will mean ``$f$ is a strand diagram that starts with $i$ strands and ends
with $j$ strands.''
\end{notation}

Given strand diagrams $f\colon i\rightarrow j$ and $g\colon
j\rightarrow k$, the \emph{concatenation} \mbox{$f\cdot g\colon i\rightarrow k$} is
obtained by attaching $g$ to the bottom of $f$. For example, if $f$ and $g$
are the strand diagrams:
\begin{center}
\includegraphics{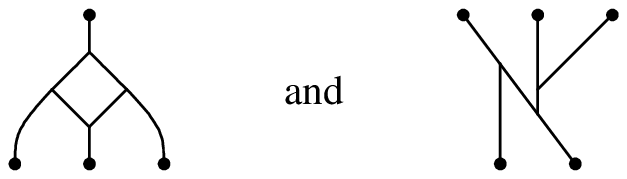}
\end{center}
then $f\cdot g$ is the strand diagram:
\begin{center}
\includegraphics{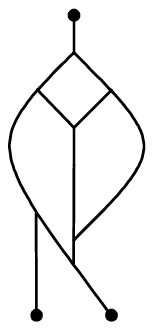}
\end{center}

\begin{proposition} Concatenation of strand diagrams is well-defined with
respect to equivalence.
\end{proposition}

If $f\colon i\rightarrow j$ and $g\colon j\rightarrow k$ are reduced
strand diagrams, define the \emph{composition} $fg\colon i\rightarrow k$ to be the
reduced strand diagram equivalent to $f\cdot g$.

\begin{proposition}  The collection of reduced strand diagrams forms a
groupoid under composition (with one object for each positive integer).
\end{proposition}

The identity morphism on $n$ is just the trivial strand diagram with
$n$ strands. Inverses are obtained by reflection across a horizontal line:
\begin{center}
\includegraphics{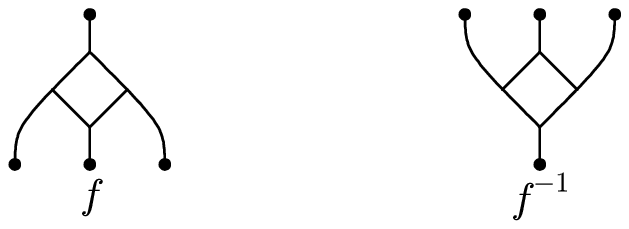}
\end{center}

The following theorem explains our interest in strand diagrams:

\begin{theorem}  The fundamental group of the groupoid of strand diagrams
is Thompson's group $F$.

That is, given any positive integer $n$, the group of all reduced
strand diagrams that begin and end with $n$ strands is isomorphic with $F$.
\end{theorem}

\begin{proof}[Informal Proof]  Define a \emph{forest} to be any reduced strand diagram that has
no merges. Observe that each forest is essentially just a finite sequence
of binary trees. We claim that every reduced strand diagram is the
concatenation of a forest and an inverse forest.

Given any picture of a reduced strand diagram, we can draw a curve with
the following properties:
\begin{enumerate}
\item  Every path from the top to the bottom of the strand diagram crosses
the curve exactly once.
\item  Every split lies above the curve, and every merge lies below the
curve.
\end{enumerate}
\begin{center}
\includegraphics{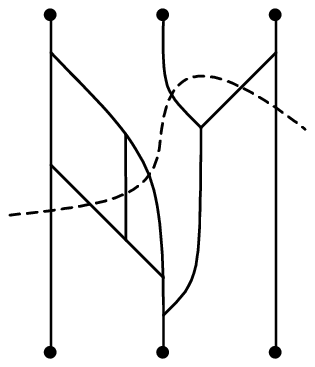}
\end{center}
Cutting along this curve gives the desired decomposition.

In particular, any reduced strand diagram that begins and end with one strand is the
concatenation of a tree and an inverse tree:
\begin{center}
\includegraphics{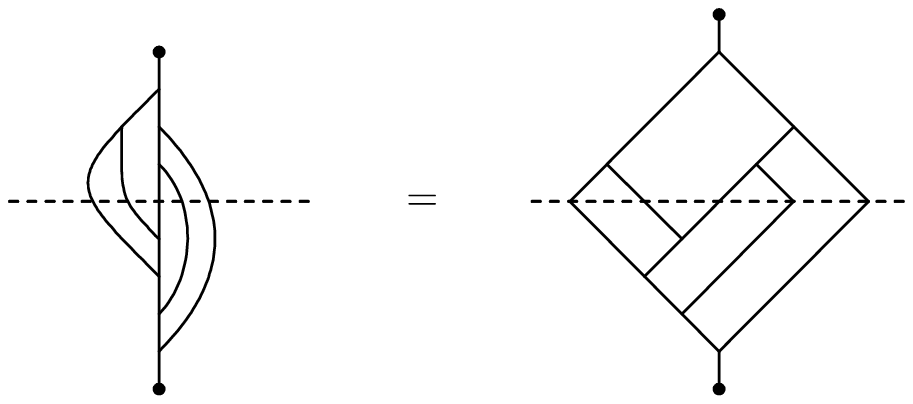}
\end{center}
This is the tree diagram for the corresponding element of
$F$.\end{proof}

Because of this theorem, we will refer to the groupoid of strand
diagrams as \emph{Thompson's groupoid $\mathcal{F}$}.

Observe that we have skirted the issue of whether composition of
reduced strand diagrams in fact corresponds to multiplication in $F$. We
will prove this in an algebraic context in the following section.

\

\section{Thompson's Groupoid}

In this section we define Thompson's groupoid algebraically and show that its
fundamental group is isomorphic with Thompson's group $F$.

First we define the \emph{category of forests $\mathcal{P}$}:
\begin{flushleft}
\textbf{Objects:} There is one object of $\mathcal{P}$ for each positive
integer.
\end{flushleft}
\begin{flushleft}
\textbf{Morphisms:}  A morphism $i\rightarrow j$ is a finite binary forest with
$i$ trees and $j$ total leaves.
\end{flushleft}
\begin{flushleft}
\textbf{Composition:}  If $f\colon i\rightarrow j$ and $g\colon j\rightarrow
k$, the composition $fg\colon i\rightarrow k$ is obtained by attaching the
roots of the trees of $g$ to the leaves of $f$ in an order-preserving way.
\end{flushleft}

Note that we are composing elements of $\mathcal{P}$ using the same ``backwards''
convention that we have been using throughout.

For $0\leq n<w$, let $x_n\colon w\rightarrow w+1$ be the forest:
\begin{center}
\includegraphics{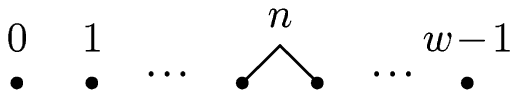}
\end{center}
Then any forest is a product of the $x_n$'s since any forest can be
obtained from a trivial forest by attaching carets. If we attach carets
from left to right, we get the \emph{normal form}:

\begin{proposition}  Every nontrivial morphism of $\mathcal{P}$ can be
expressed uniquely as:
\begin{equation*}
x_0^{a_0}x_1^{a_1}\cdots x_n^{a_n}
\end{equation*}
where $a_0,\ldots,a_n\in\mathbb{N}$ and $a_n\neq 0$.\quad\qedsymbol
\end{proposition}

\begin{corollary} The category $\mathcal{P}$ has the following
presentation:
\begin{enumerate}
\item[]\textbf{Generators:} One generator $x_n\colon w\rightarrow w+1$ for each
$0\leq n<w$.
\item[]\textbf{Relations:} One relation:
\begin{equation*}
w\overset{x_n}{\longrightarrow}w+1\overset{x_k}{\longrightarrow}w+2\quad=%
\quad w\overset{x_k}{\longrightarrow}w+1\overset{x_{n+1}}{\longrightarrow}w+2
\end{equation*}
for each $0\leq k<n<w$.
\end{enumerate}
\end{corollary}

We wish to construct $\mathcal{F}$ as the groupoid of fractions for the category
$\mathcal{P}$. We begin with a brief general discussion concerning groupoids of fractions.
See \cite{ClPr} for proofs of the statements below in the context of semigroups and groups of fractions.

\begin{definition} Let $\mathcal{C}$ be any category. A \emph{groupoid of right
fractions} for $\mathcal{C}$ is a groupoid $\mathcal{G}$ containing
$\mathcal{C}$, and having the following properties:
\begin{enumerate}
\item Every object of $\mathcal{G}$ is an object of $\mathcal{C}$.
\item Every morphism of $\mathcal{G}$ can be expressed as
$pq^{-1}$, where $p$ and $q$ are morphisms of $\mathcal{C}$.
\end{enumerate}
\end{definition}

An expression of the form $pq^{-1}$, where $p$ and $q$ are morphisms in $\mathcal{C}$, is
called a \emph{right fraction}.

\begin{proposition}
Let $\mathcal{C}$ be a category with groupoid of right fractions $\mathcal{G}$.  If
$p_1q_1^{-1}$ and $p_2q_2^{-1}$ are right fractions, then $p_1q_1^{-1} = p_2q_2^{-1}$
if and only if there exist morphisms $r_1,r_2$ in $\mathcal{C}$ making the following
diagram commute:
\begin{equation*}
\def\labelstyle{\textstyle}
\xymatrix{& \bullet \ar[dl]_{p_1} \ar[dr]^{p_2} & \\
\bullet \ar[r]^{r_1} & \bullet & \bullet \ar[l]_{r_2} \\
& \bullet \ar[ul]^{q_1} \ar[ur]_{q_2} &}
\end{equation*}
\end{proposition}

Based on this proposition, it is clear that any two groupoids of fractions for a given
category are isomorphic.  The following theorem gives necessary and sufficient conditions
for a category to have a groupoid of fractions:

\begin{theorem}  Let $\mathcal{C}$ be any category. Then $\mathcal{C}$
has a groupoid of fractions if and only if $\mathcal{C}$ has the following
properties:
\begin{enumerate}
\item ($\mathcal{C}$ is cancellative) For any morphisms such that the stated compositions exist:
\begin{equation*}
pr=qr\quad\Rightarrow\quad p=q\qquad and\qquad lp=lq\quad\Rightarrow\quad p=q
\end{equation*}
\item  ($\mathcal{C}$ has common right multiples) Given any morphisms
$p,q$ with the same domain, there exist morphisms $r,s$ such that $pr=qs$.\quad\qedsymbol
\end{enumerate}
\end{theorem}

\begin{theorem}  The category $\mathcal{P}$ has a groupoid of right fractions $\mathcal{F}$.
\end{theorem}
\begin{proof}  $\mathcal{P}$ is clearly cancellative. Next, suppose that $f$
and $g$ are any two morphisms with the same domain $w$ (so $f$ and $g$ are
forests with $w$ trees). Let $n$ be the maximum height of all of the trees
in $f$ and $g$. Then $f$ and $g$ have as a common right multiple the forest
with $w$ complete binary trees of height $n$.\end{proof}

\begin{remark}
A strand diagram picture is just a word in the $x_n$'s and
$x_n^{-1}$'s. In particular, suppose we have a strand diagram picture with
the property that all merges and splits occur at different heights. Then
each split corresponds to an instance of some $x_n$, and each merge
corresponds to an instance of some $x_n^{-1}$.

If we perform an isotopy on a strand diagram that causes the heights of
two intersections to switch, it corresponds to an application of one of the
following types of relations:
\begin{equation*}
\begin{aligned}[t]
&&\quad x_nx_k&=x_kx_{n+1}\\
&&\quad x_n^{-1}x_k&=x_kx_{n+1}^{-1}\\
&&\quad x_k^{-1}x_n&=x_{n+1}x_k^{-1}\\
&\text{or}&\quad x_k^{-1}x_n^{-1}&=x_{n+1}^{-1}x_k^{-1}
\end{aligned}
\end{equation*}
The two reductions:
\begin{center}
\includegraphics{npict70}
\end{center}
correspond to cancelling an $x_n^{-1}x_n$ or $x_nx_n^{-1}$ pair,
respectively.

This explains why the groupoid $\mathcal{F}$ constructed above is the
same as the groupoid of strand diagrams defined in section 7.1.
\end{remark}

\begin{theorem}  The fundamental group of $\mathcal{F}$ is isomorphic
with Thompson's group $F$.
\end{theorem}

\begin{proof}Let $\mathcal{PL}$ be the groupoid of closed intervals and piecewise-linear
homeomorphisms.  Define a homomorphism (functor) $\rho\colon\mathcal{P}\rightarrow\mathcal{PL}$
as follows:
\begin{enumerate}
\item $\rho(w)=[0,w]$ for any positive integer $w$.
\item If $x_n\colon w\rightarrow w+1$, then
$\rho(x_n)\colon [0,w]\rightarrow [0,w+1]$ is
the homeomorphism with slope $1$ on $[0,n]\cup [n+1,w]$ and slope $2$ on $[n,n+1]$.
\end{enumerate}
It is easy to verify that $\rho$ respects the relations in $\mathcal{P}$, and is therefore
a well-defined homomorphism.  Observe also that $\rho$ is one-to-one on morphisms.

Note that, if $f:1\rightarrow w$, then $\rho(f)$ is a homeomorphism that sends the
intervals of some dyadic subdivision of $[0,1]$ linearly onto the intervals
$[0,1],\ldots,[w-1,w]$.

Since $\mathcal{PL}$ is a groupoid, the monomorphism
$\rho\colon\mathcal{P}\rightarrow\mathcal{PL}$ extends to a monomorphism
$\rho\colon\mathcal{F}\rightarrow\mathcal{PL}$.  If $f\colon 1 \rightarrow 1$
is any morphism of $\mathcal{F}$, then $f=pq^{-1}$ for some morphisms
$p,q\colon 1 \rightarrow w$ of $\mathcal{P}$, and therefore $\rho(f)$ is
the homeomorphism $[0,1] \rightarrow [0,1]$ that sends the intervals of the
dyadic subdivision for $p$ linearly to the intervals of the dyadic subdivision
for $q$.
We conclude that the image under $\rho$ of $\pi_1(F,1)$ is precisely the group $F$.
\end{proof}

Using an argument similar to the proof of theorem 1.1.2, one can show that the image
of $\mathcal(F)$ under $\rho$ is precisely the set of piecewise-linear homeomorphisms
$f\colon [0,i] \rightarrow [0,j]$ such that:
\begin{enumerate}
\item All slopes of $f$ are powers of $2$, and
\item All breakpoints of $f$ have dyadic rational coordinates
\end{enumerate}

Theorem 7.2.8 yields an alternate derivation of the standard presentation for~$F$:

\begin{theorem}Thompson's group $F$ has presentation:
\begin{equation*}
\langle x_0,x_1,x_2,\ldots \mid x_n x_k = x_k x_{n+1}\text{ for }k<n \rangle
\end{equation*}
\end{theorem}
\begin{proof}Since $\mathcal{F}$ is the groupoid of fractions for
$\mathcal{P}$, the presentation for $\mathcal{F}$ is the same as the
presentation for $\mathcal{C}$. Therefore, $\mathcal{F}$ is generated by
elements $x_n^{(w)}\colon w\rightarrow w+1$ with relations $x_n^{(w)}x_k^{(w+1)}=x_k^{(w)}x_{n+1}^{(w+1)}$
($0\leq k<n<w$):
\begin{center}
\includegraphics{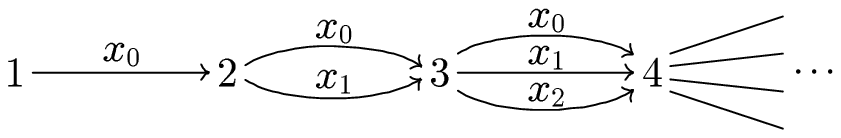}
\end{center}

To find a presentation for $\pi_1\left(\mathcal{F},1\right)$, we must
choose a spanning subtree of the graph of generators to contract. We choose the subtree
$\bigl\{ x_0^{(1)}, x_1^{(2)}, x_2^{(3)}, \ldots \bigr\}$:
\begin{center}
\includegraphics{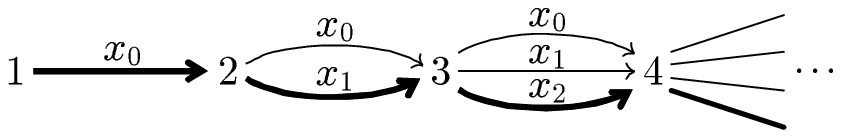}
\end{center}
Therefore, $\pi_1(\mathcal{F},1)$ is generated by elements $x_n^{(w)}$ ($n<w-1$).
The relations $x_n^{(w)}x_k^{(w+1)}=x_k^{(w)}x_{n+1}^{(w+1)}$ become:
\begin{equation*}
x_k^{(w+1)} = x_k^{(w)}
\end{equation*}
when $n=w-1$ and:
\begin{equation*}
x_n^{(w)}x_k^{(w+1)}=x_k^{(w)}x_{n+1}^{(w+1)}
\end{equation*}
for $n<w-1$.  These reduce immediately to the standard relations for $F$, where $x_n$ is the element
\begin{equation*}
x_n^{(n+2)} = x_n^{(n+3)} = x_n^{(n+4)} = \cdots\tag*{\qedhere}
\end{equation*}
\end{proof}

\

\section{Strand Diagrams and Braids}

The similarity between strand diagrams and braids is more than superficial:  it is possible
to develop some of the theory of $F$ in analogy with the development of the theory of
braid groups.  In this section, we describe a classifying space for $F$ that is
analogous to the standard classifying spaces for the braid groups.

A braid is essentially just the path of motion of $n$ points in the plane, i.e. a loop
in the configuration space on $n$ points in $\mathbb{R}^2$.  The following is well-known:

\begin{theorem}Let $B_n$ denote the braid group on $n$ strands, and let $X_n$ be the
configuration space of $n$ points in $\mathbb{R}^2$.  Then $X_n$ is a classifying space
for $\mathbb{R}^n$.
\end{theorem}
\begin{proof}See \cite{FaNe}.\end{proof}

A strand diagram represents the motion of finitely many points on the real line, with
the points allowed to split and merge in pairs.  We wish to construct the corresponding
``configuration space''.

Let $X_w$ be the collection of all $w$-tuples $(t_0,t_1,\ldots,t_{w-1})$ satisfying:
\begin{enumerate}
\item $t_0 \leq t_1 \leq \cdots \leq t_{w-1}$, and
\item $t_{i+2}-t_i \geq 1$ for all $i$.
\end{enumerate}
(The purpose of the second condition is to prevent three points from merging simultaneously.)
Let $X$ be the disjoint union of the $X_w$'s, subject to the identifications:
\begin{equation*}
(t_0,t_1,\ldots,t_{w-1}) \equiv (t_0,t_1,\ldots,t_n,t_n,\ldots,t_{w-1})
\end{equation*}

\begin{theorem} $X$ is a classifying space for $F$.
\end{theorem}

We will sketch of proof of this theorem for the remainder of this section.  Many tedious topological
details will be omitted.

For $i,j > 0$, let $[i\rightarrow j]$ denote all morphisms in Thompson's
groupoid $\mathcal{F}$ from $i$ to $j$.
For each $w$, let $\widetilde{X}_w = X_w \times [1\rightarrow w]$, where the set
$[1\rightarrow w]$ has the discrete topology.  Let $\widetilde{X}$ be the disjoint union
of the $\widetilde{X}_w$'s, subject to the identifications:
\begin{equation*}
((t_0,t_1,\ldots,t_{w-1}),f) \equiv ((t_0,t_1,\ldots,t_n,t_n,\ldots,t_{w-1}),fx_n)
\end{equation*}
There is an obvious left-action of $F$ on $\widetilde{X}$, with quotient $X$.  It is not
hard to see that this is a covering space action, so that $\widetilde{X}$ is a covering
space of $X$.  We claim that $X$ is contractible.

The plan is to exhibit an explicit contraction of the space $\widetilde{X}$.  Observe
that an element $x$ of $\widetilde{X}$ is essentially just a strand diagram that starts
with $1$ strand and ends with $w$ strands, together with specified positions for the
endpoints.  The idea is to choose our ``favorite picture'' $D$ of this strand diagram,
and then ``run the diagram backwards''.  That is, assuming $D$ has height $1$,
we will move $x$ along the path which at time $t$ is represented by the initial
segment of $D$ with height $1-t$.  The trick is to find a way of choosing
our ``favorite picture'' that varies continuously with position in $\widetilde{X}$.

However, we would first like to simplify the situation.  Let $Y$ be the subspace of
$X$ consisting of all $w$-tuples with first coordinate $0$.  To specify an element
of $Y$, we need only specify the distances between the $w$ strands:
\begin{equation*}
(0,t_1,\ldots,t_{w-1}) = [t_1,t_2-t_1,\ldots,t_{w-1}-t_{w-2}]
\end{equation*}
Notice that a tuple $[d_1,\ldots,d_{w-1}]$ specifies an element of $Y$ if and only if
each $d_n \geq 0$ and $d_n+d_{n+1} \geq 1$ for all $n$.
Let $\widetilde{Y}$ be the subspace of $\widetilde{X}$ that maps onto $Y$.
Then $\widetilde{X}$ clearly deformation retracts onto
$\widetilde{Y}$.  We will exhibit an explicit contraction of the space $\widetilde{Y}$.

Now some terminology:

\begin{definition} Suppose that $f\colon 1 \rightarrow w$.
\begin{enumerate}
\item We say that the $n$'th strand of $f$ has \emph{just merged}
if right-multiplication by $x_n$ would cancel a merge
in the strand diagram for $f$.
\item We say that strands $n$ and $n+1$ have \emph{just split}
if right-multiplication by $x_n^{-1}$ would cancel a split in the strand diagram for $f$.
\end{enumerate}
\end{definition}

For example, if $f$ is the element:
\begin{center}
\includegraphics{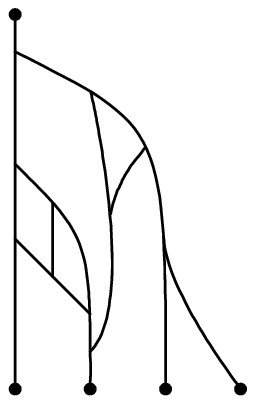}
\end{center}
then strand $1$ has just merged, and strands $2$ and $3$ have just split.

We are now ready to describe the contraction of $\widetilde{Y}$.  Suppose that
$y \in \widetilde{Y}$, with distances $[d_1,\ldots,d_{w-1}]$ and strand diagram
$f\colon 1 \rightarrow w$.  Then $y$ moves as follows:
\begin{enumerate}
\item If strand $n$ has just merged in $f$, and $d_{n-1},d_n \geq 1$, then strand
    $n$ immediately splits. That is, a $0$ is inserted between $d_{n-1}$ and $d_n$,
    and this $0$ begins increasing at unit speed.
\item If strands $n$ and $n+1$ have just split, then the distance $d_n$ decreases
    at unit speed until it reaches $0$, at which point $d_n$ is removed.
\item Otherwise, the distance $d_n$ moves toward $1$ at unit speed.
\end{enumerate}
Note that the point $y$ might take arbitrarily long to reach the basepoint of
$\widetilde{Y}$.  Therefore, the described contraction takes place during the time
interval $[0,\infty]$.\qed

\begin{example}
Let $y$ be the point $([1,0.8,1,0.6],x_0^3 x_2 x_4 x_3^{-1})$:
\begin{center}
\includegraphics{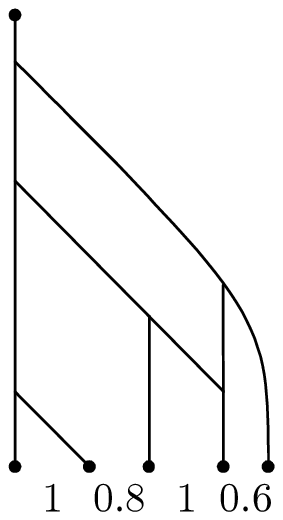}
\end{center}
Then $y$ moves towards the basepoint of $\widetilde{Y}$ as follows:
\begin{center}
\begin{tabular}{c|c|c|c}
Time & Position in $Y$ & Position in $\mathcal{F}$ & Movement \\ \hline
$0$ & $[1,.8,1,.6]$ & $x_0^3 x_2 x_4 x_3^{-1}$ &
$d_1 \rightarrow 0, d_2 \rightarrow 1, d_4\rightarrow 1$ \\
$0.2$ & $[.8,1,1,.8]$ & $x_0^3 x_2 x_4 x_3^{-1}$ &
$d_1 \rightarrow 0, d_4 \rightarrow 1$ \\
$0.4$ & $[.6,1,1,1]$ & $x_0^3 x_2 x_4 x_3^{-1}$ &
Strand $3$ splits.\\
$0.4$ & $[.6,1,1,0,1]$ & $x_0^3 x_2 x_4$ &
$d_1 \rightarrow 0, d_3 \rightarrow 0, d_4 \rightarrow 1, d_5 \rightarrow 0$ \\
$1$ & $[0,1,.4,.6,.4]$ & $x_0^3 x_2 x_4$ &
Strands $0$ and $1$ merge.\\
$1$ & $[1,.4,.6,.4]$ & $x_0^2 x_1 x_3$ &
$d_2 \rightarrow 0, d_3 \rightarrow 1, d_4 \rightarrow 0$ \\
$1.4$ & $[1,0,1,0]$ & $x_0^2 x_1 x_3$ &
Strands $1,2$ and $3,4$ merge.\\
$1.4$ & $[1,1]$ & $x_0^2$ &
$d_1 \rightarrow 0$\\
$2.4$ & $[0,1]$ & $x_0^2$ &
Strands $0$ and $1$ merge.\\
$2.4$ & $[1]$ & $x_0$ &
$d_1 \rightarrow 0$\\
$3.4$ & $[0]$ & $x_0$ &
Strands $0$ and $1$ merge\\
$3.4$ & $[]$ & identity &
Basepoint reached.\\
\end{tabular}
\end{center}

The path followed by this element can be summarized by the following diagram:
\begin{center}
\includegraphics{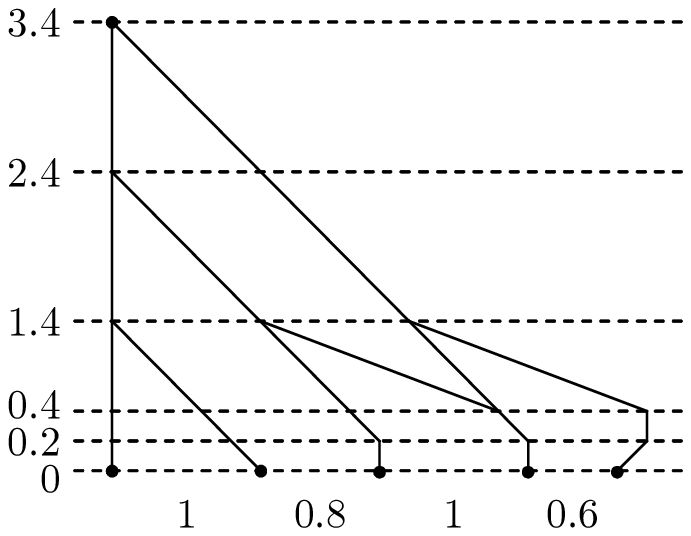}
\end{center}
\end{example}

\section{Other Thompson Groups}

There are several groups similar to $F$ that also have strand diagrams.  In this section, we will
briefly introduce each of these groups and discuss the corresponding strand diagrams, as well as
algebraic constructions of the corresponding groupoids.  As an application, we shall compute an infinite
presentation for each of these groups.  Because they are based on strand diagrams, these presentations
all admit a ``normal form'' similar to the normal form for elements of $F$.

\subsubsection{The Groups $T$ and $\widetilde{T}$}

We shall begin by discussing \emph{Thompson's Group T}, which is
a ``circular'' version of Thompson's group $F$.

Let $\pi \colon [0,1] \rightarrow S^1$ be the quotient map.  A \emph{dyadic subdivision}
of $S^1$ is any image under $\pi$ of a dyadic subdivision of $[0,1]$.  A \emph{dyadic rearrangement}
of $S^1$ is any homeomorphism $S^1 \rightarrow S^1$ that maps the intervals of one dyadic
subdivision linearly onto the intervals of another, preserving the \emph{cyclic} order
of the intervals.  For example, if $\mathcal{D}$ and $\mathcal{R}$ are the subdivisions:
\begin{center}
\includegraphics{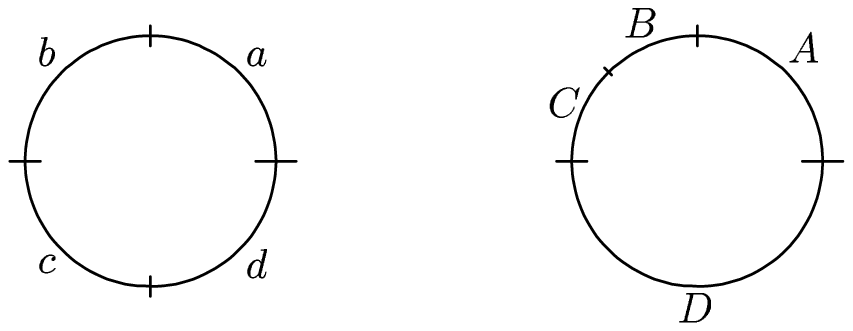}
\end{center}
then there exists a dyadic rearrangement of $S^1$ that sends the intervals $(a,b,c,d)$ linearly onto
the intervals $(B,C,D,A)$.

    The set $T$ of all dyadic rearrangements of $S^1$ forms a group under composition.  It is isomorphic
to the group of cyclic-order preserving automorphisms of a free Cantor algebra
(see section 1.6 for a definition of Cantor algebras, and see \cite{Bro} for details).
The group $T$ was introduced by Thompson, who proved that $T$ is finitely presented and simple.
(See \cite{CFP} for a published version of these results, and a thorough introduction to $T$.)
Like $F$, the group $T$ has type $\rm{F}_\infty$ (see \cite{Bro}).

    We can represent any element of $T$ by a pair of binary trees,  together with a cyclic
permutation of the leaves.  For example, the element above can be represented by the diagram:
\begin{center}
\includegraphics{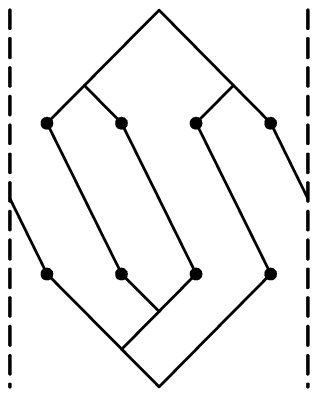}
\end{center}
This is called a \emph{tree diagram} for an element of $T$.  (It can be helpful to think of a
tree diagram as being embedded on the cylinder, with the dashed lines identified.)
The tree diagram for the above element is not reduced:
\begin{center}
\includegraphics{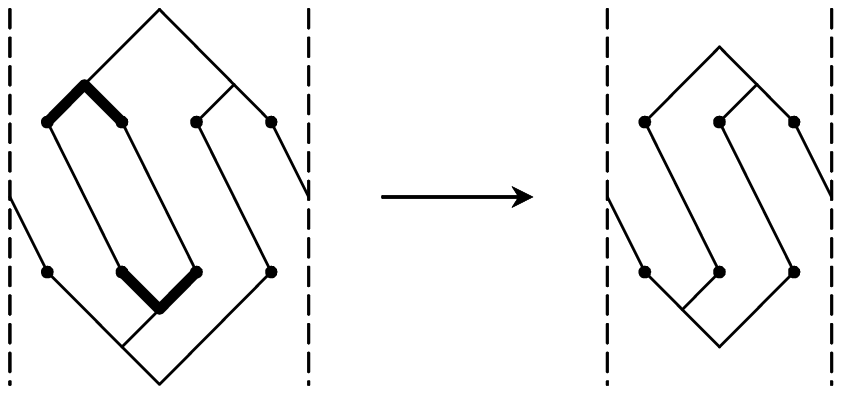}
\end{center}

A \emph{cylindrical strand diagram} is any strand diagram that is embedded on the cylinder:
\begin{center}
\includegraphics{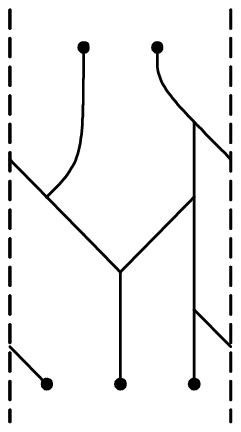}
\end{center}
Two cylindrical strand diagrams are \emph{equivalent} if one can be obtained from the other
by a sequence of:
\begin{enumerate}
\item reductions, inverse reductions, and
\item Dehn twists of the cylinder.
\end{enumerate}
(Allowing Dehn twists is necessary because a rotation of the circle
by $2\pi$ is equal to the identity in $T$.  Equivalently, we could allow the begin-points
and endpoints to move horizontally around the circle during isotopies.) The group of all
equivalence classes of cylindrical strand diagrams that start and end with one strand is
isomorphic with Thompson's group $T$.

A cylindrical strand diagram is really just a word for an element of a certain groupoid.
In particular, let $\mathcal{P}[\mathbb{Z}_w]$ be the category of ``forests plus cyclic
permutations''  obtained from $\mathcal{P}$ by attaching a copy of $\mathbb{Z}_w$ at each vertex $w$:
\begin{center}
\includegraphics{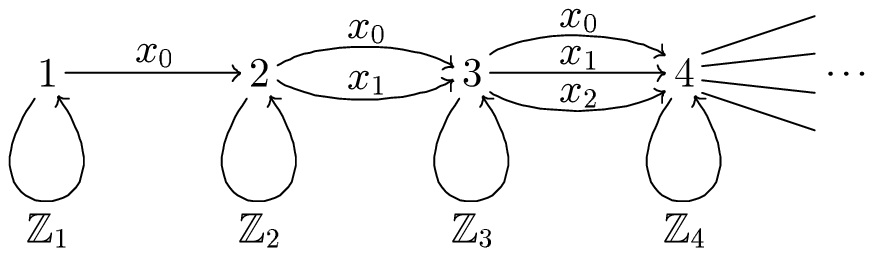}
\end{center}
Then $\mathcal{P}[\mathbb{Z}_w]$ is generated by the morphisms $x_n\colon w \rightarrow w+1$ ($n<w$) from $\mathcal{P}$
together with one morphism $\omega_w \colon w \rightarrow w$ for each $w \geq 2$ satisfying
the relations:
\begin{equation*}
\begin{array}{cl}
\omega_w ^ w = 1 & \\
\omega_w x_n = x_{n+1} \omega_{w+1} & \text{($n<w-1$)} \\
\omega_w x_{w-1} = x_0 \omega_{w+1} ^2 &
\end{array}
\end{equation*}
For example:
\begin{center}
\includegraphics{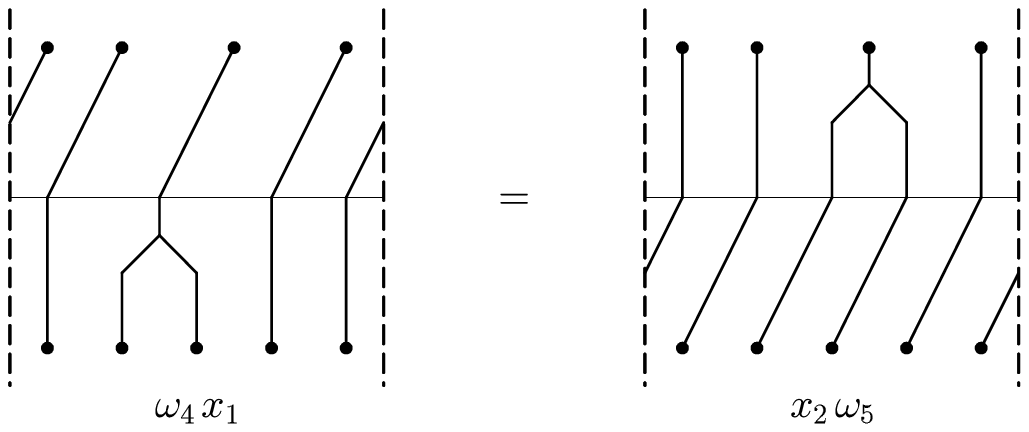}
\end{center}
and:
\begin{center}
\includegraphics{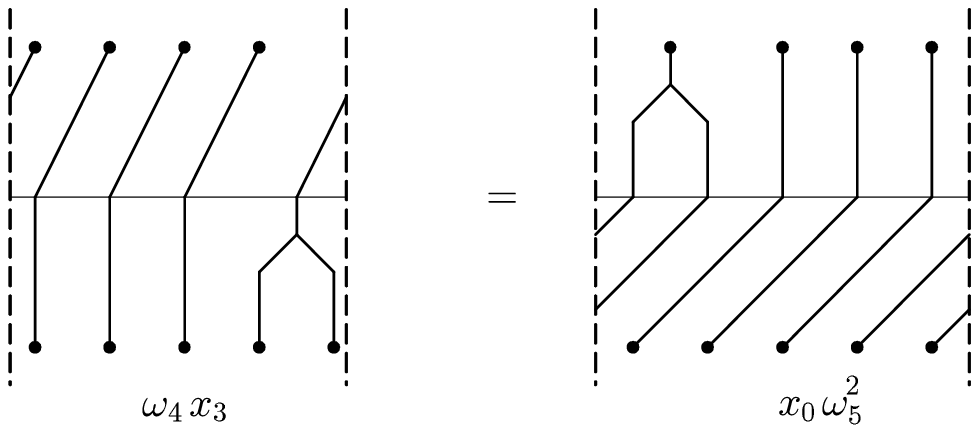}
\end{center}

It is not hard to show that $\mathcal{P}[\mathbb{Z}_n]$ is cancellative and has common right multiples, so by
theorem 7.2.5 $\mathcal{P}[\mathbb{Z}_n]$ has a groupoid of right fractions $\mathcal{T}$.

\begin{proposition}
The groupoid $\mathcal{T}$ has fundamental group $T$.\quad\qedsymbol
\end{proposition}

We can use this to calculate a presentation for $T$:

\begin{theorem}
The group $T$ is generated by elements $\{x_0,x_1,x_2,\ldots\}$ and \linebreak
$\{\omega_2,\omega_3,\omega_4,\ldots\}$, with relations
\begin{equation*}
\begin{array}{l}
\omega_n^n = 1 \\
 x_n x_k = x_k x_{n+1}\quad\text{for $n>k$}
\end{array}
\end{equation*}
and:
\begin{equation*}
\begin{aligned}
\omega_n x_k & = x_{k+1} \omega_{n+1}\quad\text{for $k<n-2$} \\
\omega_n x_{n-2} &=\omega_{n+1} \\
\omega_n &= x_0 \omega_{n+1}^2
\end{aligned}
\end{equation*}
\end{theorem}
\begin{proof}Since $\mathcal{T}$ is the groupoid of fractions of $\mathcal{P}[\mathbb{Z}_w]$, the
presentation for $\mathcal{T}$ is the same as the presentation for $\mathcal{P}[\mathbb{Z}_w]$.
Therefore, $\mathcal{T}$ is generated by the elements $x_n^{(w)}\colon w \rightarrow w+1$ and
$\omega_w \colon w \rightarrow w$ with relations:
\begin{equation*}
\begin{array}{cl}
x_n^{(w)} x_k ^{(w+1)} = x_k^{(w)} x_{n+1}^{(w+1)} \\
\omega_w ^ w = 1 & \\
\omega_w x_n^{(w)} = x_{n+1}^{(w)} \omega_{w+1} & \text{($n<w-1$)} \\
\omega_w x_{w-1}^{(w)} = x_0^{(w)} \omega_{w+1} ^2 &
\end{array}
\end{equation*}
To find a presentation for $\pi_1(\mathcal{T},1)$, we use the morphisms
$x_0^{(1)},x_1^{(2)},x_2^{(3)},\ldots$ as a spanning tree.  As in theorem 7.2.9, the
first family of relations implies that:
\begin{equation*}
x_n^{(n+2)}=x_n^{(n+3)}=x_n^{(n+4)}=\cdots
\end{equation*}
in $\pi_1(\mathcal{T},1)$ for each $n$.  If we label this element $x_n$, then the remainder
of the first family of relations reduces to:
\begin{equation*}
 x_n x_k = x_k x_{n+1}\quad\text{for $n>k$}
\end{equation*}
The third family of relations yields:
\begin{equation*}
\omega_w x_n = x_{n+1} \omega_{w+1} \quad \text{($n<w-2$)}
\end{equation*}
and:
\begin{equation*}
\omega_w x_{w-2} = \omega_{w+1}
\end{equation*}
in the case when $n=w-2$.  Finally, the fourth family of relations reduces to:
\begin{equation*}
\omega_w = x_0 \omega_{w+1} ^2\tag*{\qedhere}
\end{equation*}
\end{proof}

There is another group $\widetilde{T}$ that is similar to $T$ but torsion-free, namely
the lift of $T$ in the group of PL-homeomorphisms of $\mathbb{R}$.  It is the
fundamental group of the groupoid of right fractions of a category $\mathcal{P}[\mathbb{Z}]$
which can be constructed by attaching a copy of $\mathbb{Z}$ to each object of $\mathcal{P}$:
\begin{center}
\includegraphics{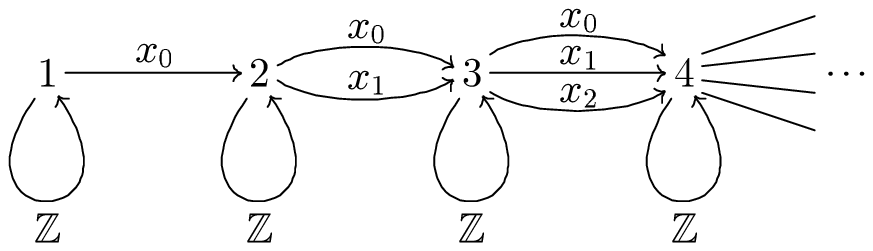}
\end{center}
The generators of the $\mathbb{Z}$'s are required to satisfy the same relations in
$\mathcal{P}[\mathbb{Z}]$ that the generators of the $\mathbb{Z}_w$'s satisfied in
$\mathcal{P}[\mathbb{Z}_w]$, excepting the relations $\omega_w^w=1$.
Elements of this groupoid can also be represented by cylindrical strand diagrams,
except that two diagrams that differ by a Dehn twist are \emph{not} equivalent.  Hence the
element:
\begin{center}
\includegraphics{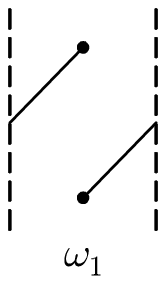}
\end{center}
is \emph{not} trivial in $\widetilde{T}$.

\begin{theorem} The group $\widetilde{T}$ is generated by the elements $\{x_0,x_1,x_2,\ldots\}$
and $\{\omega_1,\omega_2,\omega_3,\ldots\}$, with relations:
\begin{equation*}
\begin{array}{l}
x_n x_k = x_k x_{n+1} \quad \text{for $k<n$} \\
\omega_n x_k = x_{k+1} \omega_{n+1} \quad \text{for $k<n-2$} \\
\omega_n x_{n-2} = \omega_{n+1} \\
\omega_n = x_0 \omega_{n+1}^2 \quad \text{for $n>1$} \\
\omega_1 = \omega_2^2
\end{array}\tag*{\qedsymbol}
\end{equation*}
\end{theorem}

It is not hard to show that:
\begin{equation*}
\omega_1 = \omega_2^2 = \omega_3^3 = \cdots
\end{equation*}
in $\widetilde{T}$ and that this element is central.  (It corresponds to a ``rotation of the circle by an angle
of $2\pi$''.)  Therefore, the epimorphism $\widetilde{T}\twoheadrightarrow T$
has kernel $\mathbb{Z}$.

It seems likely that there is a classifying space for $\widetilde{T}$ similar to the one for $F$
constructed in the previous section, i.e. the ``configuration space'' of finitely many points
on a circle, with the points allowed to split and merge in pairs, but the details have yet to
be worked out.

\subsubsection{The Groups $V$ and $BV$}

An element of \emph{Thompson's Group $V$} is obtained by sending the intervals of some
dyadic subdivision of $[0,1]$ linearly onto the intervals of another, except that the
order of the intervals may be arbitrarily permuted.  Note that this produces bijections
$[0,1] \rightarrow [0,1]$ that are \emph{not} continuous.  (By convention, all functions
in $V$ are required to be continuous from the right.  Alternatively, one can define $V$
as a group of homeomorphisms of the Cantor set.)

The set of all elements of $V$ forms a group under composition --- it is the same as the
group called $V$ in section 1.6.  This group was introduced by Thompson along with $F$ and $V$.
He proved that $V$ is simple and finitely presented.  (See \cite{CFP} for a published version
of these results, and a thorough introduction to $V$.)  Like $F$ and $T$, the group $V$ has type
$\rm{F}_\infty$ (see \cite{Bro}).

Recall that an element of $V$ can be represented by a pair of binary forests, together with any
permutation of the leaves:
\begin{center}
\includegraphics{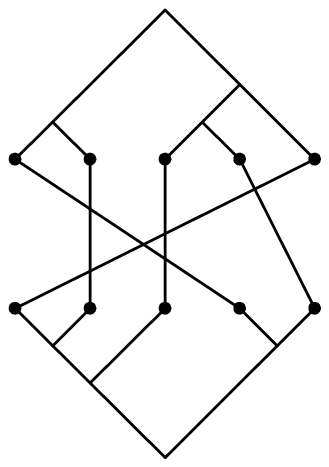}
\end{center}
This is called a \emph{tree diagram} for an element of $V$.

There is also a groupoid of strand diagrams corresponding to $V$.  An element of this groupoid is
a strand diagram with splits, merges, and crosses:
\begin{center}
\includegraphics{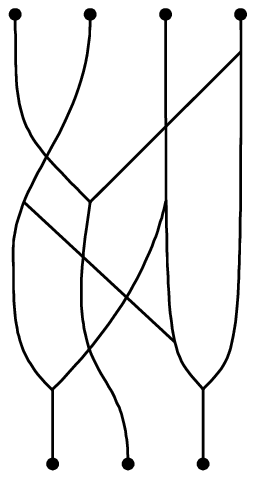}
\end{center}

Two V-strand diagrams are \emph{equivalent} if one can be obtained from the other
by a sequence of reductions, inverse reductions, and homotopies.
The group of all equivalence classes of these strand diagrams that start and end
with one strand is isomorphic with Thompson's group $V$.

It is possible to construct this groupoid algebraically as follows.  Let $\Sigma_w$ denote the
permutation group on ${0,1,\ldots,w-1}$, and let
$s_n \colon \Sigma_w \rightarrow \Sigma_{w+1}$ be the function that ``doubles''
whichever element maps to $n$.
For example, if $\sigma$ is the permutation:
\begin{center}
\includegraphics{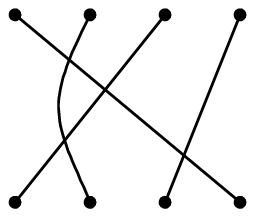}
\end{center}
then $s_0(\sigma)$ is the permutation:
\begin{center}
\includegraphics{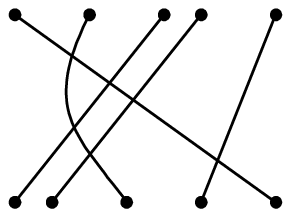}
\end{center}
Note that $s_n$ is not a homomorphism --- it is just a function from $\Sigma_w$ to
$\Sigma_{w+1}$.

Now let $\mathcal{P}[\Sigma_w]$ be the category of ``forests plus permutations''
obtained from $\mathcal{P}$ by attaching a copy of $\Sigma_w$ at each vertex $w$:
\begin{center}
\includegraphics{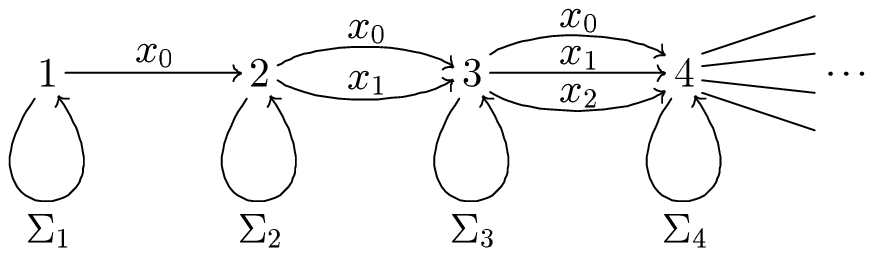}
\end{center}
with relations:
\begin{equation*}
\sigma x_n = x_{\sigma^{-1}(n)} s_n(\sigma)\qquad\text{($\sigma\in\Sigma_w$ and $n<w$)}
\end{equation*}
For example:
\begin{center}
\includegraphics{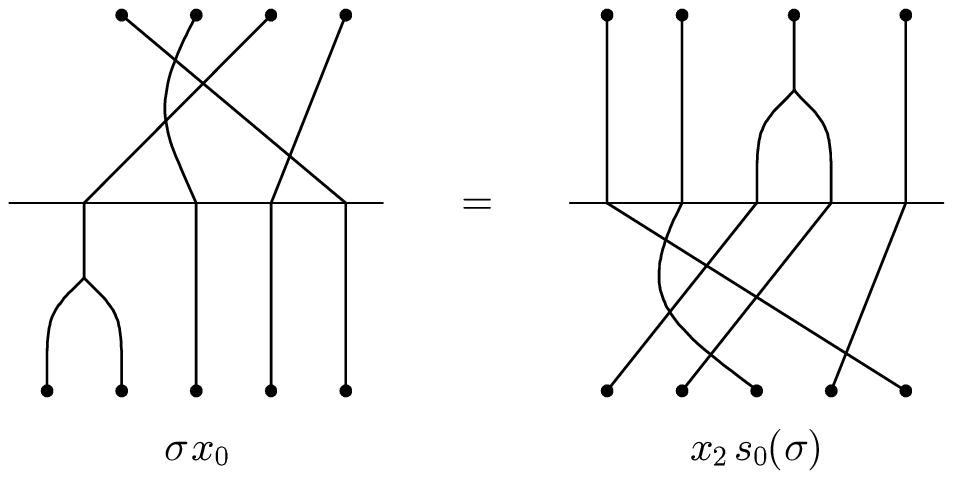}
\end{center}
It is not hard to show that $\mathcal{P}[\Sigma_w]$ is cancellative and has common right multiples.  Its groupoid
of right fractions $\mathcal{V}$ has fundamental group $V$.

To derive a presentation for $V$, recall that the symmetric group $\Sigma_w$ is generated by the adjacent transpositions
$t_1,\ldots,t_{w-1}$, with relations:
\begin{equation*}
t_n^2 =1\qquad\text{and}\qquad t_n t_{n+1} t_n = t_{n+1} t_n t_{n+1}
\end{equation*}
It is easy to check that:
\begin{equation*}
s_i(t_n)=
\begin{cases}
t_{n+1} & i < n-1 \\
t_n t_{n+1} & i = n-1 \\
t_{n+1} t_n & i = n \\
t_n & i > n
\end{cases}
\end{equation*}

\begin{theorem}The group $V$ is generated by the elements $\{x_0,x_1,x_2,\ldots\}$, \linebreak
$\{t_1,t_2,t_3,\ldots\}$, and $\{u_1,u_2,u_3,\ldots\}$, with relations:
\begin{equation*}
\begin{array}{l}
x_n x_k = x_k x_{n+1}\quad\text{for $k<n$} \\
t_n^2 = u_n^2 = 1 \\
t_n t_{n+1} t_n = t_{n+1} t_n t_{n+1},\quad t_n u_{n+1} t_n = u_{n+1} t_n u_{n+1}
\end{array}
\end{equation*}
and:
\begin{equation*}
\begin{array}{l}
t_n x_k = x_k t_{n+1},\quad u_n x_k = x_k u_{n+1}\quad\text{for $k<n-1$} \\
t_n x_{n-1} = x_n t_n t_{n+1},\quad u_n x_{n-1} = t_n u_{n+1} \\
t_n x_n = x_{n-1} t_{n+1} t_n,\quad u_n = x_{n-1} u_{n+1} t_n \\
t_k x_n = x_n t_k\quad\text{for $k<n$}
\end{array}
\end{equation*}
\end{theorem}
\begin{proof}
The groupoid $\mathcal{V}$ is generated by elements $x_n^{(w)}\colon w \rightarrow w+1$
($n<w$) and $t_n^{(w)} \colon w \rightarrow w$ ($1\leq n < w$) with relations:
\begin{equation*}
\begin{array}{ll}
x_n^{(w)} x_k^{(w+1)} = x_k^{(w)} x_{n+1}^{(w+1)} & \text{($k<n<w$)} \\
\bigl(t_n^{(w)}\bigr)^2 = 1 & \text{($1\leq n < w$)} \\
t_n^{(w)} t_{n+1}^{(w)} t_n^{(w)} = t_{n+1}^{(w)} t_n^{(w)} t_{n+1}^{(w)} & \text{($n+1<w$)} \\
t_n^{(w)} x_k^{(w)} = x_k^{(w)} t_{n+1}^{(w+1)} & \text{($k+1<n<w$)} \\
t_n^{(w)} x_{n-1}^{(w)} = x_n^{(w)} t_n^{(w+1)} t_{n+1}^{(w+1)} & \text{($n<w$)} \\
t_n^{(w)} x_n^{(w)} = x_{n-1}^{(w)} t_{n+1}^{(w+1)} t_n^{(w+1)} & \text{($n<w$)} \\
t_k^{(w)} x_n^{(w)} = x_n^{(w)} t_k^{(w+1)} & \text{($k<n<w$)} \\
\end{array}
\end{equation*}
Again, the first family of relations implies that $V$ contains a copy of $F$.  Substituting
$n=w-1$ into the last family of relations gives:
\begin{equation*}
t_k^{(k+2)}=t_k^{(k+3)}=t_k^{(k+4)} = \cdots
\end{equation*}
Let $t_k$ denote this element, and let $u_k$ denote the element $t_k^{(k+1)}$.  Then the relations
for $\mathcal{V}$ reduce to the given relations for $V=\pi_1(\mathcal{V},1)$.
\end{proof}

There is a braided version of $V$ called $BV$, introduced by Matt Brin \cite{Brin}.  It is the
fundamental group of the groupoid $\mathcal{BV}$ of right fractions of the category $\mathcal{P}[B_w]$
obtained by attaching a copy of the braid group $B_w$ to each object of $\mathcal{P}$:
\begin{center}
\includegraphics{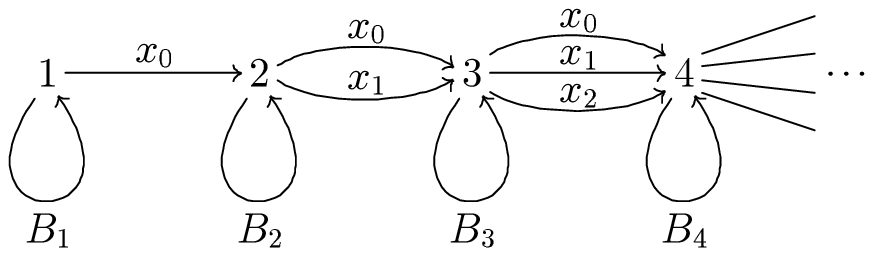}
\end{center}
This category satisfies the relations:
\begin{equation*}
bx_n = x_{b^{-1}(n)} s_n(b)\quad\text{($b\in B_w$ and $n<w$)}
\end{equation*}
where $b^{-1}(n)$ indicates the action of the braid $b^{-1}$ on $n$ via the projection
$B_w\twoheadrightarrow\Sigma_w$, and $s_n\colon B_w \rightarrow B_{w+1}$ is the function
that doubles the $n$'th strand of a braid:
\begin{center}
\includegraphics{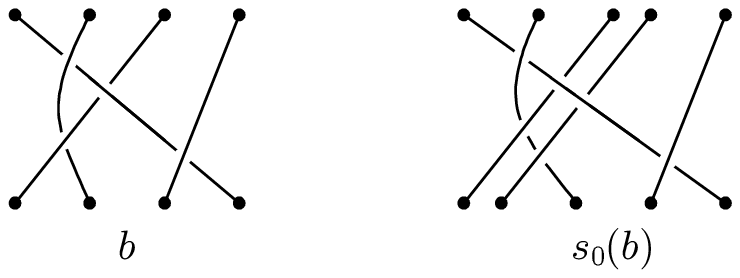}
\end{center}
Any element of the groupoid $\mathcal{BV}$ can be represented by a \emph{braided strand
diagram}:
\begin{center}
\includegraphics{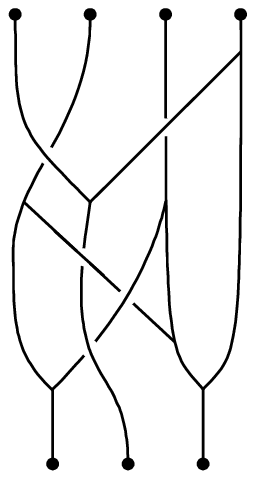}
\end{center}

\begin{theorem}The group $BV$ is generated by the elements $\{x_0,x_1,x_2,\ldots\}$,
$\{t_1,t_2,t_3,\ldots\}$, and $\{u_1,u_2,u_3,\ldots\}$, with relations:
\begin{equation*}
\begin{array}{l}
x_n x_k = x_k x_{n+1}\quad\text{for $k<n$} \\
t_n t_{n+1} t_n = t_{n+1} t_n t_{n+1},\quad t_n u_{n+1} t_n = u_{n+1} t_n u_{n+1}
\end{array}
\end{equation*}
and:
\begin{equation*}
\begin{array}{l}
t_n x_k = x_k t_{n+1},\quad u_n x_k = x_k u_{n+1}\quad\text{for $k<n-1$} \\
t_n x_{n-1} = x_n t_n t_{n+1},\quad u_n x_{n-1} = t_n u_{n+1} \\
t_n x_n = x_{n-1} t_{n+1} t_n,\quad u_n = x_{n-1} u_{n+1} t_n \\
t_k x_n = x_n t_k\quad\text{for $k<n$}
\end{array}\tag*{\qedsymbol}
\end{equation*}
\end{theorem}

\end{document}